\documentclass[12pt,reqno]{amsart}
\usepackage{amsmath}
\usepackage{amssymb}
\usepackage{verbatim}
\usepackage{mathrsfs}
\usepackage{graphicx}
\usepackage{caption}
\usepackage{subcaption}
\usepackage{epstopdf}
\usepackage{makecell}
\usepackage{tabularx}
\usepackage{longtable}
\usepackage{tikz}
\usepackage{pgfplots}
\usepackage{xcolor}
\usepackage{pgfplots}
\usepackage{pgfplotstable}
\usepackage[ruled,vlined]{algorithm2e}
\usepackage[noend]{algpseudocode}
\usepackage[toc]{appendix}
\usepackage[capitalise]{cleveref}
\usepackage[noadjust]{cite}
\usepackage{color}

\usetikzlibrary{arrows}

\usepackage[top=0.6in,bottom=0.5in,left=0.7in,right=0.7in]{geometry}
\newtheorem{lemma}{Lemma}

\newtheorem{theorem}[lemma]{Theorem}

\newtheorem{example}{Example}

\numberwithin{equation}{section}
\numberwithin{lemma}{section}
\newcommand{\N}{\mathbb{N}}    %natural numbers
\newcommand{\NN}{\mathbb{N}_0} %Natural Nonnegative numbers
\newcommand{\R}{\mathbb{R}}    %real number field

\newcommand{\nv}{\vec{n}}
\newcommand{\bo}{\mathcal{O}}

\newcommand{\Z}{\mathbb{Z}}    %real number field

\newcommand{\be}{ \begin{equation} }
\newcommand{\ee}{ \end{equation} }
\newcommand{\odd}{\operatorname{odd}}
\newcommand{\ind}{\Lambda}
\begin{document}
		
\title[Compact  Finite Difference Discretizations for  Elliptic Cross-Interface Problems]{Compact  9-Point  Finite Difference Methods with High Accuracy Order and/or  M-Matrix Property
	for  Elliptic Cross-Interface Problems}

\author{Qiwei Feng, Bin Han and Peter Minev}

	\thanks{Research supported in part by Natural Sciences and Engineering Research Council (NSERC) of Canada under grants RGPIN-2019-04276 (Bin Han), RGPIN-2017-04152 (Peter Minev),  Westgrid (www.westgrid.ca), and Compute Canada Calcul Canada (www.computecanada.ca)}

\address{Department of Mathematical and Statistical Sciences,
University of Alberta, Edmonton,\quad Alberta, Canada T6G 2G1.
\quad {\tt qfeng@ualberta.ca }\qquad {\tt bhan@ualberta.ca}
\quad{\tt minev@ualberta.ca}
}

\begin{abstract}
In this paper we develop finite difference schemes for elliptic problems with piecewise continuous coefficients that have (possibly huge) jumps across fixed internal interfaces.
In contrast with such problems involving one smooth non-intersecting interface, that have been extensively studied,
there are very few papers addressing elliptic interface problems with intersecting interfaces of coefficient jumps. It is well known that if the values of the permeability in the four subregions around a point of intersection of two such internal interfaces are all different,
the solution has a point singularity that significantly affects the accuracy of the approximation in the vicinity of the intersection point.
In the present paper we propose a fourth-order 9-point finite difference scheme on uniform Cartesian meshes for an elliptic problem whose coefficient is piecewise constant in four rectangular subdomains of the overall two-dimensional rectangular domain.
Moreover, for the special case when the  intersecting point of the two lines of coefficient jumps is a grid point, such a compact scheme, involving relatively simple formulas for computation of the stencil coefficients, can even reach sixth order of accuracy. Furthermore, we show that the resulting linear system for the special case has an M-matrix, and prove the theoretical sixth order convergence rate using the discrete maximum principle.
Our numerical experiments demonstrate the sixth (for the special case) and at least fourth (for the general case)  accuracy orders of the proposed schemes.
In the general case, we derive a compact  third-order finite difference scheme, also yielding a linear system with an  M-matrix.  In addition, using the discrete maximum principle, we prove the third order convergence rate of the scheme for the general elliptic cross-interface problem.
\end{abstract}


\maketitle

\maketitle

\pagenumbering{arabic}

\section{Introduction and problem formulation}\label{introdu:1:intersect}

	Interface problems arise in many applications such as modeling of underground waste
disposal, oil reservoirs, composite materials, and many others.
Some approaches to the solution of the elliptic interface problem with a smooth non-intersecting interface of coefficient jumps were provided by the immersed interface methods (IIM, see \cite{CFL19,LeLi94,EwingLLL99,GongLiLi08,HeLL2011,Li98,WieBube00,PanHeLi21,XiaolinZhong07,DFL20}
 and references therein) and matched interface and boundary methods (MIB, see \cite{YuZhouWei07,ZW06,FendZhao20pp109677,YuWei073D,ZZFW06}). Both methods belong to the class of the finite difference methods (FDM).
Elliptic interface problems with intersecting interfaces appear in many applications, but perhaps the most notorious example is the modeling of geological porous media flows  (see e.g. \cite{Minev2018,AliMankad2020,GuiraAusas2019,ArbogastXiao2013,ArbogastTao2013,Jaramillo2022,Kippe2008,Tahmasebi2018,Butler2020,TArbHXiao2015,Rasaei2008}).  A classical problem of this type is formulated by the Society of Petroleum Engineers, the so-called SPE10 problem (see https://www.spe.org/web/csp/datasets/set02.htm).  Here we consider a 2D simplification of this problem
that involves  interface intersections of vertical straight lines and horizontal straight lines, so that the permeability coefficient in the four subregions in the vicinity of an intersection point has different values.
Even for this relatively simple cross-interface problem, the only compact  9-point finite difference method in the literature, that we are aware of, is the scheme in \cite{AnVu2007}, that is third-order consistent, and uses special non-uniform meshes.  To our knowledge, the convergence rate of this scheme has never been proven.
We should also note here that the difficulty of the problem is usually  exacerbated if the jumps of the permeability coefficient across the different interfaces are very large (of several orders of magnitude). Details of the physical background of the elliptic interface problem can be found in \cite{Vazqu07}.

In practice, usually the permeability variation occurs on scales that are very small as compared to the size of the medium, and therefore the solution of the interface problem is highly oscillatory. This causes the appearance of the so-called pollution effect in the error of
its numerical approximation.
In order to obtain a reasonable low-order numerical solution to such problems one needs to employ a very fine, possibly nonuniform grid,  that captures the small scale features.  Therefore, the development of higher-order compact approximations can help to reduce the computational costs.
Compared to the finite element or finite volume methods, the FDM does not require the integration of highly-oscillatory or discontinuous functions. Furthermore, a compact  9-point scheme (in 2D) yields a sparse linear system that can be solved very efficiently.

Finally we should mention, that even in the relatively simple case of elliptic interface problems involving a smooth non-intersecting interface of coefficient jumps, the theoretical proof of convergence of the various proposed finite difference schemes is usually missing.
The only exceptions are presented in \cite{LiIto2001} using discrete maximum principle for a second order scheme, and \cite{FengHanMinev22} using numerically verified discrete maximum principle for a fourth order scheme.
The compact  9-point  schemes considered in the present paper possess the M-matrix property, that guarantees the discrete maximum principle for the numerical solution.  In turn, this property greatly facilitates the proof of their convergence rate.

%%%%%%%%%%%%%%%%%%%%%%%%%%%%%%%%%%%%%%%%%%%%%%%%%%%%%%%%%%%%%%%%%%%%%%%%%%%%%%%%%%%%%%%%%%%%%%%%%%%%%

In this paper we develop numerical approximations to the following elliptic cross-interface problem: {\em Given the domain $\Omega:=(l_1,l_2)\times(l_3,l_4)$ with $l_1,l_2, l_3,l_4 \in \R$, then consider:
\begin{equation} \label{intersect:3}
	\left\{ \begin{array}{llll}
		-\nabla \cdot \Big(a \nabla  u\Big)&=&f %\hskip 8mm
&\mbox{in } \Omega \setminus {\Gamma},\\
		\left[u\right]&=&\phi_p %\hskip 8mm
&\mbox{on }\Gamma_p \mbox{ for } p=1,2,3,4,\\	
		\left[a \nabla  u \cdot \vec{n}\right] &=&\psi_p %\hskip 8mm
&\mbox{on }\Gamma_p \mbox{ for } p=1,2,3,4,\\
		u&=&g %\hskip 7mm
&\mbox{on }	\partial\Omega,
	\end{array}
	\right.
\end{equation}
where the cross-interface $\Gamma$ is given by $\Gamma:=\Gamma_1 \cup \Gamma_2 \cup \Gamma_3 \cup \Gamma_4\cup\{(\xi,\zeta)\}$ with
%
%\begin{align*}
%&\Gamma_1=\{y: x=\xi, \eta < y < l_4\},\qquad
%\Gamma_2=\{y: x=\xi, l_3< y < \eta\}, \\
%&\Gamma_3=\{x: y=\eta, \xi< x < l_2\},\qquad
%\Gamma_4=\{x: y=\eta, l_1< x < \xi\}
%\end{align*}
%
\[
\Gamma_1:=\{\xi\}\times (\zeta,l_4), \qquad
\Gamma_2:=\{\xi\}\times (l_3,\zeta),\qquad
\Gamma_3:=(\xi,l_2)\times \{\zeta\}, \qquad
\Gamma_4:=(l_1,\xi)\times\{\zeta\}, \qquad
(\xi,\zeta) \in \Omega.
\]
}
As usual, that the square brackets here denote the jump of the corresponding function, i.e. for $(\xi,y)\in \Gamma_p$ with $p=1,2$ (on the vertical line of the cross-interface $\Gamma$),
\[
[u](\xi,y):=\lim_{x \to \xi^+}u(x,y)- \lim_{ x\to \xi^-}u(x,y),
\quad
[a \nabla  u \cdot \vec{n}](\xi,y):=\lim_{x \to \xi^+} a(x,y) \frac{\partial u}{\partial x}(x,y) - \lim_{x\to \xi^-} a(x,y)\frac{\partial u}{\partial x}(x,y);
\]
while for $(x,\zeta)\in \Gamma_p$ with $p=3,4$ (i.e., on the horizontal line of the cross-interface $\Gamma$),
\[
[u](x,\zeta):=\lim_{y\to \zeta^+ }u(x,y)- \lim_{ y\to \zeta^-}u(x,y),  \quad
[a \nabla  u \cdot \vec{n}](x,\zeta):=\lim_{y \to \zeta^+} a(x,y) \frac{\partial u}{\partial y}(x,y) - \lim_{y\to \zeta^-} a(x,y)\frac{\partial u}{\partial y}(x,y).
\]

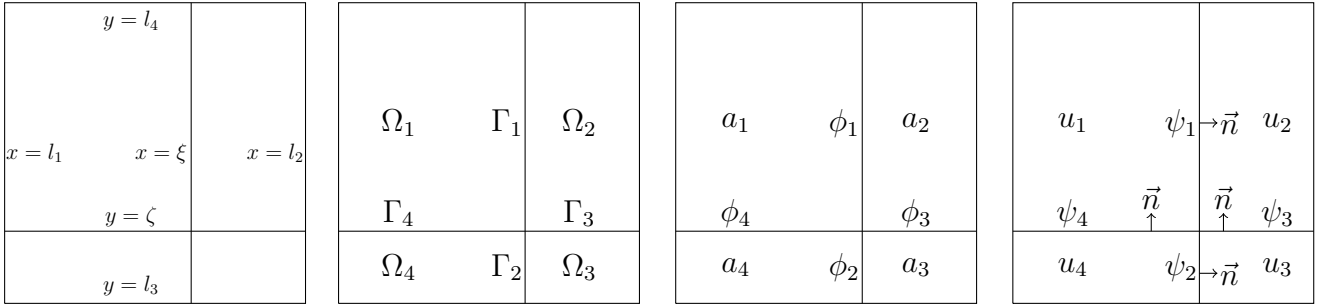
\begin{figure}[htbp]
	\centering
	\begin{subfigure}[b]{0.2\textwidth}
	\hspace{-0.9cm}
	\begin{tikzpicture}[scale = 0.8]
		\draw	(-2.5, -2.5) -- (-2.5, 2.5) -- (2.5, 2.5) -- (2.5, -2.5) --(-2.5,-2.5);	
		\draw   (0.6, -2.5) -- (0.6, 2.5);	
		\draw   (-2.5, -1.3) -- (2.5, -1.3);	
		 %%%%%%%%%%%%%%%%%%%%%%%%%%%%%%%%%%%%%%%%%%%%%%%%%
		\node (A) [scale=0.7] at (0.1,0) {$x=\xi$};
		\node (A) [scale=0.7] at (-0.4,-1.05) {$y=\zeta$};
		\node (A) [scale=0.7] at (-2,0) {$x=l_1$};
		\node (A) [scale=0.7] at (2,0) {$x=l_2$};
		\node (A) [scale=0.7] at (-0.4,2.2) {$y=l_4$};
		\node (A) [scale=0.7] at (-0.4,-2.2) {$y=l_3$};
	\end{tikzpicture}
\end{subfigure}		
	\begin{subfigure}[b]{0.2\textwidth}
		\hspace{-0.3cm}
		\begin{tikzpicture}[scale = 0.8]
			\draw	(-2.5, -2.5) -- (-2.5, 2.5) -- (2.5, 2.5) -- (2.5, -2.5) --(-2.5,-2.5);	
			\draw   (0.6, -2.5) -- (0.6, 2.5);	
			\draw   (-2.5, -1.3) -- (2.5, -1.3);	 
			\node (A) at (0.3,-1.9) {$\Gamma_2$};
			\node (A) at (0.3,0.5) {$\Gamma_1$};
			\node (A) at (1.5,-1) {$\Gamma_3$};
			\node (A) at (-1.5,-1) {$\Gamma_4$};
			 %%%%%%%%%%%%%%%%%%%%%%%%%%%%%%%%%%%%%%%%%%%%%%%%%
			\node (A) at (-1.5,0.5) {$\Omega_1$};
			\node (A) at (1.5,0.5) {$\Omega_2$};
			\node (A) at (1.5,-1.9) {$\Omega_3$};
			\node (A) at (-1.5,-1.9) {$\Omega_4$};
		\end{tikzpicture}
	\end{subfigure}
	\hspace{0.6cm}
	\begin{subfigure}[b]{0.2\textwidth}
		\begin{tikzpicture}[scale = 0.8]
	\draw	(-2.5, -2.5) -- (-2.5, 2.5) -- (2.5, 2.5) -- (2.5, -2.5) --(-2.5,-2.5);	
	\draw   (0.6, -2.5) -- (0.6, 2.5);	
	\draw   (-2.5, -1.3) -- (2.5, -1.3);	
\node (A) at (0.3,-1.9) {$\phi_2$};
\node (A) at (0.3,0.5) {$\phi_1$};
\node (A) at (1.5,-1) {$\phi_3$};
\node (A) at (-1.5,-1) {$\phi_4$};
%%%%%%%%%%%%%%%%%%%%%%%%%%%%%%%%%%%%%%%%%%%%%%%%%
\node (A) at (-1.5,0.5) {$a_1$};
\node (A) at (1.5,0.5) {$a_2$};
\node (A) at (1.5,-1.9) {$a_3$};
\node (A) at (-1.5,-1.9) {$a_4$};
\end{tikzpicture}
	\end{subfigure}
	\hspace{0.6cm}
	%%%%%%%%%%%%%%%%%%%%%%%
	\begin{subfigure}[b]{0.2\textwidth}
			\begin{tikzpicture}[scale = 0.8]
			\draw	(-2.5, -2.5) -- (-2.5, 2.5) -- (2.5, 2.5) -- (2.5, -2.5) --(-2.5,-2.5);	
			\draw   (0.6, -2.5) -- (0.6, 2.5);	
			\draw   (-2.5, -1.3) -- (2.5, -1.3);	 
			\node (A) at (0.3,-1.9) {$\psi_2$};
			\node (A) at (0.3,0.5) {$\psi_1$};
			\node (A) at (1.9,-1) {$\psi_3$};
			\node (A) at (-1.5,-1) {$\psi_4$};
			 %%%%%%%%%%%%%%%%%%%%%%%%%%%%%%%%%%%%%%%%%%%%%%%%%
			\node (A) at (-1.5,0.5) {$u_1$};
			\node (A) at (1.9,0.5) {$u_2$};
			\node (A) at (1.9,-1.9) {$u_3$};
			\node (A) at (-1.5,-1.9) {$u_4$};
			 %%%%%%%%%%%%%%%%%%%%%%%%%%%%%%%%%%%%%%%%%%%%%%%%%%%
			\node (A) at (1.1,0.5) {$\nv$};
			\draw[->] (0.6,0.5)--(0.9,0.5);
			\node (A) at (1.1,-2) {$\nv$};
			\draw[->] (0.6,-2)--(0.9,-2);
			\node (A) at (1,-0.75) {$\nv$};
			\draw[->] (1,-1.3)--(1,-1);
			\node (A) at (-0.2,-0.75) {$\nv$};
			\draw[->] (-0.2,-1.3)--(-0.2,-1);
		\end{tikzpicture}
	\end{subfigure}
	\caption{An illustration for the model elliptic cross-interface problem in \eqref{intersect:3}. }
	\label{fig:intersect:3}
\end{figure}

Note that  the interface curve $\Gamma$ divides the domain $\Omega$ into $4$ subdomains:
\[
\Omega_1:=(l_1,\xi)\times(\zeta,l_4), \quad \Omega_2:=(\xi,l_2)\times(\zeta,l_4), \quad  \Omega_3:=(\xi,l_2)\times(l_3,\zeta), \quad \Omega_4:=(l_1,\xi)\times(l_3,\zeta).
\]
See \cref{fig:intersect:3} for an illustration, where
$a_{p}:=a \chi_{\Omega_{p}}$, $f_{p}:=f \chi_{\Omega_{p}}$, and $u_{p}:=u \chi_{\Omega_{p}}$ for $p=1,2,3,4$.

To derive a compact  9-point scheme that approximates
the cross-interface problem
\eqref{intersect:3}, we assume that:
\begin{itemize}
	
	\item[(A1)] $a_{p}$ is a positive constant in $\Omega_p$.
	
	\item[(A2)] The restriction of the solution $u_{p}$ and of the source term $f_{p}$ has uniformly continuous partial derivatives of (total) orders up to seven and five, respectively, in each $\Omega_p$ for $p=1,2,3,4$.
%The source term $f$ can be continuous or discontinuous across the interfaces $\Gamma_i$ for $i=1,2,3,4$;
	
	\item[(A3)] The essentially one-dimensional functions $\phi_p$ and $\psi_p$ in \eqref{intersect:3} on the interface $\Gamma_p$ have uniformly continuous derivatives of orders up to seven and six respectively for $p=1,2,3,4$.
\end{itemize}

%%%%
%
%
%
%
%
%
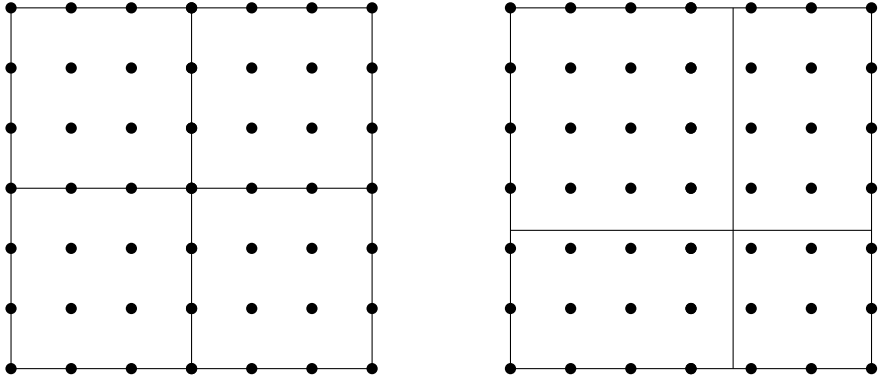
\begin{figure}[htbp]
	\centering	
	\hspace{-1cm}
	\begin{subfigure}[b]{0.3\textwidth}
		\begin{tikzpicture}[scale = 0.8]
			\draw	(-3, -3) -- (-3, 3) -- (3, 3) -- (3, -3) --(-3,-3);	
			\draw   (0, -3) -- (0, 3);	
			\draw   (-3, 0) -- (3, 0);	
		 %%%%%%%%%%%%%%%%%%%%%%%%%%%%%%%%%%%%%%%%%%%%%%%%%%%
	\node at (0,0)[circle,fill,inner sep=1.5pt,color=black]{};
	\node at (1,0)[circle,fill,inner sep=1.5pt,color=black]{};
	\node at (2,0)[circle,fill,inner sep=1.5pt,color=black]{};
	\node at (3,0)[circle,fill,inner sep=1.5pt,color=black]{};
	\node at (0,0)[circle,fill,inner sep=1.5pt,color=black]{};
	\node at (-1,0)[circle,fill,inner sep=1.5pt,color=black]{};
	\node at (-2,0)[circle,fill,inner sep=1.5pt,color=black]{};
	\node at (-3,0)[circle,fill,inner sep=1.5pt,color=black]{};
	 %%%%%%%%%%%%%%%%%%%%%%%%%%%%%%%%%%%%%%%%%%
	 %%%%%%%%%%%%%%%%%%%%%%%%%%%%%%%%%%%%%%%%%%%%%%%%%%%
	\node at (0,1)[circle,fill,inner sep=1.5pt,color=black]{};
	\node at (1,1)[circle,fill,inner sep=1.5pt,color=black]{};
	\node at (2,1)[circle,fill,inner sep=1.5pt,color=black]{};
	\node at (3,1)[circle,fill,inner sep=1.5pt,color=black]{};
	\node at (0,1)[circle,fill,inner sep=1.5pt,color=black]{};
	\node at (-1,1)[circle,fill,inner sep=1.5pt,color=black]{};
	\node at (-2,1)[circle,fill,inner sep=1.5pt,color=black]{};
	\node at (-3,1)[circle,fill,inner sep=1.5pt,color=black]{};
	 %%%%%%%%%%%%%%%%%%%%%%%%%%%%%%%%%%%%%%%%%%%%%%%%%%%
	\node at (0,2)[circle,fill,inner sep=1.5pt,color=black]{};
	\node at (1,2)[circle,fill,inner sep=1.5pt,color=black]{};
	\node at (2,2)[circle,fill,inner sep=1.5pt,color=black]{};
	\node at (3,2)[circle,fill,inner sep=1.5pt,color=black]{};
	\node at (0,2)[circle,fill,inner sep=1.5pt,color=black]{};
	\node at (-1,2)[circle,fill,inner sep=1.5pt,color=black]{};
	\node at (-2,2)[circle,fill,inner sep=1.5pt,color=black]{};
	\node at (-3,2)[circle,fill,inner sep=1.5pt,color=black]{};
	 %%%%%%%%%%%%%%%%%%%%%%%%%%%%%%%%%%%%%%%%%%%%%%%%%%%
	\node at (0,3)[circle,fill,inner sep=1.5pt,color=black]{};
	\node at (1,3)[circle,fill,inner sep=1.5pt,color=black]{};
	\node at (2,3)[circle,fill,inner sep=1.5pt,color=black]{};
	\node at (3,3)[circle,fill,inner sep=1.5pt,color=black]{};
	\node at (0,3)[circle,fill,inner sep=1.5pt,color=black]{};
	\node at (-1,3)[circle,fill,inner sep=1.5pt,color=black]{};
	\node at (-2,3)[circle,fill,inner sep=1.5pt,color=black]{};
	\node at (-3,3)[circle,fill,inner sep=1.5pt,color=black]{};
	 %%%%%%%%%%%%%%%%%%%%%%%%%%%%%%%%%%%%%%%%%%
	 %%%%%%%%%%%%%%%%%%%%%%%%%%%%%%%%%%%%%%%%%%%%%%%%%%%
	\node at (0,-1)[circle,fill,inner sep=1.5pt,color=black]{};
	\node at (1,-1)[circle,fill,inner sep=1.5pt,color=black]{};
	\node at (2,-1)[circle,fill,inner sep=1.5pt,color=black]{};
	\node at (3,-1)[circle,fill,inner sep=1.5pt,color=black]{};
	\node at (0,-1)[circle,fill,inner sep=1.5pt,color=black]{};
	\node at (-1,-1)[circle,fill,inner sep=1.5pt,color=black]{};
	\node at (-2,-1)[circle,fill,inner sep=1.5pt,color=black]{};
	\node at (-3,-1)[circle,fill,inner sep=1.5pt,color=black]{};
	 %%%%%%%%%%%%%%%%%%%%%%%%%%%%%%%%%%%%%%%%%%%%%%%%%%%
	\node at (0,-2)[circle,fill,inner sep=1.5pt,color=black]{};
	\node at (1,-2)[circle,fill,inner sep=1.5pt,color=black]{};
	\node at (2,-2)[circle,fill,inner sep=1.5pt,color=black]{};
	\node at (3,-2)[circle,fill,inner sep=1.5pt,color=black]{};
	\node at (0,-2)[circle,fill,inner sep=1.5pt,color=black]{};
	\node at (-1,-2)[circle,fill,inner sep=1.5pt,color=black]{};
	\node at (-2,-2)[circle,fill,inner sep=1.5pt,color=black]{};
	\node at (-3,-2)[circle,fill,inner sep=1.5pt,color=black]{};
	 %%%%%%%%%%%%%%%%%%%%%%%%%%%%%%%%%%%%%%%%%%%%%%%%%%%
	\node at (0,-3)[circle,fill,inner sep=1.5pt,color=black]{};
	\node at (1,-3)[circle,fill,inner sep=1.5pt,color=black]{};
	\node at (2,-3)[circle,fill,inner sep=1.5pt,color=black]{};
	\node at (3,-3)[circle,fill,inner sep=1.5pt,color=black]{};
	\node at (0,-3)[circle,fill,inner sep=1.5pt,color=black]{};
	\node at (-1,-3)[circle,fill,inner sep=1.5pt,color=black]{};
	\node at (-2,-3)[circle,fill,inner sep=1.5pt,color=black]{};
	\node at (-3,-3)[circle,fill,inner sep=1.5pt,color=black]{};
		\end{tikzpicture}
	\end{subfigure}
%%%%%%%%%%%%%%%%%%%%%%%%%%%%%%%%%%%%%%%%%%%%%%%%%%%%%%%%%%%%%%%%%%%%%%%%%%%%%%%%%%%%
	\begin{subfigure}[b]{0.3\textwidth}
			\hspace{1cm}
	\begin{tikzpicture}[scale = 0.8]
		\draw	(-3, -3) -- (-3, 3) -- (3, 3) -- (3, -3) --(-3,-3);	
		\draw   (0.7, -3) -- (0.7, 3);	
		\draw   (-3,-0.7 ) -- (3, -0.7);	
		 %%%%%%%%%%%%%%%%%%%%%%%%%%%%%%%%%%%%%%%%%%%%%%%%%%%
		\node at (0,0)[circle,fill,inner sep=1.5pt,color=black]{};
		\node at (1,0)[circle,fill,inner sep=1.5pt,color=black]{};
		\node at (2,0)[circle,fill,inner sep=1.5pt,color=black]{};
		\node at (3,0)[circle,fill,inner sep=1.5pt,color=black]{};
		\node at (0,0)[circle,fill,inner sep=1.5pt,color=black]{};
		\node at (-1,0)[circle,fill,inner sep=1.5pt,color=black]{};
		\node at (-2,0)[circle,fill,inner sep=1.5pt,color=black]{};
		\node at (-3,0)[circle,fill,inner sep=1.5pt,color=black]{};
		 %%%%%%%%%%%%%%%%%%%%%%%%%%%%%%%%%%%%%%%%%%
		 %%%%%%%%%%%%%%%%%%%%%%%%%%%%%%%%%%%%%%%%%%%%%%%%%%%
		\node at (0,1)[circle,fill,inner sep=1.5pt,color=black]{};
		\node at (1,1)[circle,fill,inner sep=1.5pt,color=black]{};
		\node at (2,1)[circle,fill,inner sep=1.5pt,color=black]{};
		\node at (3,1)[circle,fill,inner sep=1.5pt,color=black]{};
		\node at (0,1)[circle,fill,inner sep=1.5pt,color=black]{};
		\node at (-1,1)[circle,fill,inner sep=1.5pt,color=black]{};
		\node at (-2,1)[circle,fill,inner sep=1.5pt,color=black]{};
		\node at (-3,1)[circle,fill,inner sep=1.5pt,color=black]{};
		 %%%%%%%%%%%%%%%%%%%%%%%%%%%%%%%%%%%%%%%%%%%%%%%%%%%
		\node at (0,2)[circle,fill,inner sep=1.5pt,color=black]{};
		\node at (1,2)[circle,fill,inner sep=1.5pt,color=black]{};
		\node at (2,2)[circle,fill,inner sep=1.5pt,color=black]{};
		\node at (3,2)[circle,fill,inner sep=1.5pt,color=black]{};
		\node at (0,2)[circle,fill,inner sep=1.5pt,color=black]{};
		\node at (-1,2)[circle,fill,inner sep=1.5pt,color=black]{};
		\node at (-2,2)[circle,fill,inner sep=1.5pt,color=black]{};
		\node at (-3,2)[circle,fill,inner sep=1.5pt,color=black]{};
		 %%%%%%%%%%%%%%%%%%%%%%%%%%%%%%%%%%%%%%%%%%%%%%%%%%%
		\node at (0,3)[circle,fill,inner sep=1.5pt,color=black]{};
		\node at (1,3)[circle,fill,inner sep=1.5pt,color=black]{};
		\node at (2,3)[circle,fill,inner sep=1.5pt,color=black]{};
		\node at (3,3)[circle,fill,inner sep=1.5pt,color=black]{};
		\node at (0,3)[circle,fill,inner sep=1.5pt,color=black]{};
		\node at (-1,3)[circle,fill,inner sep=1.5pt,color=black]{};
		\node at (-2,3)[circle,fill,inner sep=1.5pt,color=black]{};
		\node at (-3,3)[circle,fill,inner sep=1.5pt,color=black]{};
		 %%%%%%%%%%%%%%%%%%%%%%%%%%%%%%%%%%%%%%%%%%
		 %%%%%%%%%%%%%%%%%%%%%%%%%%%%%%%%%%%%%%%%%%%%%%%%%%%
		\node at (0,-1)[circle,fill,inner sep=1.5pt,color=black]{};
		\node at (1,-1)[circle,fill,inner sep=1.5pt,color=black]{};
		\node at (2,-1)[circle,fill,inner sep=1.5pt,color=black]{};
		\node at (3,-1)[circle,fill,inner sep=1.5pt,color=black]{};
		\node at (0,-1)[circle,fill,inner sep=1.5pt,color=black]{};
		\node at (-1,-1)[circle,fill,inner sep=1.5pt,color=black]{};
		\node at (-2,-1)[circle,fill,inner sep=1.5pt,color=black]{};
		\node at (-3,-1)[circle,fill,inner sep=1.5pt,color=black]{};
		 %%%%%%%%%%%%%%%%%%%%%%%%%%%%%%%%%%%%%%%%%%%%%%%%%%%
		\node at (0,-2)[circle,fill,inner sep=1.5pt,color=black]{};
		\node at (1,-2)[circle,fill,inner sep=1.5pt,color=black]{};
		\node at (2,-2)[circle,fill,inner sep=1.5pt,color=black]{};
		\node at (3,-2)[circle,fill,inner sep=1.5pt,color=black]{};
		\node at (0,-2)[circle,fill,inner sep=1.5pt,color=black]{};
		\node at (-1,-2)[circle,fill,inner sep=1.5pt,color=black]{};
		\node at (-2,-2)[circle,fill,inner sep=1.5pt,color=black]{};
		\node at (-3,-2)[circle,fill,inner sep=1.5pt,color=black]{};
		 %%%%%%%%%%%%%%%%%%%%%%%%%%%%%%%%%%%%%%%%%%%%%%%%%%%
		\node at (0,-3)[circle,fill,inner sep=1.5pt,color=black]{};
		\node at (1,-3)[circle,fill,inner sep=1.5pt,color=black]{};
		\node at (2,-3)[circle,fill,inner sep=1.5pt,color=black]{};
		\node at (3,-3)[circle,fill,inner sep=1.5pt,color=black]{};
		\node at (0,-3)[circle,fill,inner sep=1.5pt,color=black]{};
		\node at (-1,-3)[circle,fill,inner sep=1.5pt,color=black]{};
		\node at (-2,-3)[circle,fill,inner sep=1.5pt,color=black]{};
		\node at (-3,-3)[circle,fill,inner sep=1.5pt,color=black]{};
	\end{tikzpicture}
\end{subfigure}
	\caption{An illustration for uniform Cartesian grids of the model problem in \eqref{intersect:3}}
	\label{compact:intersect:model}
\end{figure}

The remainder of the paper is organized as follows.
In \cref{Intersect:subsec:regular}, we derive a compact  9-point scheme with sixth order of consistency for interior grid points in \cref{Intersect:thm:regular}.
For grid points near the interface, as illustrated by \cref{compact:intersect:model}, we have two cases:

\textbf{Case 1:} If the  point $(\xi,\zeta)$ of intersection of the interfaces is a grid point (see the left panel of \cref{compact:intersect:model}),
we derive in \cref{Intersect:subsec:irregular} a compact 9-point  scheme that has a seventh order of consistency at every grid point lying on the cross-interface.
The stencil coefficients are given in
\cref{Intersect:thm:cross1,Intersect:thm:side1}.
In \cref{convergence:section1} we prove that this scheme is sixth-order accurate, using the discrete maximum principle satisfied by it. The results of some numerical experiments, demonstrating the sixth-order convergence rate of the scheme, are presented in \cref{numerical:section1}.
%Once the interface intersecting point $(\xi,\zeta)$ is a grid point of a mesh for some mesh size $h$, it is easy to observe that $(\xi,\zeta)$ remains as a grid point for all finer meshes using the mesh sizes $2^{-j} h$ with $j\in \N$, because these meshes are nested.

\textbf{Case 2:} If  $(\xi,\zeta)$ is not a grid point  (see the right panel in \cref{compact:intersect:model}), we derive
in  \cref{Intersect:subsec:irregular} a compact  9-point  scheme with fourth order of consistency for every grid point neighboring the interface.
Next we show in \cref{convergence:section2}
that this scheme does not satisfy the M-matrix property.
Subsequently, we obtain  a compact scheme satisfying the M-matrix property, with a consistency order three at grid points neighboring the interface points except for the vicinity of the intersection point $(\xi, \zeta)$, and  order two at grid points neighboring $(\xi,\zeta)$.
Since the M-matrix property immediately implies that the scheme satisfies a discrete maximum principle, this allows us to prove that the overall convergence rate of the scheme is of order three.
In \cref{numerical:section2} we provide some numerical results that seem to suggest that the scheme given in \cref{Intersect:thm:regular,theorem:side3:w,Intersect:thm:cross1:w1w2}, that does not satisfy a discrete maximum principle, is  fifth-order accurate.

In \cref{Intersect:sec:Conclu},  we summarize the main contributions of this paper. Finally, in \cref{Intersect:sec:proofs}  %\cref{Intersect:sec:proofs},
we present the
proofs for the results stated in \cref{Intersect:sec:sixord,error:analysis:intersect:sec}.

%%%%%%%%%%%%%%%%%%%%%%%%%%%%%%%%%%
\section{High order compact  9-point  schemes using uniform Cartesian grids}
\label{Intersect:sec:sixord}

In this section, we present some compact finite difference schemes on uniform Cartesian grids for the elliptic cross-interface problem in \eqref{intersect:3}. To improve readability, the technical proofs of the results stated in this section are deferred to \cref{Intersect:sec:proofs}.

We start by introducing a uniform Cartesian mesh on the domain:
\[
\Omega:=(l_1,l_2)\times (l_3,l_4), \quad \mbox{with} \quad  l_4-l_3=N_0 (l_2-l_1) \quad \mbox{for some positive integer }  N_0,
\]
containing the grid points:
\[
x_i:=l_1+i h, \quad i=0,\ldots,N_1, \quad \text{and} \quad y_j:=l_3+j h, \quad j=0,\ldots,N_2, \quad h:=\tfrac{l_2-l_1}{N_1}=\tfrac{l_4-l_3}{N_2},
\]
where $N_1$ is a positive integer and $N_2:=N_0 N_1$.
We also define $(u_h)_{i,j}$ to be the value of the numerical approximation  $u_h$ of the exact solution $u$ of the elliptic cross-interface problem \eqref{intersect:3}, at the grid point $(x_i, y_j)$.
For stencil coefficients $\{C_{k,\ell}\}_{k,\ell=-1,0,1}$ with $C_{k,\ell}\in \R$ in the compact 9-point stencil centered at a grid point $(x_i,y_j)$, the discrete operator $\mathcal{L}_h$ acting on $u_h$ is defined to be:
\be \label{Luh}
\begin{split}
\mathcal{L}_h u_h:=\sum_{k=-1}^1\sum_{\ell=-1}^1  C_{k,\ell}(u_{h})_{i+k,j+\ell}.
\end{split}
\ee
Similarly, the action of the discrete operator $\mathcal{L}_h$  on the exact solution $u$ is given by:
\be \label{Lu}
\begin{split}
\mathcal{L}_h u :=	 \sum_{k=-1}^1\sum_{\ell=-1}^1  C_{k,\ell}u(x_i+kh,y_j+\ell h).
\end{split}
\ee

\subsection{Compact 9-point stencils at interior points}
\label{Intersect:subsec:regular}
The following compact  FDM with a consistency order six for \eqref{intersect:3} at interior points is well known in the literature (e.g., see \cite{SDKS13,FengHanMinev21}).
\begin{theorem}\label{Intersect:thm:regular}
	Consider  $(x_i,y_j)\in \Omega$ with all 9 points $(x_i\pm h,y_j\pm h)\in \overline{ \Omega_p}$ for some $p\in \{1,2,3,4\}$.  Assume that $u_{p}:=u \chi_{\Omega_{p}}$ and $f_{p}:=f \chi_{\Omega_{p}}$ have uniformly continuous partial derivatives of (total) orders up to seven and five, respectively, in $\Omega_p$.
Let the discrete operator $\mathcal{L}_h$ be defined in \eqref{Luh} with the stencil coefficients
\[
C_{0,0}=20,\quad C_{0,-1}=C_{0,1}=C_{-1,0}=C_{1,0}=-4,\quad
C_{1,1}=C_{1,-1}=C_{-1,1}=C_{-1,-1}=-1.
\]
Then the compact 9-point finite difference scheme
$h^{-2}\mathcal{L}_h u_h=\frac{1}{a(x_i,y_j)} F$ with
\[
F:=6f{(x_i, y_j)}+\tfrac{h^2}{2} \left[  \tfrac{\partial^{2} f}{ \partial^2 x}(x_i,y_j) +\tfrac{\partial^{2} f}{ \partial^2 y}(x_i,y_j)  \right]+\tfrac{h^4}{60}\left[ \tfrac{\partial^{4} f}{ \partial^4 x}(x_i,y_j) +\tfrac{\partial^{4} f}{ \partial^4 y}(x_i,y_j)   \right]+\tfrac{h^4}{15}\tfrac{\partial^{4} f}{ \partial^2 x \partial^2 y}(x_i,y_j),
\]
has a sixth order of consistency at the interior grid point $(x_i,y_j)$ for $-\nabla \cdot (a \nabla  u)=f$.
% i.e., $h^{-2}\mathcal{L}_h(u_h-u)=\bo(h^6)$ as $h\to 0$.
\end{theorem}

\subsection{Compact 9-point stencils at grid points on the interface}
\label{Intersect:subsec:irregular}

In this subsection, we now discuss how to find a compact  FDM of a consistency order seven at grid points $(x_i,y_j)$ lying on the cross-interface (i.e., Case 1 in \cref{introdu:1:intersect}, see the left panel of \cref{compact:intersect:model}). Note that all interfaces $\Gamma_1,\ldots,\Gamma_4$ are open intervals and hence they do not contain the intersection point $(\xi,\zeta)$.

In order to devise the explicit formulas for the coefficients of the compact scheme, we need some notations and definitions. Recall that
$a_{p}:=a \chi_{\Omega_{p}}$ is a positive constant, $f_{p}:=f \chi_{\Omega_{p}}$, and $u_{p}:=u \chi_{\Omega_{p}}$ for $p=1,2,3,4$.
For $(x_i^*,y_j^*) \in  \overline{\Omega_p}$ with $p=1,2,3,4$,
using only the information of $u_p$ in $\overline{\Omega_p}$,
we can define
\[
u_{p}^{(m,n)}:=\frac{\partial^{m+n} u_p}{ \partial^m x \partial^n y}(x_i^*,y_j^*),\qquad
f_{p}^{(m,n)}:=\frac{\partial^{m+n} f_p}{ \partial^m x \partial^n y}(x_i^*,y_j^*).
\]
For $(x_i^*,y_j^*)=(\xi, y_j^*)\in \Gamma_p$ with $p=1,2$, we define
$\phi_p^{(n)}:=\frac{d^{n} \phi_p(\xi, y)}{d^n y}|_{y=y_j^*}$ and $\psi_p^{(n)}:=\frac{d^{n} \psi_p(\xi, y)}{ d^n y}|_{y=y_j^*}$,
while for $(x_i^*,y_j^*)=(x_i^*,\zeta)\in {\Gamma_p}$ with $p=3,4$, we similarly define
$\phi_p^{(n)}:=\frac{d^{n} \phi_p(x,\zeta)}{ d^n x}|_{x=x_i^*}$ and $\psi_p^{(n)}:=\frac{d^{n} \psi_p(x,\zeta)}{ d^n x}|_{x=x_i^*}$.
Note that $\phi_p,\psi_p$ in \eqref{intersect:3} are essentially 1D functions defined on the line segment $\Gamma_p$.

Define $\NN:=\N\cup\{0\}$ and for $M\in \NN$,
we define the following index subsets of $\NN^2$ as follows:
\be \label{Intersect:Sk}
\ind_{M}:=\{(m,n)\in \N_0^2 \; : \; m+n\le M\},\quad \ind_M^{1}:=\{(m,n)\in \ind_M \; :\; m=0,1\},\quad \ind_M^{2}:=\ind_M\backslash \ind_M^{1}.
\ee
The illustrations for $\ind_{7}^{1}$ and $\ind_{7}^{2}$
%$\ind_{7}^{H, 1}$, $\ind_{7}^{H, 2}$
are shown in \cref{Intersect:fig:umnV}
%,Intersect:fig:umnH}
in \cref{Intersect:sec:proofs}.
We shall also define the following bivariate polynomials (which will be used in our compact FDMs later):
\be \label{Intersect:GmnV}
G_{M,m,n}(x,y):=
\sum_{\ell=0}^{\lfloor \frac{n}{2}\rfloor}
\frac{(-1)^\ell x^{m+2\ell} y^{n-2\ell}}{(m+2\ell)!(n-2\ell)!},\qquad  H_{M,m,n}(x,y):=\sum_{\ell=1}^{1+\lfloor \frac{n}{2}\rfloor} \frac{(-1)^{\ell} x^{m+2\ell} y^{n-2\ell+2}}{(m+2\ell)!(n-2\ell+2)!},
\ee
where $\lfloor x\rfloor$ is the floor function representing the largest integer less than or equal to $x\in \R$.

To state our compact FDMs later, we shall use some auxiliary polynomials of $h$. For $w\in \R$, $M\in \NN$ and $m,n\in \NN$, we define the univariate polynomials of $h$ with the parameter $w$ for the vertical interface $\Gamma_1$ or $\Gamma_2$ to be
\be \label{GHw:v}
\begin{split}
&G^{w}_{M,m,n}:= \sum_{\ell=-1}^1 C_{-1,\ell} G_{M,m,n}(wh-h, \ell h),\\
&H^{-,w}_{M,m,n}:=\sum_{\ell=-1}^1 C_{-1,\ell}  H_{M,m,n}(wh-h, \ell h), \quad
H^{+,w}_{M,m,n}:= \sum_{k=0}^1 \sum_{\ell=-1}^1 C_{k,\ell}  H_{M,m,n}(wh+kh, \ell h).
\end{split}
\ee
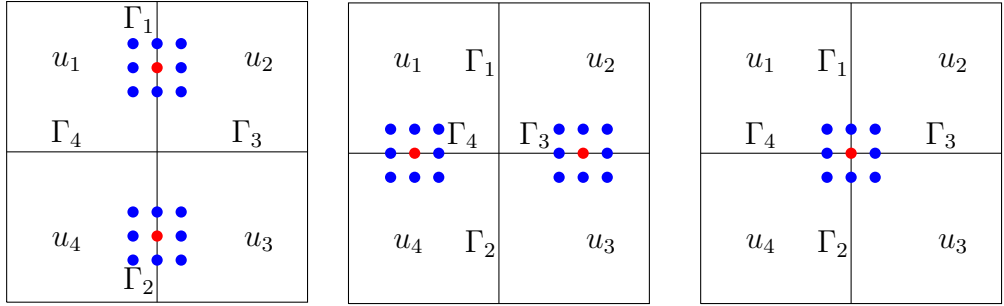
\begin{figure}[tbhp]
	\centering	
	\hspace{1.2cm}
	\begin{subfigure}[b]{0.2\textwidth}
		\hspace{-1.8cm}
		\begin{tikzpicture}[scale = 0.8]
			\draw	(-2.5, -2.5) -- (-2.5, 2.5) -- (2.5, 2.5) -- (2.5, -2.5) --(-2.5,-2.5);	
			\draw   (0, -2.5) -- (0, 2.5);	
			\draw   (-2.5, 0) -- (2.5, 0);	
			\node (A) at (-1.5,1.5) {$u_1$};
			\node (A) at (1.7,1.5) {$u_2$};
			\node (A) at (1.7,-1.5) {$u_3$};
			\node (A) at (-1.5,-1.5) {$u_4$};
			 %%%%%%%%%%%%%%%%%%%%%%%%%%%%%%%%%%%%%%%%%%%%%%%%%%%
			\node at (0,1)[circle,fill,inner sep=1.5pt,color=blue]{};
			\node at (0,1.4)[circle,fill,inner sep=1.5pt,color=red]{};
			\node at (0,1.8)[circle,fill,inner sep=1.5pt,color=blue]{};	
			\node at (0.4,1)[circle,fill,inner sep=1.5pt,color=blue]{};
			\node at (0.4,1.4)[circle,fill,inner sep=1.5pt,color=blue]{};
			\node at (0.4,1.8)[circle,fill,inner sep=1.5pt,color=blue]{};	
			\node at (-0.4,1)[circle,fill,inner sep=1.5pt,color=blue]{};
			\node at (-0.4,1.4)[circle,fill,inner sep=1.5pt,color=blue]{};
			\node at (-0.4,1.8)[circle,fill,inner sep=1.5pt,color=blue]{};	
			 %%%%%%%%%%%%%%%%%%%%%%%%%%%%%%%%%%%%%%%%%%%%%%%%
			\node at (0,-1)[circle,fill,inner sep=1.5pt,color=blue]{};
			\node at (0,-1.4)[circle,fill,inner sep=1.5pt,color=red]{};
			\node at (0,-1.8)[circle,fill,inner sep=1.5pt,color=blue]{};	
			\node at (0.4,-1)[circle,fill,inner sep=1.5pt,color=blue]{};
			\node at (0.4,-1.4)[circle,fill,inner sep=1.5pt,color=blue]{};
			\node at (0.4,-1.8)[circle,fill,inner sep=1.5pt,color=blue]{};	
			\node at (-0.4,-1)[circle,fill,inner sep=1.5pt,color=blue]{};
			\node at (-0.4,-1.4)[circle,fill,inner sep=1.5pt,color=blue]{};
			\node at (-0.4,-1.8)[circle,fill,inner sep=1.5pt,color=blue]{};	
			 %%%%%%%%%%%%%%%%%%%%%%%%%%%%%%%%%%%%%%%%%%%%%%%%%
			\node (A) at (-0.3,-2.14) {$\Gamma_2$};
			\node (A) at (-0.3,2.2) {$\Gamma_1$};
			\node (A) at (1.5,0.3) {$\Gamma_3$};
			\node (A) at (-1.5,0.3) {$\Gamma_4$};
		\end{tikzpicture}
	\end{subfigure}
	 %%%%%%%%%%%%%%%%%%%%%%%%%%%%%%%%%%%%%
	\hspace{-1cm}
	\begin{subfigure}[b]{0.2\textwidth}
		\begin{tikzpicture}[scale = 0.8]
			\draw	(-2.5, -2.5) -- (-2.5, 2.5) -- (2.5, 2.5) -- (2.5, -2.5) --(-2.5,-2.5);	
			\draw   (0, -2.5) -- (0, 2.5);	
			\draw   (-2.5, 0) -- (2.5, 0);	
			\node (A) at (-1.5,1.5) {$u_1$};
			\node (A) at (1.7,1.5) {$u_2$};
			\node (A) at (1.7,-1.5) {$u_3$};
			\node (A) at (-1.5,-1.5) {$u_4$};
			 %%%%%%%%%%%%%%%%%%%%%%%%%%%%%%%%%%%%%%%%%%%%%%%%%%%
			\node at (-1,0)[circle,fill,inner sep=1.5pt,color=blue]{};
			\node at (-1.4,0)[circle,fill,inner sep=1.5pt,color=red]{};
			\node at (-1.8,0)[circle,fill,inner sep=1.5pt,color=blue]{};	
			\node at (-1,0.4)[circle,fill,inner sep=1.5pt,color=blue]{};
			\node at (-1.4,0.4)[circle,fill,inner sep=1.5pt,color=blue]{};
			\node at (-1.8,0.4)[circle,fill,inner sep=1.5pt,color=blue]{};	
			\node at (-1,-0.4)[circle,fill,inner sep=1.5pt,color=blue]{};
			\node at (-1.4,-0.4)[circle,fill,inner sep=1.5pt,color=blue]{};
			\node at (-1.8,-0.4)[circle,fill,inner sep=1.5pt,color=blue]{};				
			\node at (1,0)[circle,fill,inner sep=1.5pt,color=blue]{};
			\node at (1.4,0)[circle,fill,inner sep=1.5pt,color=red]{};
			\node at (1.8,0)[circle,fill,inner sep=1.5pt,color=blue]{};	
			\node at (1,0.4)[circle,fill,inner sep=1.5pt,color=blue]{};
			\node at (1.4,0.4)[circle,fill,inner sep=1.5pt,color=blue]{};
			\node at (1.8,0.4)[circle,fill,inner sep=1.5pt,color=blue]{};	
			\node at (1,-0.4)[circle,fill,inner sep=1.5pt,color=blue]{};
			\node at (1.4,-0.4)[circle,fill,inner sep=1.5pt,color=blue]{};
			\node at (1.8,-0.4)[circle,fill,inner sep=1.5pt,color=blue]{};	
			 %%%%%%%%%%%%%%%%%%%%%%%%%%%%%%%%%%%%%%%%%%%%%%%%%%%
			 %%%%%%%%%%%%%%%%%%%%%%%%%%%%%%%%%%%%%%%%%%%%%%%%%
			\node (A) at (-0.3,-1.5) {$\Gamma_2$};
			\node (A) at (-0.3,1.5) {$\Gamma_1$};
			\node (A) at (0.6,0.3) {$\Gamma_3$};
			\node (A) at (-0.6,0.3) {$\Gamma_4$};
		\end{tikzpicture}
	\end{subfigure}
	 %%%%%%%%%%%%%%%%%%%%%%%%%%%%%%%%%%%%%%%%%%%%%%%%%%%%%%%%%%%%%%%%%%%%%%%%%%%%%%%%%%%%%%%%%
	\hspace{0.8cm}
	\begin{subfigure}[b]{0.2\textwidth}
	\begin{tikzpicture}[scale = 0.8]
		\draw	(-2.5, -2.5) -- (-2.5, 2.5) -- (2.5, 2.5) -- (2.5, -2.5) --(-2.5,-2.5);	
		\draw   (0, -2.5) -- (0, 2.5);	
		\draw   (-2.5, 0) -- (2.5, 0);	
		\node (A) at (-1.5,1.5) {$u_1$};
		\node (A) at (1.7,1.5) {$u_2$};
		\node (A) at (1.7,-1.5) {$u_3$};
		\node (A) at (-1.5,-1.5) {$u_4$};
		 %%%%%%%%%%%%%%%%%%%%%%%%%%%%%%%%%%%%%%%%%%%%%%%%%%%
		\node at (0,0.4)[circle,fill,inner sep=1.5pt,color=blue]{};
		\node at (0,0)[circle,fill,inner sep=1.5pt,color=red]{};
		\node at (0,-0.4)[circle,fill,inner sep=1.5pt,color=blue]{};	
		\node at (0.4,0.4)[circle,fill,inner sep=1.5pt,color=blue]{};
		\node at (0.4,0)[circle,fill,inner sep=1.5pt,color=blue]{};
		\node at (0.4,-0.4)[circle,fill,inner sep=1.5pt,color=blue]{};
		\node at (-0.4,0.4)[circle,fill,inner sep=1.5pt,color=blue]{};
		\node at (-0.4,0)[circle,fill,inner sep=1.5pt,color=blue]{};
		\node at (-0.4,-0.4)[circle,fill,inner sep=1.5pt,color=blue]{};
		 %%%%%%%%%%%%%%%%%%%%%%%%%%%%%%%%%%%%%%%%%%
		 %%%%%%%%%%%%%%%%%%%%%%%%%%%%%%%%%%%%%%%%%%%%%%%%%
		\node (A) at (-0.3,-1.5) {$\Gamma_2$};
		\node (A) at (-0.3,1.5) {$\Gamma_1$};
		\node (A) at (1.5,0.3) {$\Gamma_3$};
		\node (A) at (-1.5,0.3) {$\Gamma_4$};
	\end{tikzpicture}
\end{subfigure}
	\caption{ First panel:  compact  9-point  stencils in \cref{Intersect:thm:side1} with $(x_i,y_j)=(x_i^*,y_j^*)\in \Gamma_1 \cup \Gamma_2$. Second panel:  compact  9-point  stencils with $(x_i,y_j)=(x_i^*,y_j^*)\in \Gamma_3 \cup \Gamma_4$. Third panel: compact  9-point  stencil in \cref{Intersect:thm:cross1} with $(x_i,y_j)=(x_i^*,y_j^*)=(\xi,\zeta)$. The grid point $(x_i,y_j)$ is indicated by the red color.
}
	\label{fig:intersect:side1}
\end{figure}

We first consider the compact  discretization at grid points lying on the vertical interface line $\Gamma_1$ or $\Gamma_2$.  The modification
of the scheme corresponding to the horizontal interfaces $\Gamma_3$ or $\Gamma_4$ is straightforward and we will only briefly mention it afterwards.

\begin{theorem}\label{Intersect:thm:side1}
Consider a grid point  $(x_{i},y_{j})$ such that $(x_i,y_j)\in \Gamma_1$ (see the first panel of \cref{fig:intersect:side1}).  Assume that $u_{p}:=u \chi_{\Omega_{p}}$
and $f_{p}:=f \chi_{\Omega_{p}}$ have uniformly continuous partial derivatives of (total) orders up to seven and five, respectively, in each $\Omega_p$ for $p=1,2$. Also assume that
the essentially one-dimensional functions $\phi_1$ and $\psi_1$  on the interface $\Gamma_1$ have uniformly continuous derivatives of orders up to seven and six, respectively.
Let $(x_i^*,y_j^*):=(x_i,y_j)$ and
$\mathcal{L}_h$ be the discrete operator in \eqref{Luh} with the stencil coefficients
\be \label{C:Gamma1}
\begin{split}
&C_{1,0}=-4,\quad C_{1,-1}=C_{1,1}=-1,\quad
C_{-1,-1}=C_{-1,1}=-\alpha,\quad
 C_{0,-1}=C_{0,1}=-2(1+\alpha),\quad\\
&C_{-1,0}=-4\alpha,\quad C_{0,0}=10(1+\alpha),\quad \alpha:=a_1/a_2>0.
\end{split}
\ee
Then the compact 9-point finite difference scheme $h^{-1}\mathcal{L}_h u_h=h^{-1}F$
has a seventh order of consistency  at  the grid point $(x_i,y_j)\in \Gamma_1$,
%i.e., $h^{-1}\mathcal{L}_h(u_h-u)=\bo(h^7)$ as $h\to 0$,
where
\be\label{intersect:them1:irre:side1}
F:=\frac{1}{a_1}\sum_{(m,n)\in \ind_{5}} f_1^{(m,n)}
H^{-,0}_{7,m,n} +\frac{1}{a_2}\sum_{(m,n)\in \ind_{5}} f_2^{(m,n)}
H^{+,0}_{7,m,n}-\sum_{n=0}^{7} \phi_1^{(n)} G^0_{7,0,n}
-\frac{1}{a_1} \sum_{n=0}^{6} \psi_1^{(n)}  G^0_{7,1,n}
\ee
and $H^{-,0}_{7,m,n}, H^{+,0}_{7,m,n},
G^0_{7,0,n},
G^0_{7,1,n}$ are defined in \eqref{GHw:v}.
\end{theorem}

For $(x_i,y_j)\in \Gamma_2$  (see the first panel of \cref{fig:intersect:side1}), we can obtain the compact  9-point  finite difference scheme with a consistency order seven  by using $\alpha:={a_4}/{a_3}$, replacing the subscript $1$ by $2$ for $\Gamma_1$, $\phi_1$, and $\psi_1$, and replacing  all the subscripts $1$ and $2$ by $4$ and $3$ for $f_1$ and $f_2$, respectively, in \cref{Intersect:thm:side1}.

Similarly to \cref{Intersect:thm:side1},
we can specify the scheme for grid points lying on the horizontal line interface $\Gamma_3$ (see the second panel of \cref{fig:intersect:side1})
by modifying the stencil coefficients to:%
%
%\begin{theorem}\label{Intersect:thm:side3}
%Consider a grid point  $(x_{i},y_{j})$ such that $(x_i,y_j)\in \Gamma_3$ (see the second panel of \cref{fig:intersect:side1}).
%Let $(x_{i}^*,y_{j}^*):=(x_i,y_j)$ and
%$\mathcal{L}_h$ be the discrete operator in \eqref{Luh} with the stencil coefficients
%
\be \label{C:Gamma3}
\begin{split}
&C_{0,1}=-4,\quad C_{-1,1}=C_{1,1}=-1,\quad
C_{-1,-1}=C_{1,-1}=-\alpha,\quad
 C_{-1,0}=C_{1,0}=-2(1+\alpha),\quad\\
&C_{0,-1}=-4\alpha,\quad C_{0,0}=10(1+\alpha),\quad \alpha:=a_3/a_2>0.
\end{split}
\ee
The right hand side vector in this case is given by:
\be\label{F:Gamma:3}
F:=\frac{1}{a_3}\sum_{(m,n)\in \ind_{5}} f_3^{(m,n)}
\tilde{H}^{-,0}_{7,m,n} +\frac{1}{a_2} \sum_{(m,n)\in \ind_{5}} f_2^{(m,n)}
\tilde{H}^{+,0}_{7,m,n}-\sum_{m=0}^{7} \phi_3^{(m)}  \tilde{G}^0_{7,m,0} -\frac{1}{a_3} \sum_{m=0}^{6} \psi_3^{(m)}  \tilde{G}^0_{7,m,1}
\ee
%and $\tilde{H}^{-,0}_{7,m,n}, \tilde{H}^{+,0}_{7,m,n}, \tilde{G}^0_{7,m,0}, \tilde{G}^0_{7,m,1}$ are defined in \eqref{GHw:h}.
%\end{theorem}
where the polynomials for the horizontal interface are defined as:
\be \label{GHw:h}
\begin{split}
&\tilde{G}^{w}_{M,m,n}:=\sum_{k=-1}^1 C_{k,-1} G_{M,n,m}(wh-h, k h),\\
&\tilde{H}^{-,w}_{M,m,n}:=\sum_{k=-1}^1 C_{k,-1}  H_{M,n,m}(wh-h, k h), \quad
\tilde{H}^{+,w}_{M,m,n}:= \sum_{k=-1}^1 \sum_{\ell=0}^1 C_{k,\ell}  H_{M,n,m}(wh+\ell h, k h).
\end{split}
\ee

For	$(x_i,y_j)\in \Gamma_4$ (see the second panel of \cref{fig:intersect:side1}),
we can obtain the scheme by using $\alpha=a_4/a_1$, replacing the subscript $3$ by $4$ for $\Gamma_3$, $\phi_3$, and $\psi_3$, and replacing  all the subscripts $2$ and $3$ by $1$ and $4$ for $f_2$ and $f_3$, respectively in \eqref{C:Gamma3}--\eqref{F:Gamma:3}.

Finally, we handle the case when the intersecting interface point $(\xi,\zeta)$ is a grid point.

\begin{theorem}\label{Intersect:thm:cross1}
Consider the grid point  $(x_i,y_j)=(\xi,\zeta)$
(see the third panel of \cref{fig:intersect:side1}).  Assume that $u_{p}:=u \chi_{\Omega_{p}}$ and $f_{p}:=f \chi_{\Omega_{p}}$ have uniformly continuous partial derivatives of (total) orders up to seven and five, respectively, in each $\Omega_p$ for $p=1,2,3,4$. Also assume that
	the essentially one-dimensional functions $\phi_p$ and $\psi_p$  on the interface $\Gamma_p$ have uniformly continuous derivatives of orders up to seven and six, respectively, for $p=1,2,3,4$.
Let $(x_{i}^*,y_{j}^*):=(\xi,\zeta)$ and
$\mathcal{L}_h$ be the discrete operator in \eqref{Luh} with the stencil coefficients
	\be\label{Intersect:CKL:cross1}
	\begin{aligned}	
		&     C_{-1,1}:=\frac{-a_1^2(a_2+a_3)}{a_2^2(a_1+a_4)},    \qquad &
		&      C_{0,1}:=\frac{-2(a_1+a_2)}{a_2},   \qquad &
		& C_{1,1}:=-1,\\
		 %%%%%%%%%%%%%%%%%%%%%%%%%%%%%%%%%%%%%%%%%%
		&  C_{-1,0}:=\frac{-2a_1(a_2+a_3)}{a_2^2},   \qquad &
		&    C_{0,0}:=\frac{5(a_2+a_3)(a_1+a_2)}{a_2^2},  \qquad &
		&    C_{1,0}:=\frac{-2(a_2+a_3)}{a_2},   \\
		 %%%%%%%%%%%%%%%%%%%%%%%%%%%%%%%%%%%%%%%
		&  C_{-1,-1}:=\frac{-a_1a_4(a_2+a_3)}{a_2^2(a_1+a_4)},\qquad  &
& C_{0,-1}:= \frac{-2a_3(a_1+a_2)}{a_2^2} ,\qquad  &
&  C_{1,-1} :=\frac{-a_3}{a_2}.
\end{aligned}
\ee
Then the compact 9-point finite difference scheme $h^{-1}\mathcal{L}_h u_h=h^{-1}F$
has a seventh order of consistency at the grid point $(x_i,y_j)$ (i.e., $(\xi,\zeta)$),
where
\be\label{intersect:them1:irre:cross1}
\begin{aligned}
F:=	& \sum_{(n,m)\in \ind_{7}^{1}} F^{1,0,0}_{7,m,n}
		+ \sum_{(m,n)\in \ind_{5}} F^{0,0}_{7,m,n}  +   \sum_{(n,m)\in \ind_{7}^{1}} \Phi^{1,0,0}_{7,m,n}+ \sum_{m=0}^{7}  \Phi^{0,0}_{7,m}
	+  \sum_{(n,m)\in \ind_{7}^{1}} \Psi^{1,0,0}_{7,m,n}+ \sum_{m=0}^{6}  \Psi^{0,0}_{7,m}
	\end{aligned}
	\ee
and $F^{1,0,0}_{7,m,n}$, $F^{0,0}_{7,m,n}$, $\Phi^{1,0,0}_{7,m,n}$, $\Phi^{0,0}_{7,m}$, $\Psi^{1,0,0}_{7,m,n}$, $\Psi^{0,0}_{7,m}$ are defined in \eqref{intersect:FH1mn:Cross}--\eqref{intersect:PsiH1Mmn:Cross}.
\end{theorem}

\subsection{Compact 9-point stencils at grid points neighboring the interface}\label{neighboring:interface}

In this subsection, we derive a compact scheme with a consistency order four  for every grid point neighboring the interface $\Gamma:=\Gamma_1 \cup \Gamma_2 \cup \Gamma_3 \cup \Gamma_4\cup\{(\xi,\zeta)\}$ (i.e., Case 2 in \cref{introdu:1:intersect}, see the right panel of \cref{compact:intersect:model}).
%To present our compact   FDM, we shall also use the notations and definitions in \cref{Intersect:subsec:irregular}.

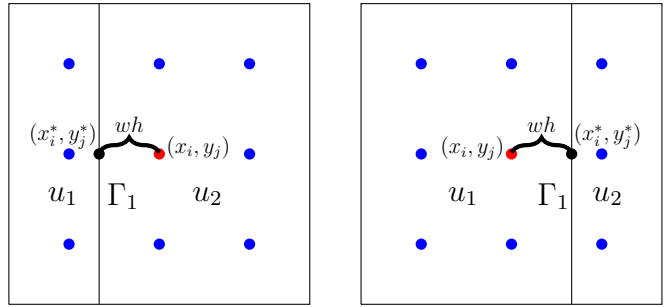
\begin{figure}[htbp]
	\centering	
	 %%%%%%%%%%%%%%%%%%%%%%%%%%%%%%%%%%%%%%%%%%%%%%%%%%%%%%%%%%%%%%%%%%%%%%%%%%%%%%%%%%%%%%%%%
	\hspace{-0.8cm}
	\begin{subfigure}[b]{0.2\textwidth}
		\begin{tikzpicture}[scale = 0.8]
			\draw	(-2.5, -2.5) -- (-2.5, 2.5) -- (2.5, 2.5) -- (2.5, -2.5) --(-2.5,-2.5);		
			\draw   (-1, -2.5) -- (-1, 2.5);
			 %%%%%%%%%%%%%%%%%%%%%%%%%%%%%%%%%%%%%%%%%%%%%%%%%%%
			\node at (0,1.5)[circle,fill,inner sep=1.5pt,color=blue]{};
			\node at (0,0)[circle,fill,inner sep=1.5pt,color=red]{};
			\node at (0,-1.5)[circle,fill,inner sep=1.5pt,color=blue]{};	
			\node at (1.5,1.5)[circle,fill,inner sep=1.5pt,color=blue]{};
			\node at (1.5,0)[circle,fill,inner sep=1.5pt,color=blue]{};
			\node at (1.5,-1.5)[circle,fill,inner sep=1.5pt,color=blue]{};
			\node at (-1.5,1.5)[circle,fill,inner sep=1.5pt,color=blue]{};
			\node at (-1.5,0)[circle,fill,inner sep=1.5pt,color=blue]{};
			\node at (-1.5,-1.5)[circle,fill,inner sep=1.5pt,color=blue]{};
			 %%%%%%%%%%%%%%%%%%%%%%%%%%%%%%%%%%%%%%%%%%
			%%%	\draw   (0,0) -- (0,-1);
			\node (A)[scale=0.7] at (-0.5,0.5) {$wh$};
			 \draw[decorate,decoration={brace,mirror,amplitude=2mm},xshift=0pt,yshift=10pt,ultra thick] (0,-0.33) -- node [black,midway,yshift=0.6cm]{} (-1,-0.33);
			%%%%%%%%%%%%%%%%%%
			\node at (-1,0)[circle,fill,inner sep=1.5pt,color=black]{};
			\node (A)[scale=0.7] at (0.65,0.1) {$(x_i,y_j)$};
			\node (A)[scale=0.7] at (-1.6,0.3) {$(x_i^*,y_j^*)$};
			%%%%%%%%%%%%%%%%
			\node (A) at (-1.6,-0.7) {$u_1$};
			\node (A) at (0.8,-0.7) {$u_2$};
			\node (A) at (-0.6,-0.7) {$\Gamma_1$};
		\end{tikzpicture}
	\end{subfigure}
	%%%%%%%%%%%%%%%%%%%%%
	\hspace{0.8cm}
	\begin{subfigure}[b]{0.2\textwidth}
		\begin{tikzpicture}[scale = 0.8]
			\draw	(-2.5, -2.5) -- (-2.5, 2.5) -- (2.5, 2.5) -- (2.5, -2.5) --(-2.5,-2.5);		
			\draw   (1, -2.5) -- (1, 2.5);
			 %%%%%%%%%%%%%%%%%%%%%%%%%%%%%%%%%%%%%%%%%%%%%%%%%%%
			\node at (0,1.5)[circle,fill,inner sep=1.5pt,color=blue]{};
			\node at (0,0)[circle,fill,inner sep=1.5pt,color=red]{};
			\node at (0,-1.5)[circle,fill,inner sep=1.5pt,color=blue]{};	
			\node at (1.5,1.5)[circle,fill,inner sep=1.5pt,color=blue]{};
			\node at (1.5,0)[circle,fill,inner sep=1.5pt,color=blue]{};
			\node at (1.5,-1.5)[circle,fill,inner sep=1.5pt,color=blue]{};
			\node at (-1.5,1.5)[circle,fill,inner sep=1.5pt,color=blue]{};
			\node at (-1.5,0)[circle,fill,inner sep=1.5pt,color=blue]{};
			\node at (-1.5,-1.5)[circle,fill,inner sep=1.5pt,color=blue]{};
			 %%%%%%%%%%%%%%%%%%%%%%%%%%%%%%%%%%%%%%%%%%
			%%%	\draw   (0,0) -- (0,-1);
			\node (A)[scale=0.7] at (0.5,0.5) {$wh$};
			 \draw[decorate,decoration={brace,amplitude=2mm},xshift=0pt,yshift=10pt,ultra thick] (0,-0.33) -- node [black,midway,yshift=0.6cm]{} (1,-0.33);
			%%%%%%%%%%%%%%%%%%
			\node at (1,0)[circle,fill,inner sep=1.5pt,color=black]{};
			\node (A)[scale=0.7] at (-0.65,0.1) {$(x_i,y_j)$};
			\node (A)[scale=0.7] at (1.6,0.3) {$(x_i^*,y_j^*)$};
			%%%%%%%%%%%%%%%%
			\node (A) at (1.6,-0.7) {$u_2$};
			\node (A) at (-0.8,-0.7) {$u_1$};
			\node (A) at (0.7,-0.7) {$\Gamma_1$};
		\end{tikzpicture}
	\end{subfigure}
	\caption{  Left: compact  9-point  stencil in \cref{theorem:side3:w} with $(x_i^*,y_j^*):=(x_i-wh,y_j)\in \Gamma_1$ and $0<w<1$. Right: compact 9-point stencil with $(x_i^*,y_j^*):=(x_i+wh,y_j)\in \Gamma_1$ and $0<w<1$. The grid point $(x_i,y_j)$ is indicated by the red color, while the interface point $(x_i^*,y_j^*)\in \Gamma_1$ is indicated by the black color and is not a grid point.}
	\label{fig:intersect:side1:w}
\end{figure}

If the grid point $(x_i,y_j)\not\in \Gamma$ is in the vicinity of the interface line $\Gamma_1$ (see the first  panel of \cref{fig:intersect:side1:w}),
the compact scheme with a consistency order four is given in the theorem below.

\begin{theorem}\label{theorem:side3:w}
	Consider a grid point $(x_{i},y_{j})\not\in \Gamma$ such that
$(x_i^*,y_j^*):=(x_i-wh,y_j)\in \Gamma_1$ with $0<w<1$, $(\xi,\zeta)\not\in (x_i,y_j)+(-h,h)^2$ (see the first panel of \cref{fig:intersect:side1:w}).  Assume that $u_{p}:=u \chi_{\Omega_{p}}$ and $f_{p}:=f \chi_{\Omega_{p}}$ have uniformly continuous partial derivatives of (total) orders up to four and two, respectively, in each $\Omega_p$ for $p=1,2$. Also assume that
	the essentially one-dimensional functions $\phi_1$ and $\psi_1$  on the interface $\Gamma_1$ have uniformly continuous derivatives of orders up to four and three, respectively.
Let  $\mathcal{L}_h$ be the discrete operator in \eqref{Luh} with the stencil coefficients
{\small{
\be\label{Intersect:CKL:side3:w}
		\begin{aligned}
			&  C_{-1,1}:=
			-(r_4\alpha^2+r_5\alpha)/
			\beta,\quad  &
			& C_{0,1}:= -1,\quad  &
			&  C_{1,1} :=-(r_6+t_1\alpha^2
			+r_7\alpha)/\beta,\\
			 %%%%%%%%%%%%%%%%%%%%%%%%%%%%%%%%%%%
			&  C_{-1,0}:= -(r_8\alpha^2 +r_{9}\alpha)/\beta,   \quad &
			&  C_{0,0}:=-(s_{1} +s_{2}\alpha^2
			+s_{3}\alpha)/\beta,  \quad &
			&  C_{1,0}:=-(r_{10} +r_{11}\alpha^2
			+r_{12}\alpha)/\beta,   \\
			 %%%%%%%%%%%%%%%%%%%%%%%%%%%%%%%%%%%%%%%
			& C_{-1,-1}:=C_{-1,1},    \quad &
			& C_{0,-1}:=C_{0,1},   \quad &
			& C_{1,-1}:=C_{1,1},
		\end{aligned}
		\ee
}}
where $\alpha:=a_1/a_2>0$,
$\beta:=r_1 +r_2\alpha^2+r_3\alpha> 0$ and
{\tiny{
		\be\label{r1tor12:w:plus}
		\begin{split}
			&   r_1:=   (2w+1)^2(w+2)(w-1)^2,\     r_2:=   4w^5-4w^4+5w^3+6w^2-5w+2,\
r_3:=   -8w^5+6w^3-2w^2+4,  \\
& r_4:=    4w^4-4w^3+w^2+1,\     r_5:=    -4w^4+4w^3-w^2+1,\   r_6:=  -(2w+1)^2(w-1)^3,\    r_7:=  8w^5-20w^4+14w^3-3w^2+1,   \\
& r_8:=       -8w^4+8w^3+10w^2-6w+4,\       r_{9}:=       8w^4-8w^3-10w^2+6w+4,\       r_{10}:=       8w^5-16w^4+14w^3-8w^2-2w+4,    \\
& r_{11}:=       8w^5-24w^4+38w^3-22w^2+8w,\     r_{12}:=       -16w^5+40w^4-52w^3+30w^2-6w+4,\      t_1:=  -4w^5+12w^4-13w^3+8w^2-w,   \\
&    s_{1}:=       -8w^5-8w^4+10w^3+26w^2-10w-10,\   s_{2}:=       -8w^5+8w^4-22w^3-18w^2+10w-10,\    s_{3}:=       16w^5+12w^3-8w^2-20.
		\end{split}
		\ee}}
Then the compact 9-point finite difference  scheme $h^{-1}\mathcal{L}_h u_h=h^{-1}F$
has a fourth order of consistency  at the grid point $(x_i,y_j) \not\in \Gamma$ with $(\xi,\zeta)\not \in (x_i,y_j)+(-h,h)^2$,
where
	\be\label{LhGamma1uh:w:side1}
\begin{aligned}
F:=\frac{1}{a_1} \sum_{(m,n)\in \ind_{2}} f_1^{(m,n)} H^{-,w}_{4,m,n} +\frac{1}{a_2} \sum_{(m,n)\in \ind_{2}} f_2^{(m,n)}
	 H^{+,w}_{4,m,n}- \sum_{n=0}^{4} \phi_1^{(n)} G^w_{4,0,n}-\frac{1}{a_1} \sum_{n=0}^{3} \psi_1^{(n)}  G^w_{4,1,n}
\end{aligned}
\ee
and	 $H^{-,w}_{4,m,n}, H^{+,w}_{4,m,n}, G^w_{4,0,n}, G^w_{4,1,n}$ are defined in \eqref{GHw:v}.
Moreover, up to a multiplicative constant for normalization, the stencil coefficients $\{  C_{k,\ell}\}_{k,\ell=-1,0,1}$ in \eqref{Intersect:CKL:side3:w} are unique.	
\end{theorem}

If a grid point $(x_i,y_j)$ satisfies
$(x_i^*,y_j^*):=(x_i+wh,y_j)\in \Gamma_1$  (see the second panel of \cref{fig:intersect:side1:w}) or $(x_i^*,y_j^*):=(x_i\pm wh,y_j)\in \Gamma_2$  or
$(x_i^*,y_j^*):=(x_i,y_j\pm wh)\in \Gamma_3, \Gamma_4$ with $0<w<1$, then the compact  9-point scheme with a consistency order four at the grid point $(x_i,y_j)\not\in \Gamma$ can be obtained similarly.

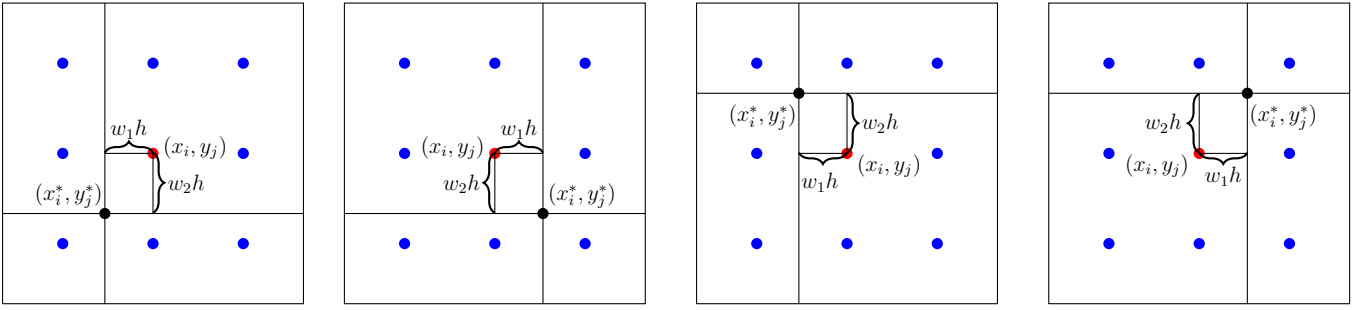
\begin{figure}[htbp]
	\centering	
	\hspace{1.2cm}
	\begin{subfigure}[b]{0.2\textwidth}
		\hspace{-1.8cm}
		\begin{tikzpicture}[scale = 0.8]
			\draw	(-2.5, -2.5) -- (-2.5, 2.5) -- (2.5, 2.5) -- (2.5, -2.5) --(-2.5,-2.5);	
			\draw   (-0.8, -2.5) -- (-0.8, 2.5);	 
			\draw   (-2.5, -1) -- (2.5, -1);
			 %%%%%%%%%%%%%%%%%%%%%%%%%%%%%%%%%%%%%%%%%%%%%%%%%%%
			\node at (0,1.5)[circle,fill,inner sep=1.5pt,color=blue]{};
			\node at (0,0)[circle,fill,inner sep=1.5pt,color=red]{};
			\node at (0,-1.5)[circle,fill,inner sep=1.5pt,color=blue]{};	
			\node at (1.5,1.5)[circle,fill,inner sep=1.5pt,color=blue]{};
			\node at (1.5,0)[circle,fill,inner sep=1.5pt,color=blue]{};
			\node at (1.5,-1.5)[circle,fill,inner sep=1.5pt,color=blue]{};
			\node at (-1.5,1.5)[circle,fill,inner sep=1.5pt,color=blue]{};
			\node at (-1.5,0)[circle,fill,inner sep=1.5pt,color=blue]{};
			\node at (-1.5,-1.5)[circle,fill,inner sep=1.5pt,color=blue]{};
			 %%%%%%%%%%%%%%%%%%%%%%%%%%%%%%%%%%%%%%%%%%
			\draw   (0,0) -- (-0.8,0);
			\node (A)[scale=0.7] at (-0.4,0.4) {$w_1h$};
			 \draw[decorate,decoration={brace,mirror,amplitude=1.5mm},xshift=0pt,yshift=10pt,thick] (0,-0.34) -- node [black,midway,yshift=0.6cm]{} (-0.8,-0.34);
			 %%%%%%%%%%%%%%%%%%%%%%%%%%%%%%%%%%%
			\draw   (0,0) -- (0,-1);
			\node (A)[scale=0.7] at (0.55,-0.5) {$w_2h$};
			 \draw[decorate,decoration={brace,amplitude=1.5mm},xshift=0pt,yshift=10pt, thick] (0,-0.34) -- node [black,midway,yshift=0.6cm]{} (0,-1.34);
			 %%%%%%%%%%%%%%%%%%%%%%%%%%%%%%%%%%%%
			\node at (-0.8,-1)[circle,fill,inner sep=1.5pt,color=black]{};
			\node (A)[scale=0.7] at (0.7,0.1) {$(x_i,y_j)$};
			\node (A)[scale=0.7] at (-1.4,-0.7) {$(x_i^*,y_j^*)$};
		\end{tikzpicture}
	\end{subfigure}
	 %%%%%%%%%%%%%%%%%%%%%%%%%%%%%%%%%%%%%
	\hspace{-1cm}
	\begin{subfigure}[b]{0.2\textwidth}
		\begin{tikzpicture}[scale = 0.8]
			\draw	(-2.5, -2.5) -- (-2.5, 2.5) -- (2.5, 2.5) -- (2.5, -2.5) --(-2.5,-2.5);	
			\draw   (0.8, -2.5) -- (0.8, 2.5);	
			\draw   (-2.5, -1) -- (2.5, -1);
			 %%%%%%%%%%%%%%%%%%%%%%%%%%%%%%%%%%%%%%%%%%%%%%%%%%%
			\node at (0,1.5)[circle,fill,inner sep=1.5pt,color=blue]{};
			\node at (0,0)[circle,fill,inner sep=1.5pt,color=red]{};
			\node at (0,-1.5)[circle,fill,inner sep=1.5pt,color=blue]{};	
			\node at (1.5,1.5)[circle,fill,inner sep=1.5pt,color=blue]{};
			\node at (1.5,0)[circle,fill,inner sep=1.5pt,color=blue]{};
			\node at (1.5,-1.5)[circle,fill,inner sep=1.5pt,color=blue]{};
			\node at (-1.5,1.5)[circle,fill,inner sep=1.5pt,color=blue]{};
			\node at (-1.5,0)[circle,fill,inner sep=1.5pt,color=blue]{};
			\node at (-1.5,-1.5)[circle,fill,inner sep=1.5pt,color=blue]{};
			 %%%%%%%%%%%%%%%%%%%%%%%%%%%%%%%%%%%%%%%%%%
			\draw   (0,0) -- (0.8,0);
			\node (A)[scale=0.7] at (0.4,0.4) {$w_1h$};
			 \draw[decorate,decoration={brace,amplitude=1.5mm},xshift=0pt,yshift=10pt, thick] (0,-0.34) -- node [black,midway,yshift=0.6cm]{} (0.8,-0.34);
			 %%%%%%%%%%%%%%%%%%%%%%%%%%%%%%%%%%%
			\draw   (0,0) -- (0,-1);
			\node (A)[scale=0.7] at (-0.55,-0.5) {$w_2h$};
			 \draw[decorate,decoration={brace,mirror,amplitude=1.5mm},xshift=0pt,yshift=10pt, thick] (0,-0.34) -- node [black,midway,yshift=0.6cm]{} (0,-1.34);
			 %%%%%%%%%%%%%%%%%%%%%%%%%%%%%%%%%%%%
			\node at (0.8,-1)[circle,fill,inner sep=1.5pt,color=black]{};
			\node (A)[scale=0.7] at (-0.7,0.1) {$(x_i,y_j)$};
			\node (A)[scale=0.7] at (1.45,-0.7) {$(x_i^*,y_j^*)$};
		\end{tikzpicture}
	\end{subfigure}
	 %%%%%%%%%%%%%%%%%%%%%%%%%%%%%%%%%%%%%%%%%%%%%%%%%%%%%%%%%%%%%%%%%%%%%%%%%%%%%%%%%%%%%%%%%
	\hspace{0.8cm}
	\begin{subfigure}[b]{0.2\textwidth}
		\begin{tikzpicture}[scale = 0.8]
			\draw	(-2.5, -2.5) -- (-2.5, 2.5) -- (2.5, 2.5) -- (2.5, -2.5) --(-2.5,-2.5);	
			\draw   (-0.8, -2.5) -- (-0.8, 2.5);	 
			\draw   (-2.5, 1) -- (2.5, 1);
			 %%%%%%%%%%%%%%%%%%%%%%%%%%%%%%%%%%%%%%%%%%%%%%%%%%%
			\node at (0,1.5)[circle,fill,inner sep=1.5pt,color=blue]{};
			\node at (0,0)[circle,fill,inner sep=1.5pt,color=red]{};
			\node at (0,-1.5)[circle,fill,inner sep=1.5pt,color=blue]{};	
			\node at (1.5,1.5)[circle,fill,inner sep=1.5pt,color=blue]{};
			\node at (1.5,0)[circle,fill,inner sep=1.5pt,color=blue]{};
			\node at (1.5,-1.5)[circle,fill,inner sep=1.5pt,color=blue]{};
			\node at (-1.5,1.5)[circle,fill,inner sep=1.5pt,color=blue]{};
			\node at (-1.5,0)[circle,fill,inner sep=1.5pt,color=blue]{};
			\node at (-1.5,-1.5)[circle,fill,inner sep=1.5pt,color=blue]{};
			 %%%%%%%%%%%%%%%%%%%%%%%%%%%%%%%%%%%%%%%%%%
			\draw   (0,0) -- (-0.8,0);
			\node (A)[scale=0.7] at (-0.45,-0.48) {$w_1h$};
			 \draw[decorate,decoration={brace,amplitude=1.5mm},xshift=0pt,yshift=10pt, thick] (0,-0.34) -- node [black,midway,yshift=0.6cm]{} (-0.8,-0.34);
			 %%%%%%%%%%%%%%%%%%%%%%%%%%%%%%%%%%%
			\draw   (0,0) -- (0,1);
			\node (A)[scale=0.7] at (0.55,0.5) {$w_2h$};
			 \draw[decorate,decoration={brace,mirror,amplitude=1.5mm},xshift=0pt,yshift=10pt, thick] (0,-0.34) -- node [black,midway,yshift=0.6cm]{} (0,0.64);
			 %%%%%%%%%%%%%%%%%%%%%%%%%%%%%%%%%%%%
            \node at (-0.8,1)[circle,fill,inner sep=1.5pt,color=black]{};
            \node (A)[scale=0.7] at (0.7,-0.2) {$(x_i,y_j)$};
            \node (A)[scale=0.7] at (-1.4,0.6) {$(x_i^*,y_j^*)$};			
		\end{tikzpicture}
	\end{subfigure}
	%%%%%%%%%%%%%%%%%%%%%
	\hspace{0.8cm}
	\begin{subfigure}[b]{0.2\textwidth}
		\begin{tikzpicture}[scale = 0.8]
			\draw	(-2.5, -2.5) -- (-2.5, 2.5) -- (2.5, 2.5) -- (2.5, -2.5) --(-2.5,-2.5);	
			\draw   (0.8, -2.5) -- (0.8, 2.5);	
			\draw   (-2.5, 1) -- (2.5, 1);
			 %%%%%%%%%%%%%%%%%%%%%%%%%%%%%%%%%%%%%%%%%%%%%%%%%%%
			\node at (0,1.5)[circle,fill,inner sep=1.5pt,color=blue]{};
			\node at (0,0)[circle,fill,inner sep=1.5pt,color=red]{};
			\node at (0,-1.5)[circle,fill,inner sep=1.5pt,color=blue]{};	
			\node at (1.5,1.5)[circle,fill,inner sep=1.5pt,color=blue]{};
			\node at (1.5,0)[circle,fill,inner sep=1.5pt,color=blue]{};
			\node at (1.5,-1.5)[circle,fill,inner sep=1.5pt,color=blue]{};
			\node at (-1.5,1.5)[circle,fill,inner sep=1.5pt,color=blue]{};
			\node at (-1.5,0)[circle,fill,inner sep=1.5pt,color=blue]{};
			\node at (-1.5,-1.5)[circle,fill,inner sep=1.5pt,color=blue]{};
			 %%%%%%%%%%%%%%%%%%%%%%%%%%%%%%%%%%%%%%%%%%
			\draw   (0,0) -- (0.8,0);
			\node (A)[scale=0.7] at (0.4,-0.38) {$w_1h$};
			 \draw[decorate,decoration={brace,mirror,amplitude=1.5mm},xshift=0pt,yshift=10pt, thick] (0,-0.34) -- node [black,midway,yshift=0.6cm]{} (0.8,-0.34);
			 %%%%%%%%%%%%%%%%%%%%%%%%%%%%%%%%%%%
			\draw   (0,0) -- (0,1);
			\node (A)[scale=0.7] at (-0.6,0.5) {$w_2h$};
			 \draw[decorate,decoration={brace,amplitude=1.5mm},xshift=0pt,yshift=10pt, thick] (0,-0.34) -- node [black,midway,yshift=0.6cm]{} (0,0.64);
			 %%%%%%%%%%%%%%%%%%%%%%%%%%%%%%%%%%%%
			\node at (0.8,1)[circle,fill,inner sep=1.5pt,color=black]{};
			\node (A)[scale=0.7] at (-0.7,-0.2) {$(x_i,y_j)$};
			\node (A)[scale=0.7] at (1.4,0.6) {$(x_i^*,y_j^*)$};	
		\end{tikzpicture}
	\end{subfigure}
	\caption{ Compact  9-point stencils neighboring the intersecting point $(\xi,\zeta)$ in \cref{Intersect:thm:cross1:w1w2}: $(x_{i}^*,y_{j}^*):= (\xi,\zeta)=(x_i\pm w_1h,y_j\pm w_2h)$ with $0<w_1<1$ and $0<w_2<1$.
The grid point $(x_i,y_j)$ is indicated by the red color, while the intersecting interface point $(x_i^*,y_j^*)$ (i.e., $(\xi,\zeta)$) is indicated by the black color and is not a grid point.}
	\label{fig:intersect:cross:point:w}
\end{figure}

If the grid point $(x_i,y_j)\not\in \Gamma$ is situated as shown in the first  panel of  \cref{fig:intersect:cross:point:w},
a compact scheme of a consistency order four can be specified as in the following theorem.

\begin{theorem}\label{Intersect:thm:cross1:w1w2}
Consider a grid point $(x_{i},y_{j})\not\in \Gamma$ such that $(x_{i}^*,y_{j}^*):=(\xi,\zeta)=(x_i-w_1h,y_j-w_2h)$
with $0<w_1,w_2<1$ (see the first panel of \cref{fig:intersect:cross:point:w}).   Assume that $u_{p}:=u \chi_{\Omega_{p}}$ and $f_{p}:=f \chi_{\Omega_{p}}$ have uniformly continuous partial derivatives of (total) orders up to four and two, respectively, in each $\Omega_p$ for $p=1,2,3,4$. Also assume that
	the essentially one-dimensional functions $\phi_p$ and $\psi_p$  on the interface $\Gamma_p$ have uniformly continuous derivatives of orders up to four and three, respectively, for $p=1,2,3,4$.
Let $\mathcal{L}_h$ be the discrete operator in \eqref{Luh} with the stencil coefficients
$\{{C}_{k,\ell} \}_{k,\ell=-1,0,1}$ (up to a multiplicative constant for normalization) being uniquely determined by solving the linear system specified by \eqref{sum:Ckl} and \eqref{solve:Imn:last}  with $M=4$.
Then the compact 9-point finite difference scheme $h^{-1}\mathcal{L}_h u_h=h^{-1}F$ has a fourth order of consistency  at the grid point $(x_i,y_j)\not \in \Gamma$ with $(\xi,\zeta)\in (x_i,y_j)+(-h,h)^2$,
where
\[
\begin{aligned}
F:=& \sum_{(n,m)\in \ind_{4}^{1}} \tilde{F}^{1,w_1,w_2}_{4,m,n}+ \sum_{(m,n)\in \ind_{2}} \tilde{F}^{w_1,w_2}_{4,m,n}
+ \sum_{(n,m)\in \ind_{4}^{1}} \tilde{\Phi}^{1,w_1,w_2}_{4,m,n}+ \sum_{m=0}^{4}  \tilde{\Phi}^{w_1,w_2}_{4,m}
+ \sum_{(n,m)\in \ind_{4}^{1}} \tilde{\Psi}^{1,w_1,w_2}_{4,m,n}+ \sum_{m=0}^{3}  \tilde{\Psi}^{w_1,w_2}_{4,m}
\end{aligned}
\]
and the quantities $\tilde{F}^{1,w_1,w_2}_{4,m,n}$, $\tilde{F}^{w_1,w_2}_{4,m,n}$, $\tilde{\Phi}^{1,w_1,w_2}_{4,m,n}$, $\tilde{\Phi}^{w_1,w_2}_{4,m}$, $\tilde{\Psi}^{1,w_1,w_2}_{4,m,n}$, $\tilde{\Psi}^{w_1,w_2}_{4,m}$ are given in
\eqref{Imn:Final:one}--\eqref{intersect:PsiMm:Cross:tilde}.
\end{theorem}

Compact  9-point  finite difference schemes with a consistency order four for the
other three cases (see the second to fourth panels of \cref{fig:intersect:cross:point:w}) can be obtained similarly.
Note that the solution to the linear system in  \eqref{sum:Ckl} and \eqref{solve:Imn:last} with $M=4$ can be readily obtained by direct computations or using symbolic software.

%%%%%%%%%%%%%%%%%%%%%%%%%%%%%%%%%%%%%%%%%%%%
\section{M-matrix property and convergence analysis}\label{error:analysis:intersect:sec}

An M-matrix is a non-singular matrix with non-positive off-diagonal entries and positive diagonal entries such that all row sums are non-negative with at least one row sum being positive.
In this section, we consider only the compact  9-point schemes resulting in a linear system with an M-matrix. It allows us to prove their convergence rate. For the sake of readability, the proofs of the results stated in this section are provided in \cref{Intersect:sec:proofs}.

For a scheme with stencil coefficients $\{C_{k,\ell}\}_{k,\ell=-1,0,1}$, let us introduce
the following sign condition:
\be\label{intersect:sign:condition}
\begin{cases}
	C_{k,\ell}>0, &\quad \mbox{if} \quad (k,\ell)=(0,0),\\
	C_{k,\ell}\le  0, &\quad \mbox{if} \quad (k,\ell)\ne(0,0),
\end{cases}
\ee
and summation condition:
\be\label{intersect:sum:condition}
\sum_{k=-1}^1\sum_{\ell=-1}^1 C_{k,\ell}=0.
\ee
 Under suitable boundary conditions (such as Dirichlet boundary conditions), it is well known that the sign condition \eqref{intersect:sign:condition} and summation condition \eqref{intersect:sum:condition} together guarantee the resulting coefficient matrix to be an M-matrix \cite{LiZhang2020}.

\subsection{Case 1: The interfaces intersection point $(\xi,\zeta)$ is a grid point}\label{convergence:section1}
\begin{theorem}\label{intersect:thm:convergence}
Consider the elliptic cross-interface problem  in \eqref{intersect:3} with the  interfaces intersection point $(\xi,\zeta)$ being a grid point.   Assume that $u_{p}:=u \chi_{\Omega_{p}}$ and $f_{p}:=f \chi_{\Omega_{p}}$ have uniformly continuous partial derivatives of (total) orders up to seven and five, respectively, in each $\Omega_p$ for $p=1,2,3,4$. Also assume that
	the essentially one-dimensional functions $\phi_p$ and $\psi_p$  on the interface $\Gamma_p$ have uniformly continuous derivatives of orders up to seven and six, respectively, for $p=1,2,3,4$.
Then the matrix of the linear system resulting from the compact 9-point scheme given by \cref{Intersect:thm:regular,Intersect:thm:side1,Intersect:thm:cross1}, including the corresponding modifications for the other parts of the interface $\Gamma$, is an M-matrix.
Consequently, under the above assumptions,
the scheme
	is sixth-order accurate, i.e. there exists a positive constant $C$, independent of $h$, such that
	\be\label{intersect:order6:formula}
	\|u-u_h\|_\infty \le C h^6,
	\ee
	where $u$ is the exact solution of \eqref{intersect:3}, and $u_h$ is its numerical approximation.
\end{theorem}

\subsection{Case 2: The  interfaces intersection point $(\xi,\zeta)$ is not a grid point}\label{convergence:section2}

We first consider
the compact 9-point  scheme, with a consistency order four at a grid point $(x_{i},y_{j})\not\in \Gamma$, in \cref{theorem:side3:w}, where
	$(x_i^*,y_j^*):=(x_i-wh,y_j)\in \Gamma_1$ with $0<w<1$  (see the first panel of \cref{fig:intersect:side1:w}).
The summation condition \eqref{intersect:sum:condition} for $\{C_{k,\ell}\}_{k,\ell=-1,0,1}$ in \eqref{Intersect:CKL:side3:w} can be directly verified.
	From \eqref{r1tor12:w:plus}, we can check that $r_{p}>0$ for all $p=1,2,\dots,12$ and $w\in(0,1)$, while $s_{p}<0$ for all $p=1,2,3$ and $w\in(0,1)$.
	Therefore, the stencil coefficients
	$C_{k,\ell}, k,\ell=-1,0,1$ in \eqref{Intersect:CKL:side3:w} satisfy:
	\[
	\begin{cases}
		C_{k,\ell}>0, &\quad \mbox{if} \quad (k,\ell)=(0,0),\\
		C_{k,\ell}< 0, &\quad \mbox{if} \quad (k,\ell)\ne\{(0,0), (1,-1),(1,1)\},
	\end{cases}
	\]
for any positive $a_1,a_2$, and $w\in(0,1)$.
Since we know from \cref{theorem:side3:w} that all the coefficients $C_{k,\ell},k,\ell=-1,0,1$ in \eqref{Intersect:CKL:side3:w} are unique after normalization, we observe from the above inequalities that the FDM in \cref{theorem:side3:w} satisfies the sign condition in \eqref{intersect:sign:condition} if and only if
\be \label{CCC}
C_{1,-1}\le 0 \quad \mbox{ and }\quad C_{1,1}\le 0.
\ee
%%%%%%
By \eqref{Intersect:CKL:side3:w} and \eqref{r1tor12:w:plus}, we have
\[
C_{1,-1}=C_{1,1}=\frac{-(r_6
+r_7\alpha+t_1\alpha^2)}{r_1 +r_2\alpha^2+r_3\alpha}, \quad \alpha=\frac{a_1}{a_2}>0,
\]
\[
 t_1=  -4w^5+12w^4-13w^3+8w^2-w,
\]
and $r_{p}>0$ for all $p=1,2,3,6,7$ and $w\in (0,1)$.
We can easily check that $t_1 > 0$ for $w\in[0.162,1)$ and $t_1 < 0$ for $w\in(0,0.161]$. So  $w \in[0.162,1)$ is required to achieve \eqref{CCC} for all positive $a_1, a_2$.
When $(x_i-h,y_j)$  (see the left panel in \cref{fig:intersect:side1:w}) or $(x_i+h,y_j)$  (see the right panel in \cref{fig:intersect:side1:w}) is the center point,    $(1-w) \in[0.162,1)$ is also necessary.
Thus the conditions in \eqref{CCC} only hold for all positive $a_1, a_2$ if $w\in[0.162,0.838]$.
However, for $w\in (0,1)$ outside $[0.162,0.838]$, one can always find particular positive $a_1,a_2$ so that \eqref{CCC} fails.
In order to achieve the sign condition for proving convergence, it is necessary to lower the consistency order.

Now we propose the following compact  scheme with a consistency order three such that $\{C_{k,\ell}\}_{k,\ell=-1,0,1}$  satisfies the sign condition \eqref{intersect:sign:condition} and the summation condition \eqref{intersect:sum:condition} for all positive $a_1,a_2$ and $w\in(0,1)$.

\begin{theorem}\label{theorem:side1:w:order3}
	Consider a grid point $(x_{i},y_{j})\not\in \Gamma$ such that
$(x_i^*,y_j^*):=(x_i-wh,y_j)\in \Gamma_1$ with $0<w<1$, $(\xi,\zeta)\not\in (x_i,y_j)+(-h,h)^2$ (see the first panel of \cref{fig:intersect:side1:w}).  Assume that $u_{p}:=u \chi_{\Omega_{p}}$ and $f_{p}:=f \chi_{\Omega_{p}}$ have uniformly continuous partial derivatives of (total) orders up to three and one, respectively, in each $\Omega_p$ for $p=1,2$. Also assume that
	the essentially one-dimensional functions $\phi_1$ and $\psi_1$  on the interface $\Gamma_1$ have uniformly continuous derivatives of orders up to three and two, respectively.
Let $\mathcal{L}_h$ be the discrete operator in \eqref{Luh} with the stencil coefficients
	{\footnotesize{
		\be\label{max:sign:Ckl:side}
		\begin{aligned}	
			&
			C_{0,0}:=1,\quad  C_{-1,1}=C_{-1,-1}:=
			 [(r_4\alpha^2+r_5\alpha)\rho+
			 t_1\alpha^2+t_2\alpha]/\beta,\\
			& C_{1,1} =C_{1,-1} :=\rho, \quad   C_{0,1}=C_{0,-1}:=
			 [(s_1+s_2\alpha^2+s_3\alpha)
			 \rho+t_3+t_4\alpha^2+t_5\alpha]/\beta,\\
			 %%%%%%%%%%%%%%%%%%%%%%%%%%%%%%%%%%%%%%%%%%%%%%%%%%%%%
			& C_{-1,0} :=[(s_4\alpha^2+s_5\alpha)\rho+
			 t_6\alpha^2+t_7\alpha]/\beta, \quad	 C_{1,0}:=[(s_6+s_7\alpha^2+s_8
			 \alpha)\rho+t_8+t_9\alpha^2+t_{10}\alpha]/\beta,
		\end{aligned}
		\ee
}}
where $\alpha:=a_1/a_2>0$, $\beta:=    r_1 +r_2\alpha^2+r_3\alpha> 0$,
$\rho \in \R$ is a free parameter, and
{\tiny{
		\be\label{Ckl:side:splus}
		\begin{aligned}	
			&r_1:=12(w^3-w^2-w+1),\
r_2:=4(w^3+3w^2+2w),\
r_3:=	4(-4w^3+w+3),\  r_4:=4(2w^3+w+3),\  	 r_5:=4(-2w^3-w+3),  \\
%%%%%%%%%%%%%%%%%%%%%%%%%%%%%%%%%%%%%%%%%%
&s_1:= 6(-w^2+2w-1),\
s_2:=-6(w^2+2w+1),\
s_3:=12(w^2-1),\  s_4:=4(-4w^3+w-3),\   s_5:=4(4w^3-w-3),  \\
&  s_6:=12(-2w^3+3w^2-1),\   s_7:=4(-2w^3-3w^2-w),\   s_8:=4(8w^3-6w^2+w-3),\   t_1:=-4w^3+6w^2-2w,  \\
& t_2:=4w^3-6w^2+2w,\   t_3:=-6w^3+9w^2-3,\   t_4:=-2w^3-3w^2-w,\   t_5:=8w^3-6w^2+w-3,\   	 t_6:=8w^3-12w^2-2w, \\
& t_7:=-8w^3+12w^2+2w-6,\    t_8:= 6(-w^2+2w-1),\     t_9:=-6w^2,\     t_{10}:=12(w^2-w).
		\end{aligned}
		\ee	}}
Then the compact 9-point finite difference scheme $h^{-1}\mathcal{L}_h u_h=h^{-1}F$
has a third order of consistency  at the grid point $(x_i,y_j)\not \in \Gamma$ with $(\xi,\zeta)\not\in (x_i,y_j)+(-h,h)^2$,
where
	\be\label{LhGamma1uh:order3}
	\begin{aligned}
F:=	&\frac{1}{a_1}\sum_{(m,n)\in \ind_{1}} f_1^{(m,n)}H^{-,w}_{3,m,n}+\frac{1}{a_2}
\sum_{(m,n)\in \ind_{1}} f_2^{(m,n)}
H^{+,w}_{3,m,n}-\sum_{n=0}^{3} \phi_1^{(n)}  G^w_{3,0,n}-\frac{1}{a_1}\sum_{n=0}^{2} \psi_1^{(n)}  G^w_{3,1,n}
	\end{aligned}
	\ee
and $H^{-,w}_{3,m,n}, H^{+,w}_{3,m,n}, G^w_{3,0,n}, G^w_{3,1,n}$ are defined in \eqref{GHw:v}.
Furthermore, for all positive $a_1,a_2$ and $w\in(0,1)$, $\{C_{k,\ell}\}_{k,\ell=-1,0,1}$ in \eqref{max:sign:Ckl:side} always satisfies
the summation condition \eqref{intersect:sum:condition} for all $\rho\in \R$, while $\{C_{k,\ell}\}_{k,\ell=-1,0,1}$ in \eqref{max:sign:Ckl:side} satisfies the sign condition \eqref{intersect:sign:condition} if and only if $\rho$ belongs to the following nonempty interval:
\be \label{rho:range}
\left[\max\left( -\frac{t_3+t_4\alpha^2+t_5\alpha}{s_1+s_2\alpha^2+s_3\alpha}, \  -\frac{t_6\alpha^2+t_7\alpha}{s_4\alpha^2+s_5\alpha}, \  -\frac{t_8+t_9\alpha^2+t_{10}\alpha}{s_6+s_7\alpha^2+s_8\alpha},  \right),\
\min\left(0,   \ -\frac{t_1\alpha^2+t_2\alpha}{r_4\alpha^2+r_5\alpha} \right)\right].
\ee
\end{theorem}

For the range of the parameter $\rho$ to achieve the M-matrix property,
in fact, the nonempty interval in \eqref{rho:range} always contains the following subintervals:
\be\label{rho:interval}
\begin{cases}
[-0.2, -0.018], &\text{if } \alpha \in (0,1] \text{ and } w\in(0,1/2],\\
[-0.04, 0], &\text{if } \alpha \in (0,1] \text{ and } w\in[1/2,1),\\	
0,  &\text{if } \alpha \in [1,+\infty) \text{ and } w\in(0,1/2],\\
[-0.1, -0.012], &\text{if } \alpha \in [1,+\infty) \text{ and } w\in[1/2,1).
\end{cases}
\ee

In order to produce an overall scheme satisfying the M-matrix property, it is also necessary to lower the order of consistency  for grid points near the interfaces intersection point $(\xi,\zeta)$.
Such a scheme is given in the next theorem.

\begin{theorem}\label{Intersect:thm:cross1:w1w2:order2}
Consider a grid point $(x_{i},y_{j})\not\in \Gamma$ such that $(x_{i}^*,y_{j}^*):=(\xi,\zeta)=(x_i-w_1h,y_j-w_2h)$
with $0<w_1,w_2<1$ (see the first panel of \cref{fig:intersect:cross:point:w}).  Assume that $u_{p}:=u \chi_{\Omega_{p}}$ and $f_{p}:=f \chi_{\Omega_{p}}$ have uniformly continuous partial derivatives of (total) orders up to two and zero, respectively, in each $\Omega_p$ for $p=1,2,3,4$. Also assume that
	the essentially one-dimensional functions $\phi_p$ and $\psi_p$  on the interface $\Gamma_p$ have uniformly continuous derivatives of orders up to two and one, respectively, for $p=1,2,3,4$.
Let $\mathcal{L}_h$ be the discrete operator in \eqref{Luh} with the stencil coefficients
	{\small{		 \be\label{max:sign:Ckl:intersect}
			\begin{aligned}	
				&
C_{-1,0}:=-[a_{1}a_{2}r_{2}+a_{1}a_{3}r_{1}]/\beta,\quad  C_{0,-1}:=-[a_{1}a_{3}r_{4}+a_{2}a_{3}r_{3}]/\beta,\\ &C_{0,1}:=-[a_{1}a_{2}r_{8}+a_{1}a_{3}r_{7}+a_{2}^2r_{6}+a_{2}a_{3}r_{5}]/\beta, \quad
C_{1,0}:=-[a_{1}a_{2}r_{12}+a_{1}a_{3}r_{11}+a_{2}^2r_{10}+a_{2}a_{3}r_{9}]/\beta,\\
&C_{0,0}:=1,\quad C_{-1,-1}=C_{-1,1}=C_{1,-1}=C_{1,1}:=0,
			\end{aligned}
			\ee
}}
where $\beta:=-( a_{1}a_{2}s_{4}+a_{1}a_{3}s_{3}+a_{2}^2s_{2}+a_{2}a_{3}s_{1})> 0$ and
{\tiny{
\be\label{Ckl:cross:tplus}
\begin{aligned}	
			&r_{1} := 2w_{2}^2-w_{2}+1,\   r_{2} := -2w_{2}^2+w_{2}+1,\   r_{3} := -2w_{1}^2+w_{1}+1,\   r_{4} := 2w_{1}^2-w_{1}+1,\   r_{5} := -w_{2}(2w_{1}^2-w_{1}-1), \\
&  r_{6} := (w_{2}-1)(2w_{1}^2-w_{1}-1),\   r_{7} := w_{2}(2w_{1}^2-w_{1}+1),\   r_{8} := -(w_{2}-1)(2w_{1}^2-w_{1}+1),\   r_{9} := -(2w_{2}^2-w_{2}+1)(w_{1}-1), \\
&  r_{10} := (2w_{2}^2-w_{2}-1)(w_{1}-1),\    r_{11} := (2w_{2}^2-w_{2}+1)w_{1},\   r_{12} := -(2w_{2}^2-w_{2}-1)w_{1},\    s_{1} := 2(w_{1}-1)(w_{1}w_{2}+w_{2}^2+w_{1}+1), \\
&  s_{2} := 2(1-w_{2})(w_{1}-1)(w_{1}+w_{2}+1),\    s_{3} := -2(w_{2}+1)w_{1}^2-2(w_{2}^2-w_{2})w_{1}-2w_{2}^2-2,\   s_{4} := 2(w_{2}-1)(w_{1}^2+w_{1}w_{2}+w_{2}+1).
		\end{aligned}
		\ee}}
Then the compact 9-point finite difference scheme $h^{-1}\mathcal{L}_h u_h=h^{-1}F$ has a second order of consistency  at the grid point $(x_i,y_j)\not \in \Gamma$ with $(\xi,\zeta)\in (x_i,y_j)+(-h,h)^2$,
where
\[
F:= \sum_{(n,m)\in \ind_{2}^{1}} F^{1,w_1,w_2}_{2,m,n}
	+ \sum_{(m,n)\in \ind_{0}} F^{w_1,w_2}_{2,m,n} %\\&
+ \sum_{(n,m)\in \ind_{2}^{1}} \Phi^{1,w_1,w_2}_{2,m,n}+\sum_{m=0}^{2}  \Phi^{w_1,w_2}_{2,m}
	+ \sum_{(n,m)\in \ind_{2}^{1}} \Psi^{1,w_1,w_2}_{2,m,n}+ \sum_{m=0}^{1}  \Psi^{w_1,w_2}_{2,m}
\]
and $F^{1,w_1,w_2}_{2,m,n}$, $F^{w_1,w_2}_{2,m,n}$, $\Phi^{1,w_1,w_2}_{2,m,n}$, $\Phi^{w_1,w_2}_{2,m}$, $\Psi^{1,w_1,w_2}_{2,m,n}$, $\Psi^{w_1,w_2}_{2,m}$ are defined in  \eqref{intersect:FH1mn:Cross}--\eqref{intersect:PsiH1Mmn:Cross}.
Furthermore, $\{C_{k,\ell}\}_{k,\ell=-1,0,1}$ in \eqref{max:sign:Ckl:intersect}  satisfy both, the sign condition \eqref{intersect:sign:condition},  and the summation condition \eqref{intersect:sum:condition} for all positive $a_1,a_2,a_3,a_4$ and $(w_1,w_2)\in(0,1)^2$.
\end{theorem}

The convergence rate of the overall scheme, that satisfies the M-matrix property, is claimed in the next theorem.
\begin{theorem}\label{intersect:thm:convergence:order3}
	Consider the elliptic cross-interface problem  in \eqref{intersect:3} and assume that the intersection point $(\xi,\zeta)$ is not a grid point.  Assume that $u_{p}:=u \chi_{\Omega_{p}}$ and $f_{p}:=f \chi_{\Omega_{p}}$ have uniformly continuous partial derivatives of (total) orders up to four and two, respectively, in each $\Omega_p$ for $p=1,2,3,4$. Also assume that the essentially one-dimensional functions $\phi_p$ and $\psi_p$  on the interface $\Gamma_p$ have uniformly continuous derivatives of orders up to four and three respectively for $p=1,2,3,4$.
Then the matrix of the linear system resulting from the compact 9-point scheme given by \cref{Intersect:thm:regular,theorem:side1:w:order3,Intersect:thm:cross1:w1w2:order2} with $\rho$ inside the nonempty interval in \eqref{rho:range}, including the corresponding modifications for the other parts of the interface $\Gamma$, is an M-matrix. Consequently, under above assumptions,
	the  scheme
	is third-order accurate, i.e., there exists a positive constant $C$, independent of $h$, such that:
\be\label{intersect:order6:formula:order3}
	\|u-u_h\|_\infty \le C h^3,
	\ee
	where $u$ is the exact solution of \eqref{intersect:3}, and $u_h$ is its numerical approximation.
	%%
	%where $u$ is the exact solution of \eqref{intersect:3}, and $u_h$ is the numerical solution  computed by  compact  FDM in \cref{Intersect:thm:regular,theorem:side1:w:order3,Intersect:thm:cross1:w1w2:order2} of \eqref{intersect:3}.
\end{theorem}

\section{Numerical experiments}\label{Intersect:sec:Numeri}
%%
%	
%
%Let $\Omega:=(l_1,l_2)\times(l_3,l_4)$ with
%$l_4-l_3=N_0(l_2-l_1)$ for some positive integer $N_0$. For a given $J\in \NN$, we define
%$h:=(l_2-l_1)/N_1$ with $N_1:=2^J$ and let
%$x_i:=l_1+ih$ and
%$y_j:=l_3+jh$ for $i=0,1,\dots,N_1$ and $j=0,1,\dots,N_2$ with $N_2:=N_0N_1$.
%Let
%$u$ be the exact solution of \eqref{intersect:3} and $(u_{h})_{i,j}$ be the numerical solution at the grid point $(x_i, y_j)$ using the mesh size $h$.

Let us now choose  $\Omega=(0,1)^2$ and $N_1=N_2:=2^J$ for some $J\in \NN$.
To quantify the accuracy of the various schemes in the numerical examples presented below we use the relative error in the $l_2$ norm:
$\epsilon_h := \frac{\|u_{h}-u\|_{2}}{\|u\|_{2}}$, where:
\begin{align*}
	\|u_{h}-u\|_{2}^2:= h^2&\sum_{i=0}^{N_1}\sum_{j=0}^{N_2} \left((u_h)_{i,j}-u(x_i,y_j)\right)^2, \qquad \ \|u\|_{2}^2:=h^2 \sum_{i=0}^{N_1}\sum_{j=0}^{N_2} \left(u(x_i,y_j)\right)^2,
\end{align*}
as well as its infinity  norm:
\[
\|u_h-u\|_\infty
:=\max_{0\le i\le N_1, 0\le j\le N_2} \left|(u_h)_{i,j}-u(x_i,y_j)\right|.
\]

The theoretical rate of convergence of the various schemes presented above is verified on a set of numerical examples.
\subsection{Compact  9-point  scheme in \cref{Intersect:thm:regular,Intersect:thm:side1,Intersect:thm:cross1} (if $(\xi,\zeta)$ is a grid point):
}\label{numerical:section1}
\begin{figure}[htbp]
	\centering
	\begin{subfigure}[b]{0.3\textwidth}
		\hspace{0.5cm}
		\begin{tikzpicture}[scale = 4.35]
			\draw	(0, 0) -- (0, 1) -- (1, 1) -- (1, 0) --(0,0);	
			\draw   (1/2, 0) -- (1/2, 1);	
			\draw   (0, 1/2) -- (1, 1/2);	
			 %%%%%%%%%%%%%%%%%%%%%%%%%%%%%%%%%%%%%%%%%%%%%%%%%
			\node (A) at (1/4,3/4) {$a_1=10^{-5}$};
			\node (A) at (3/4,3/4) {$a_2=10^{5}$};
			\node (A) at (3/4,1/4) {$a_3=10^{-5}$};
			\node (A) at (1/4,1/4) {$a_4=10^{5}$};
		\end{tikzpicture}
	\end{subfigure}
	\begin{subfigure}[b]{0.3\textwidth}
		 \includegraphics[width=5.5cm,height=4.5cm]{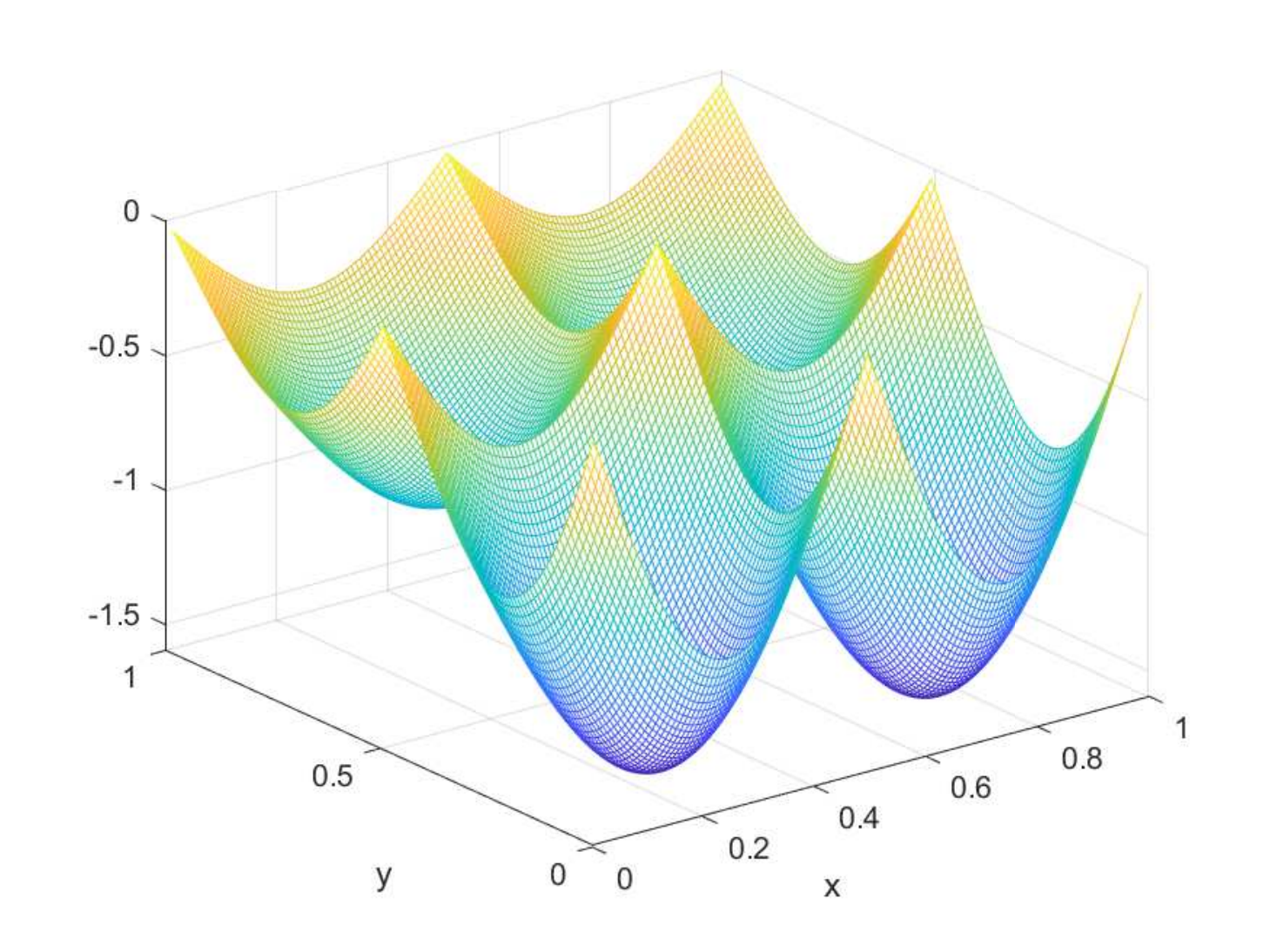}
	\end{subfigure}
	\begin{subfigure}[b]{0.3\textwidth}
		 \includegraphics[width=5.5cm,height=4.5cm]{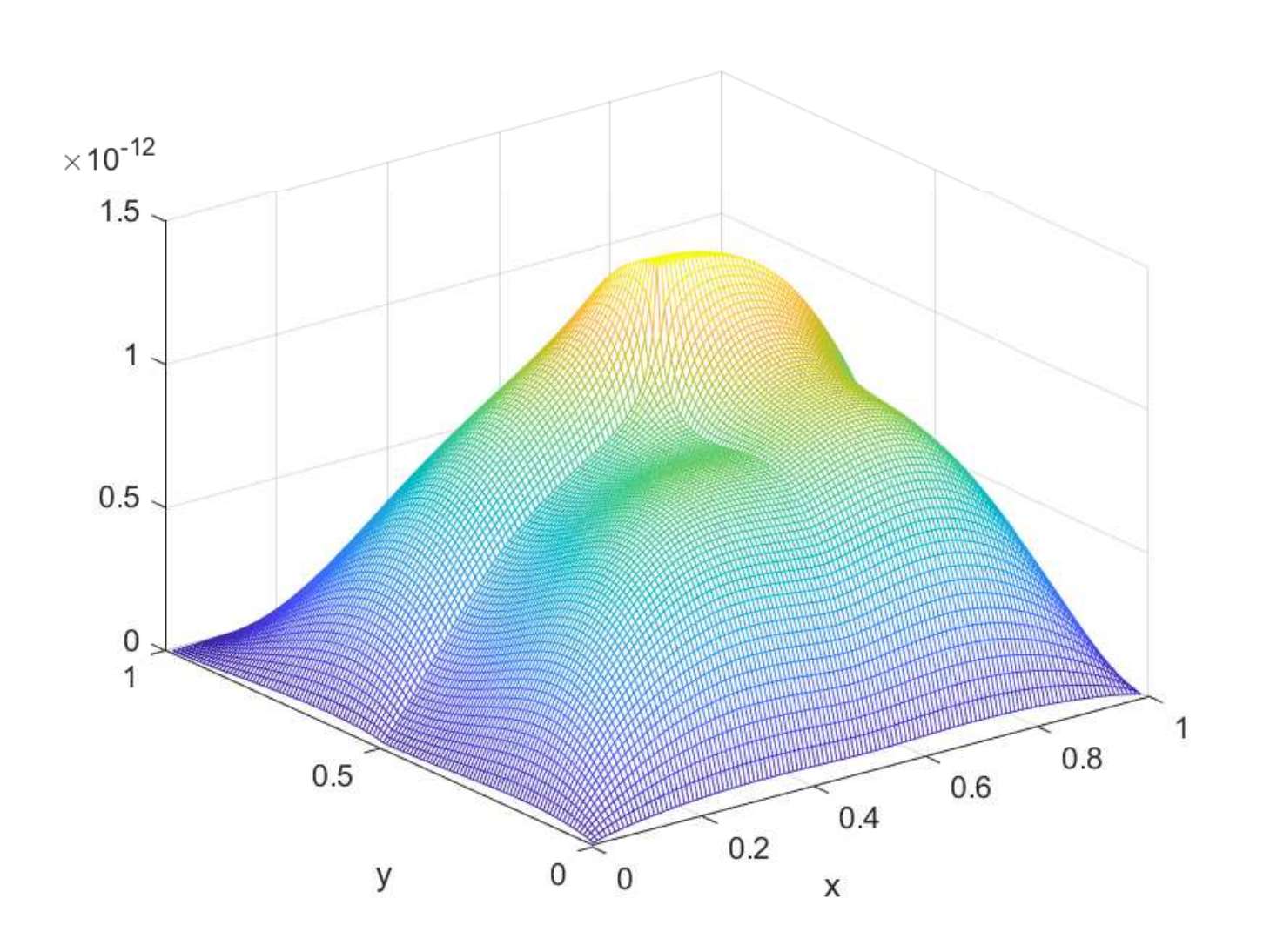}
	\end{subfigure}
	\caption
	{\cref{Intersect:ex1}:  the coefficient $a(x,y)$ (left),  the numerical solution $(u_h)_{i,j}$ (middle), and  $|(u_h)_{i,j}-u(x_i,y_j)|$ (right), at all grid points $(x_i,y_j)$ on $\Omega$ with $h=2^{-7}$,  for the scheme in \cref{Intersect:thm:regular,Intersect:thm:side1,Intersect:thm:cross1}.}
	\label{Intersect:Numeri:fig1}
\end{figure}

\begin{table}[htbp]
	\caption{Performance in \cref{Intersect:ex1,Intersect:ex2,Intersect:ex3}  of the scheme given in \cref{Intersect:thm:regular,Intersect:thm:side1,Intersect:thm:cross1} on uniform Cartesian meshes with $h=2^{-J}$. Note that $\epsilon_h := \frac{\|u_{h}-u\|_{2}}{\|u\|_{2}}$. }
	\centering
	\setlength{\tabcolsep}{0.3mm}{
		 \begin{tabular}{c|c|c|c|c|c|c|c|c|c|c|c|c}
			\hline
			\multicolumn{1}{c|}{}  &
			 \multicolumn{4}{c|}{\cref{Intersect:ex1}}  &
			 \multicolumn{4}{c|}{\cref{Intersect:ex2}}  &
			 \multicolumn{4}{c}{\cref{Intersect:ex3}}  \\
			\hline
			$J$
			&  $\epsilon_h $
			&order &  $\|u_{h}-u\|_{\infty}$ & order &  $\epsilon_h$
			&order &  $\|u_{h}-u\|_{\infty}$ & order &  $\epsilon_h $
			&order &  $\|u_{h}-u\|_{\infty}$ & order \\
			\hline
2   &1.046E-03   &   &1.537E-03   &   &1.236E-05   &   &2.416E-06   &   &   &   &   &\\
3   &1.303E-05   &6.3   &2.225E-05   &6.1   &2.093E-07   &5.9   &3.930E-08   &5.9   &   &   &   &\\
4   &1.887E-07   &6.1   &3.396E-07   &6.0   &3.488E-09   &5.9   &6.340E-10   &6.0   &2.185E-02   &   &2.185E+02   &\\
5   &2.852E-09   &6.0   &5.268E-09   &6.0   &5.673E-11   &5.9   &1.027E-11   &5.9   &2.761E-04   &6.3   &2.761E+00   &6.3\\
6   &4.380E-11   &6.0   &8.216E-11   &6.0   &9.043E-13   &6.0   &1.598E-13   &6.0   &4.079E-06   &6.1   &4.079E-02   &6.1\\
7   &7.888E-13   &5.8   &1.424E-12   &5.9   &   &   &   &   &6.283E-08   &6.0   &6.283E-04   &6.0\\
8   &   &   &   &   &   &   &   &   &9.783E-10   &6.0   &9.783E-06   &6.0\\		
			\hline
	\end{tabular}}
	\label{Intersect:Numeri:table1}
\end{table}	
%
%
%%
%%%%%%%%%%%%%%%%%%%%%%%%%%%%%%%%%%%%%%%%%%%%%
\begin{example}\label{Intersect:ex1}
	\normalfont
	Let $\Omega=(0,1)^2$ and $(\xi,\zeta)=(1/2,1/2)$.
	The functions  in \eqref{intersect:3} are given by
	%%%
	\begin{align*}
		&a_{1}=a \chi_{\Omega_{1}}=10^{-5},\qquad  a_{2}=a \chi_{\Omega_{2}}=10^{5}, \qquad a_{3}=a \chi_{\Omega_{3}}=10^{-5},\qquad  a_{4}=a \chi_{\Omega_{4}}=10^{5},\\
		& {  u_{1}=u \chi_{\Omega_{1}}=-\sin(2\pi x)\exp(-y)-\sin(2\pi(-y+1))\exp(-y),}\\
		& {  u_{2}=u \chi_{\Omega_{2}}=-\sin(2\pi (-x+1))\exp(-y)-\sin(2\pi(-y+1))\exp(-y),}\\
		& {  u_{3}=u \chi_{\Omega_{3}}=-\sin(2\pi (-x+1))\exp(-y)-\sin(2\pi y)\exp(-y),}\\
		& {  u_{4}=u \chi_{\Omega_{4}}=-\sin(2\pi x)\exp(-y)-\sin(2\pi y)\exp(-y),}
	\end{align*}
and	$g, f_{p}$, $\phi_p$, $\psi_p$ for $p=1,2,3,4$ in \eqref{intersect:3} are obtained by plugging the above functions into \eqref{intersect:3}.
	Note that $\phi_p=0$ for $p=1,2,3,4$.
	The numerical results are presented in \cref{Intersect:Numeri:table1} and \cref{Intersect:Numeri:fig1}.	
\end{example}

%%%%%%%%%%%%%%%%%%%%%%%%%%%%%%%%%%%%%%%%%%%%
%%%%%%%%%%%%%%%%%%%%%%%%%%%%%%%%%%%%%%%%%%%%%
%%%%%%%%%%%%%%%%%%%%%%%%%%%%%%%%%%%%%%%%%%%%%
%%%%%%%%%%%%%%%%%%%%%%%%%%%%%%%%%%%%%%%%%%%%%

\begin{figure}[htbp]
	\centering
	\begin{subfigure}[b]{0.3\textwidth}
		\hspace{0.5cm}
		\begin{tikzpicture}[scale = 4.35]
	\draw	(0, 0) -- (0, 1) -- (1, 1) -- (1, 0) --(0,0);	
	\draw   (1/2, 0) -- (1/2, 1);	
	\draw   (0, 1/2) -- (1, 1/2);	
	 %%%%%%%%%%%%%%%%%%%%%%%%%%%%%%%%%%%%%%%%%%%%%%%%%
	\node (A) at (1/4,3/4) {$a_1=10^{7}$};
	\node (A) at (3/4,3/4) {$a_2=10^{-3}$};
	\node (A) at (3/4,1/4) {$a_3=10^{4}$};
	\node (A) at (1/4,1/4) {$a_4=10^{-6}$};
\end{tikzpicture}
	\end{subfigure}
	\begin{subfigure}[b]{0.3\textwidth}
		 \includegraphics[width=5.5cm,height=4.5cm]{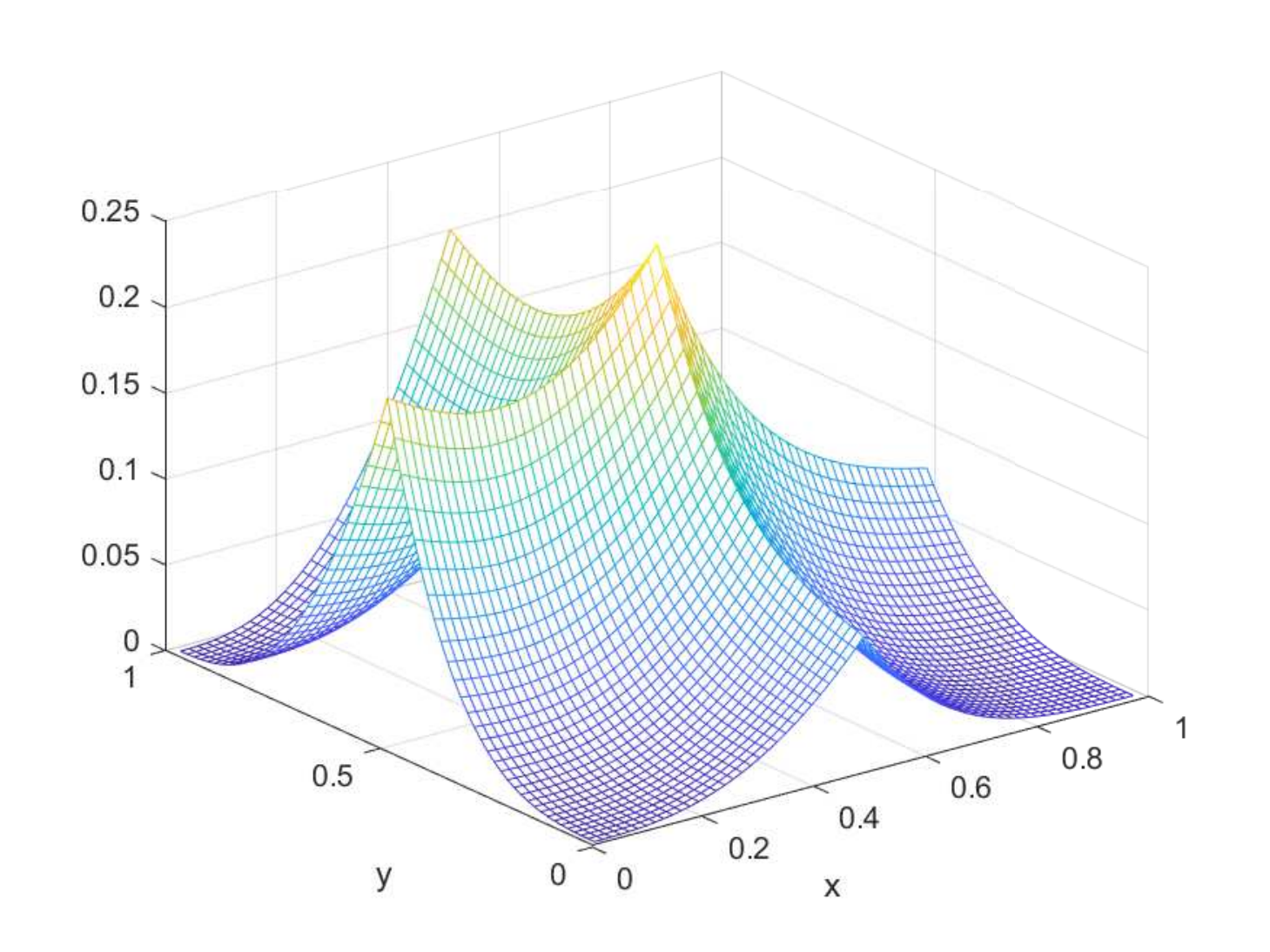}
	\end{subfigure}
	\begin{subfigure}[b]{0.3\textwidth}
		 \includegraphics[width=5.5cm,height=4.5cm]{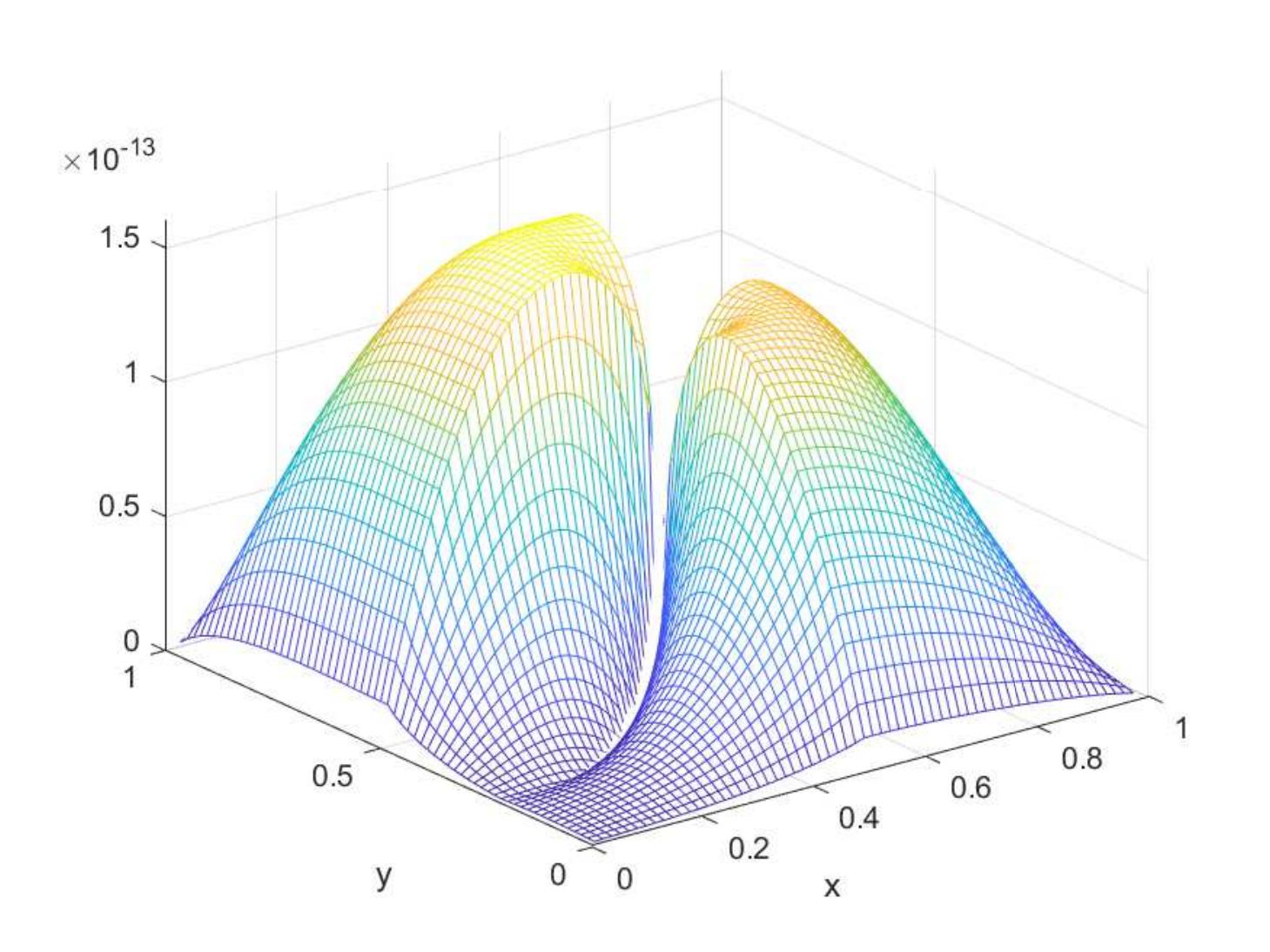}
	\end{subfigure}
	\caption
	{\cref{Intersect:ex2}:  the coefficient $a(x,y)$ (left),  the numerical solution $(u_h)_{i,j}$ (middle), and  $|(u_h)_{i,j}-u(x_i,y_j)|$ (right), at all grid points $(x_i,y_j)$ on $\Omega$ with $h=2^{-6}$,  for the scheme in \cref{Intersect:thm:regular,Intersect:thm:side1,Intersect:thm:cross1}.}
	\label{Intersect:Numeri:fig2}
\end{figure}

\begin{example}\label{Intersect:ex2}
	\normalfont
	Let $\Omega=(0,1)^2$  and $(\xi,\zeta)=(1/2,1/2)$.
	The functions  in \eqref{intersect:3} are given by
	%%%
	\begin{align*}
		&a_{1}=a \chi_{\Omega_{1}}=10^{7},\qquad  a_{2}=a \chi_{\Omega_{2}}=10^{-3}, \qquad a_{3}=a \chi_{\Omega_{3}}=10^{4},\qquad  a_{4}=a \chi_{\Omega_{4}}=10^{-6},\\
		&  {  u_{1}=u \chi_{\Omega_{1}}=(x^3+(1-y)^3)\exp(-x+y),\qquad   u_{2}=u \chi_{\Omega_{2}}=((1-x)^3+(1-y)^3)\exp(-x+y),}\\
		&  {  u_{3}=u \chi_{\Omega_{3}}=((1-x)^3+y^3)\exp(-x+y),\qquad  u_{4}=u \chi_{\Omega_{4}}=(x^3+y^3)\exp(-x+y),}
	\end{align*}
and	$g, f_{p}$, $\phi_p$, $\psi_p$ for $p=1,2,3,4$ in \eqref{intersect:3} are obtained by plugging the above functions into \eqref{intersect:3}. Again, $\phi_p=0$ for $p=1,2,3,4$.
	The numerical results are presented in \cref{Intersect:Numeri:table1} and \cref{Intersect:Numeri:fig2}.	
\end{example}

%%%%%%%%%%%%%%%%%%%%%%%%%%%%%%%%%%%%%%%%%%%%
%%%%%%%%%%%%%%%%%%%%%%%%%%%%%%%%%%%%%%%%%%%%%
%%%%%%%%%%%%%%%%%%%%%%%%%%%%%%%%%%%%%%%%%%%%%
%%%%%%%%%%%%%%%%%%%%%%%%%%%%%%%%%%%%%%%%%%%%%

\begin{figure}[htbp]
	\centering
	\begin{subfigure}[b]{0.3\textwidth}
		\hspace{0.5cm}
		\begin{tikzpicture}[scale = 4.35]
	\draw	(0, 0) -- (0, 1) -- (1, 1) -- (1, 0) --(0,0);	
	\draw   (1/4, 0) -- (1/4, 1);	
	\draw   (0, 1/8) -- (1, 1/8);	
	 %%%%%%%%%%%%%%%%%%%%%%%%%%%%%%%%%%%%%%%%%%%%%%%%%
	\node (A)[scale = 0.8] at (0.13,0.5) {$10^{-4}$};
	\node (A) at (0.7,0.5) {$a_2=10^{5}$};
	\node (A)[scale = 0.7] at (0.7,0.08) {$a_3=2\times 10^{-4}$};
	\node (A)[scale = 0.7] at (0.13,0.08) {$10^{6}$};
\end{tikzpicture}
	\end{subfigure}
	\begin{subfigure}[b]{0.3\textwidth}
		 \includegraphics[width=5.5cm,height=4.5cm]{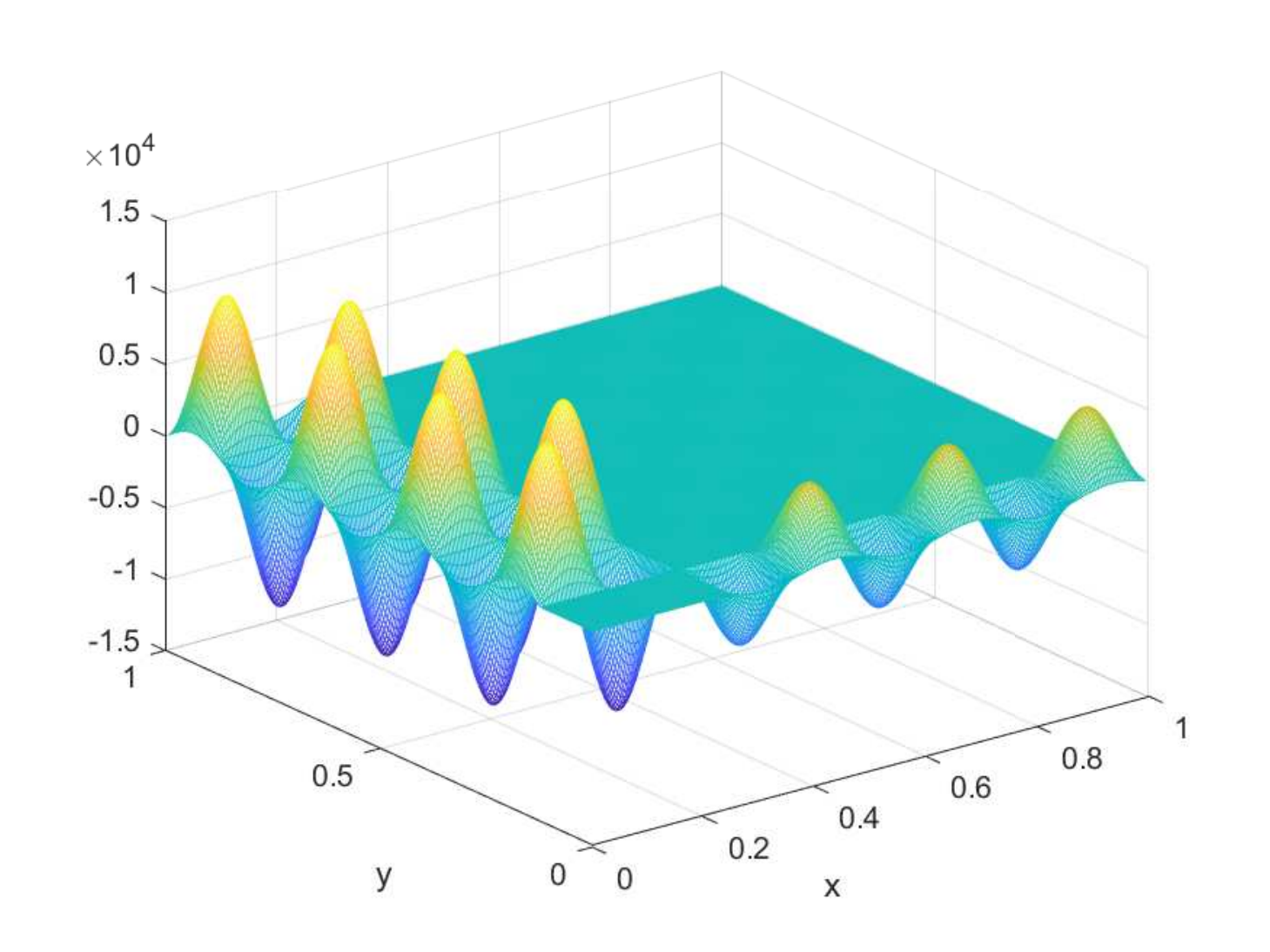}
	\end{subfigure}
	\begin{subfigure}[b]{0.3\textwidth}
		 \includegraphics[width=5.5cm,height=4.5cm]{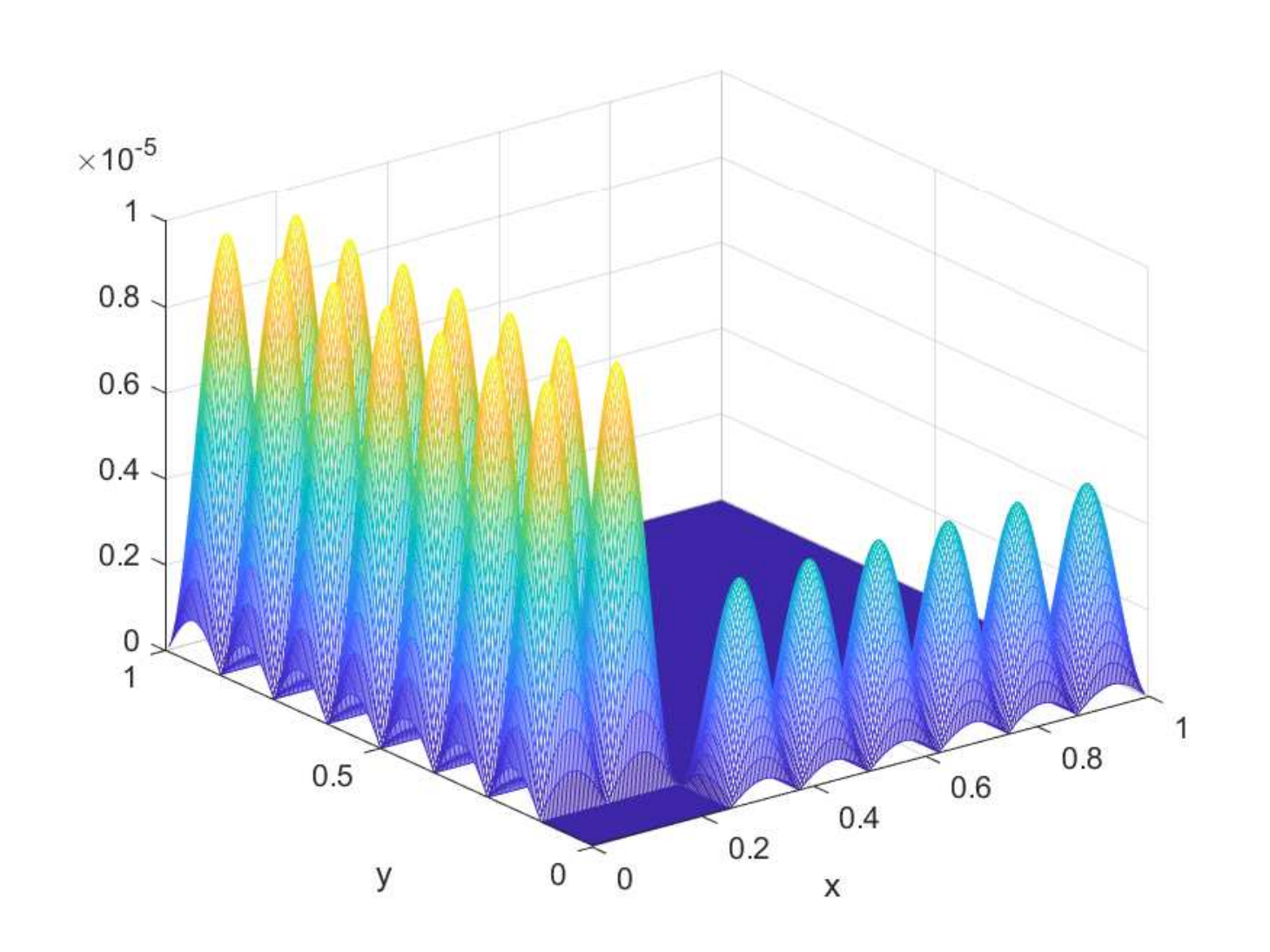}
	\end{subfigure}
	\caption
	{\cref{Intersect:ex3}: the coefficient $a(x,y)$ (left),  $-(u_h)_{i,j}$ (middle), and  $|(u_h)_{i,j}-u(x_i,y_j)|$ (right), at all grid points $(x_i,y_j)$ on $\Omega$ with $h=2^{-8}$,  for the scheme  in \cref{Intersect:thm:regular,Intersect:thm:side1,Intersect:thm:cross1}.}
	\label{Intersect:Numeri:fig3}
\end{figure}

\begin{example}\label{Intersect:ex3}
	\normalfont
	Let $\Omega=(0,1)^2$ and $(\xi,\zeta)=(1/4,1/8)$.
	The functions  in \eqref{intersect:3} are given by
	%%%
	\begin{align*}
		&a_{1}=a \chi_{\Omega_{1}}=10^{-4},\qquad  a_{2}=a \chi_{\Omega_{2}}=10^{5}, \qquad a_{3}=a \chi_{\Omega_{3}}=2\times 10^{-4},\qquad  a_{4}=a \chi_{\Omega_{4}}=10^{6},\\
		&  {  u_{1}=u \chi_{\Omega_{1}}=\sin(8\pi x)\sin(8\pi y)/a_1,\qquad u_{2}=u \chi_{\Omega_{2}}=\sin(8\pi x)\sin(8\pi y)/a_2,}\\
		&  {  u_{3}=u \chi_{\Omega_{3}}=\sin(8\pi x)\sin(8\pi y)/a_3,\qquad u_{4}=u \chi_{\Omega_{4}}=\sin(8\pi x)\sin(8\pi y)/a_4,}
	\end{align*}
and	functions $g, f_{p}$, $\phi_p$, $\psi_p$ for $p=1,2,3,4$ in \eqref{intersect:3} are obtained by plugging the above functions into \eqref{intersect:3}.
	Note that $\phi_p=0$ for $p=1,2,3,4$.
	The numerical results are presented in \cref{Intersect:Numeri:table1} and \cref{Intersect:Numeri:fig3}.	
\end{example}
%%%%%%%%%%%%%%%%%%%%%%%%%%%%%%%%%%%%%%%%%%%%
%%%%%%%%%%%%%%%%%%%%%%%%%%%%%%%%%%%%%%%%%%%%%
%%%%%%%%%%%%%%%%%%%%%%%%%%%%%%%%%%%%%%%%%%%%%	
%
%
%
%
%
\subsection{Compact  9-point  scheme given  in \cref{Intersect:thm:regular,theorem:side3:w,Intersect:thm:cross1:w1w2}\label{numerical:section2} (if $(\xi,\zeta)$ is not a grid point):}

\begin{table}[htbp]
	\caption{Performance in \cref{Intersect:ex4,Intersect:ex5,Intersect:ex6}  of the scheme in \cref{Intersect:thm:regular,theorem:side3:w,Intersect:thm:cross1:w1w2} on uniform Cartesian meshes with $h=2^{-J}$.  Note that $\epsilon_h := \frac{\|u_{h}-u\|_{2}}{\|u\|_{2}}$.}
	\centering
	\setlength{\tabcolsep}{0.1mm}{
		 \begin{tabular}{c|c|c|c|c|c|c|c|c|c|c|c|c}
			\hline
			\multicolumn{1}{c|}{}  &
			 \multicolumn{4}{c|}{\cref{Intersect:ex4}}  &
			 \multicolumn{4}{c|}{\cref{Intersect:ex5}}  &
			 \multicolumn{4}{c}{\cref{Intersect:ex6}}  \\
			\hline
			$J$
			&  $\epsilon_h$
			&order &  $\|u_{h}-u\|_{\infty}$ & order &  $\epsilon_h$
			&order &  $\|u_{h}-u\|_{\infty}$ & order &  $\epsilon_h$
			&order &  $\|u_{h}-u\|_{\infty}$ & order \\
			\hline
3   &   &   &  &   &   &   &   &   &2.916E-03   &   &7.344E+01   &\\
4   &1.236E-02   &   &4.022E-02   &   &5.236E-01   &   &9.973E-01   &   &6.873E-05   &5.4   &1.943E+00   &5.2\\
5   &4.082E-04   &4.9   &1.317E-03   &4.9   &2.057E-02   &4.7   &3.724E-02   &4.7   &1.605E-06   &5.4   &5.972E-02   &5.0\\
6   &9.528E-06   &5.4   &3.036E-05   &5.4   &1.261E-03   &4.0   &2.335E-03   &4.0   &1.199E-07   &3.7   &3.844E-03   &4.0\\
7   &3.116E-07   &4.9   &9.905E-07   &4.9   &5.273E-05   &4.6   &1.006E-04   &4.5   &3.058E-09   &5.3   &1.242E-04   &5.0\\
8   &1.293E-08   &4.6   &3.660E-08   &4.8   &1.674E-06   &5.0   &3.207E-06   &5.0   &6.623E-11   &5.5   &4.067E-06   &4.9\\			
			\hline
	\end{tabular}}
	\label{Intersect:Numeri:table2}
\end{table}

\begin{figure}[htbp]
	\centering
	\begin{subfigure}[b]{0.3\textwidth}
		\hspace{0.5cm}
		\begin{tikzpicture}[scale = 4.35]
			\draw	(0, 0) -- (0, 1) -- (1, 1) -- (1, 0) --(0,0);	
			\draw   (pi/5, 0) -- (pi/5, 1);	
			\draw   (0, pi/8) -- (1, pi/8);	
			 %%%%%%%%%%%%%%%%%%%%%%%%%%%%%%%%%%%%%%%%%%%%%%%%%
			\node (A)[scale = 0.9] at (0.3,0.7) {$a_1=10^{4}$};
			\node (A)[scale = 0.9] at (0.81,0.7) {$a_2=10^{-4}$};
			\node (A)[scale = 0.9] at (0.81,0.2) {$a_3=10^{4}$};
			\node (A)[scale = 0.9] at (0.3,0.2) {$a_4=10^{-4}$};
		\end{tikzpicture}
	\end{subfigure}
	\begin{subfigure}[b]{0.3\textwidth}
		 \includegraphics[width=5.5cm,height=4.5cm]{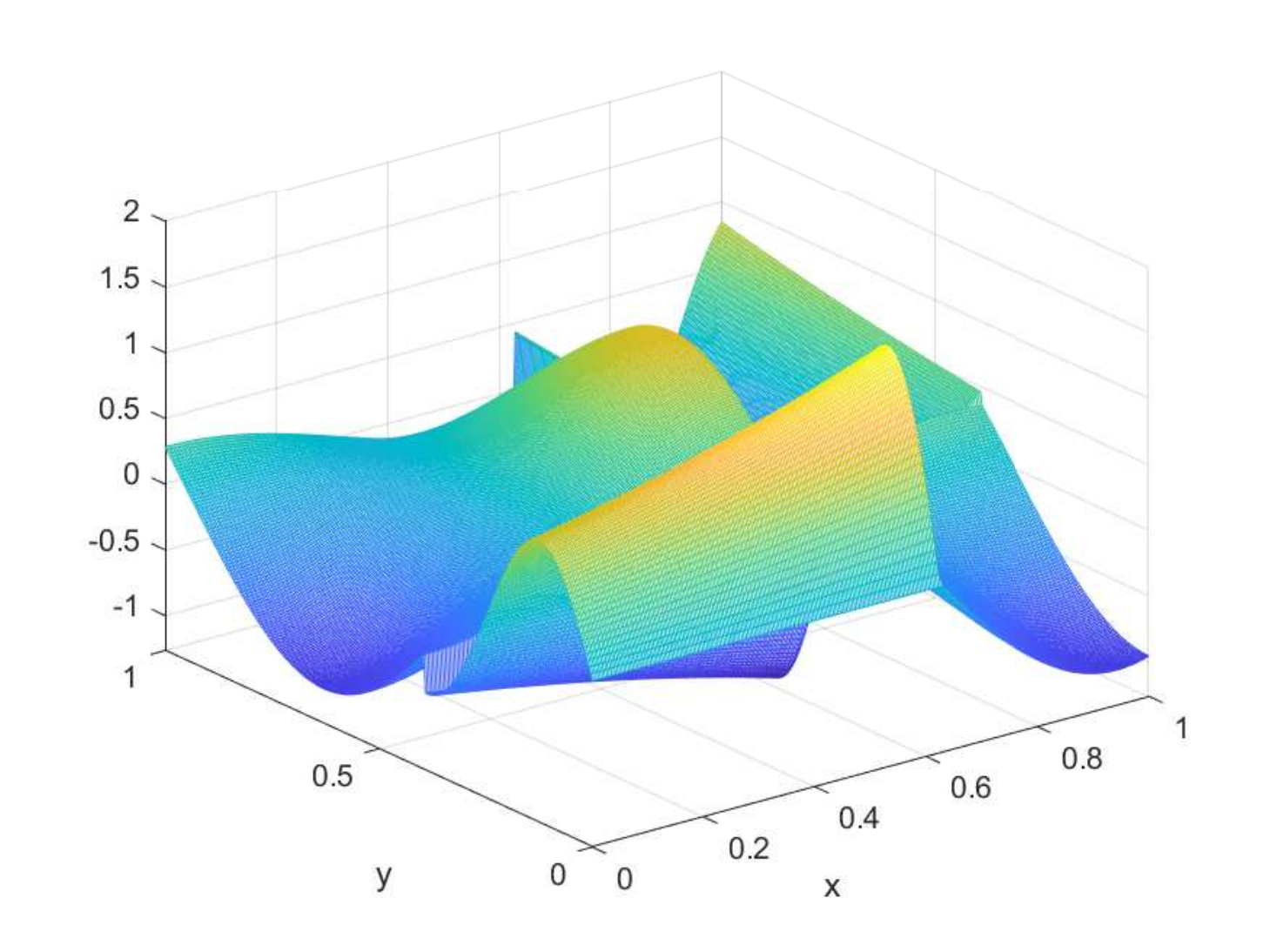}
	\end{subfigure}
	\begin{subfigure}[b]{0.3\textwidth}
		 \includegraphics[width=5.5cm,height=4.5cm]{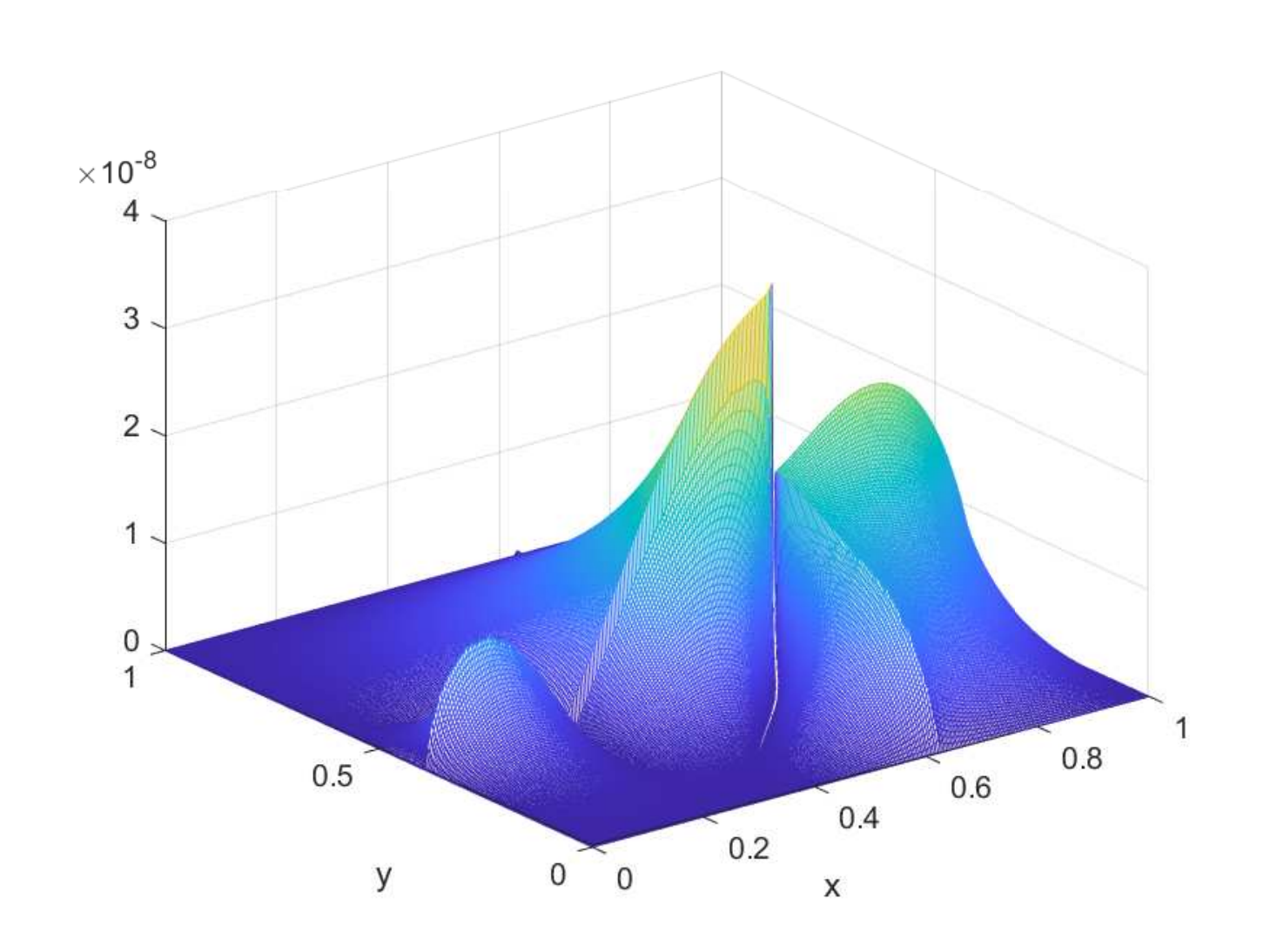}
	\end{subfigure}
	\caption
	{\cref{Intersect:ex4}:  the coefficient $a(x,y)$ (left),  the numerical solution $(u_h)_{i,j}$ (middle), and $|(u_h)_{i,j}-u(x_i,y_j)|$ (right), at all grid points $(x_i,y_j)$ on $\Omega$ with $h=2^{-8}$,  for the scheme  in \cref{Intersect:thm:regular,theorem:side3:w,Intersect:thm:cross1:w1w2}.}
	\label{Intersect:Numeri:fig4}
\end{figure}

%%%%%%%%%%%%%%%%%%%%%%%%%%%%%%%%%%%%%%%%%%%%%
\begin{example}\label{Intersect:ex4}
	\normalfont
	Let $\Omega=(0,1)^2$ and $(\xi,\zeta)=(\pi/5,\pi/8)$.
	The functions  in \eqref{intersect:3} are given by
	%%%
	\begin{align*}
		&a_{1}=a \chi_{\Omega_{1}}=10^{4},\qquad  a_{2}=a \chi_{\Omega_{2}}=10^{-4}, \qquad a_{3}=a \chi_{\Omega_{3}}=10^{4},\qquad  a_{4}=a \chi_{\Omega_{4}}=10^{-4},\\
		& {   u_{1}=u \chi_{\Omega_{1}}=\cos(5x)\cos(5y),\qquad  u_{2}=u \chi_{\Omega_{2}}=\cos(12x)\exp(y-x),}\\
		&  {  u_{3}=u \chi_{\Omega_{3}}=\sin(5x)\cos(5y),\qquad  u_{4}=u \chi_{\Omega_{4}}=\sin(12y)\exp(x-y),}
	\end{align*}
and $g, f_{p}$, $\phi_p$, $\psi_p$ for $p=1,2,3,4$ in \eqref{intersect:3} are obtained by plugging the above functions into \eqref{intersect:3}.
	Note that $\phi_p \ne 0$  for $p=1,2,3,4$.
	The numerical results are presented in \cref{Intersect:Numeri:table2} and \cref{Intersect:Numeri:fig4}.	
\end{example}

%%%%%%%%%%%%%%%%%%%%%%%%%%%%%%%%%%%%%%%%%%%%
%%%%%%%%%%%%%%%%%%%%%%%%%%%%%%%%%%%%%%%%%%%%%
%%%%%%%%%%%%%%%%%%%%%%%%%%%%%%%%%%%%%%%%%%%%%
%%%%%%%%%%%%%%%%%%%%%%%%%%%%%%%%%%%%%%%%%%%%%

\begin{figure}[htbp]
	\centering
	\begin{subfigure}[b]{0.3\textwidth}
		\hspace{0.5cm}
		\begin{tikzpicture}[scale = 4.35]
			\draw	(0, 0) -- (0, 1) -- (1, 1) -- (1, 0) --(0,0);	
			\draw   (pi/4, 0) -- (pi/4, 1);	
			\draw   (0, pi/10) -- (1, pi/10);	
			 %%%%%%%%%%%%%%%%%%%%%%%%%%%%%%%%%%%%%%%%%%%%%%%%%
			\node (A)[scale = 0.9] at (0.35,0.67) {$a_1=10^{4}$};
			\node (A)[scale = 0.9] at (0.9,0.67) {$a_2$};
			\node (A)[scale = 0.9] at (0.9,0.57) {$10^{-6}$};
			\node (A)[scale = 0.9] at (0.9,0.19) {$a_3$};
			\node (A)[scale = 0.9] at (0.9,0.09) {$10^{5}$};
			\node (A)[scale = 0.9] at (0.38,0.17) {$a_4=10^{-5}$};
		\end{tikzpicture}
	\end{subfigure}
	\begin{subfigure}[b]{0.3\textwidth}
		 \includegraphics[width=5.5cm,height=4.5cm]{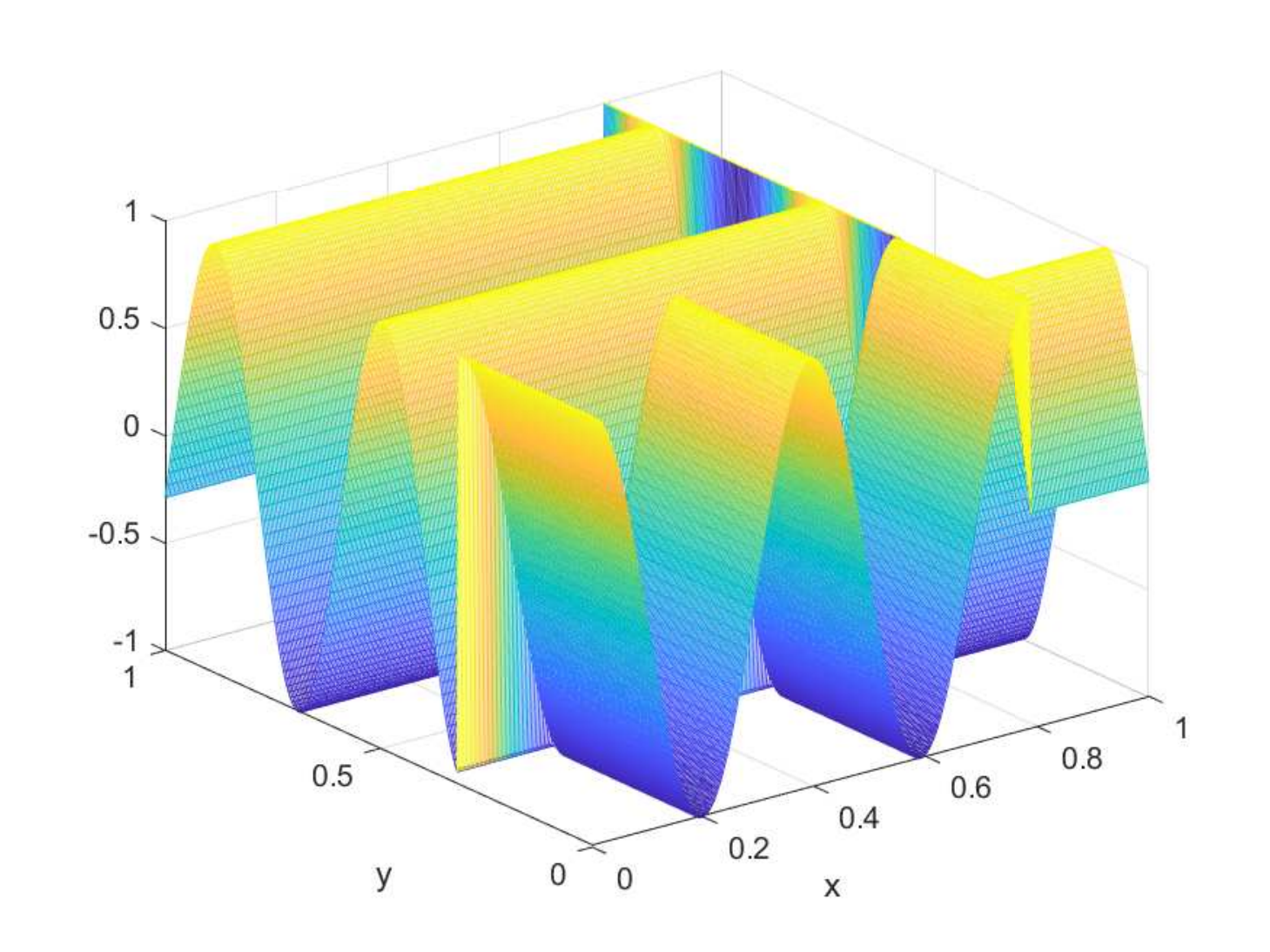}
	\end{subfigure}
	\begin{subfigure}[b]{0.3\textwidth}
		 \includegraphics[width=5.5cm,height=4.5cm]{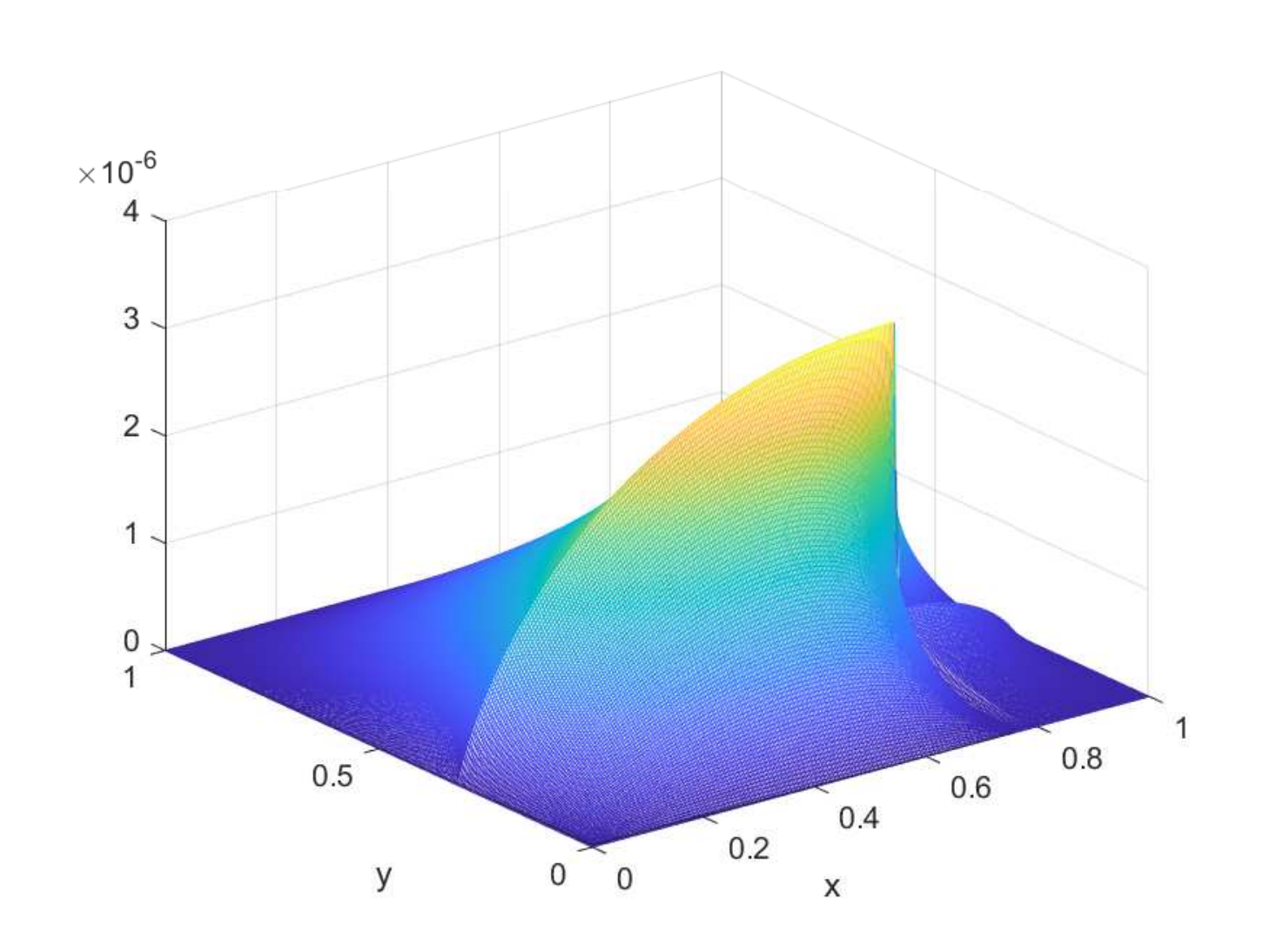}
	\end{subfigure}
	\caption
	{\cref{Intersect:ex5}:  the coefficient $a(x,y)$ (left),  the numerical solution $(u_h)_{i,j}$ (middle), and  $|(u_h)_{i,j}-u(x_i,y_j)|$ (right), at all grid points $(x_i,y_j)$ on $\Omega$ with $h=2^{-8}$,  for the scheme in \cref{Intersect:thm:regular,theorem:side3:w,Intersect:thm:cross1:w1w2}.}
	\label{Intersect:Numeri:fig5}
\end{figure}

%%%%%%%%%%%%%%%%%%%%%%%%%%%%%%%%%%%%%%%%%%%%%
\begin{example}\label{Intersect:ex5}
	\normalfont
	Let $\Omega=(0,1)^2$ and $(\xi,\zeta)=(\pi/4,\pi/10)$.
	The functions  in \eqref{intersect:3} are given by
	%%%
	\begin{align*}
		&a_{1}=a \chi_{\Omega_{1}}=10^{4},\qquad  a_{2}=a \chi_{\Omega_{2}}=10^{-6}, \qquad a_{3}=a \chi_{\Omega_{3}}=10^{5},\qquad  a_{4}=a \chi_{\Omega_{4}}=10^{-5},\\
		& {  u_{1}=u \chi_{\Omega_{1}}=\sin(16y)
		,\quad  u_{2}=u \chi_{\Omega_{2}}=\cos(16x), \quad
		u_{3}=u \chi_{\Omega_{3}}=\sin(16y)
		,\quad  u_{4}=u \chi_{\Omega_{4}}=\cos(16x),}
	\end{align*}
and	$g, f_{p}$, $\phi_p$, $\psi_p$ for $p=1,2,3,4$ in \eqref{intersect:3} are obtained by plugging the above functions into \eqref{intersect:3}.
		Note that $\phi_p \ne 0$  for $p=1,2,3,4$.
	The numerical results are presented in \cref{Intersect:Numeri:table2} and \cref{Intersect:Numeri:fig5}.	
\end{example}

%%%%%%%%%%%%%%%%%%%%%%%%%%%%%%%%%%%%%%%%%%%%
%%%%%%%%%%%%%%%%%%%%%%%%%%%%%%%%%%%%%%%%%%%%%
%%%%%%%%%%%%%%%%%%%%%%%%%%%%%%%%%%%%%%%%%%%%%
%%%%%%%%%%%%%%%%%%%%%%%%%%%%%%%%%%%%%%%%%%%%%

\begin{figure}[htbp]
	\centering
	\begin{subfigure}[b]{0.3\textwidth}
		\hspace{0.5cm}
		\begin{tikzpicture}[scale = 4.35]
			\draw	(0, 0) -- (0, 1) -- (1, 1) -- (1, 0) --(0,0);	
			\draw   (pi/6, 0) -- (pi/6, 1);	
			\draw   (0, pi/8) -- (1, pi/8);	
			 %%%%%%%%%%%%%%%%%%%%%%%%%%%%%%%%%%%%%%%%%%%%%%%%%
			\node (A)[scale = 0.9] at (0.3,0.7) {$a_1=10^{-4}$};
			\node (A)[scale = 0.9] at (0.78,0.7) {$a_2=10^{5}$};
			\node (A)[scale = 0.9] at (0.78,0.2) {$a_3=10^{-4}$};
			\node (A)[scale = 0.9] at (0.3,0.2) {$a_4=10^{6}$};
		\end{tikzpicture}
	\end{subfigure}
	\begin{subfigure}[b]{0.3\textwidth}
		 \includegraphics[width=5.5cm,height=4.5cm]{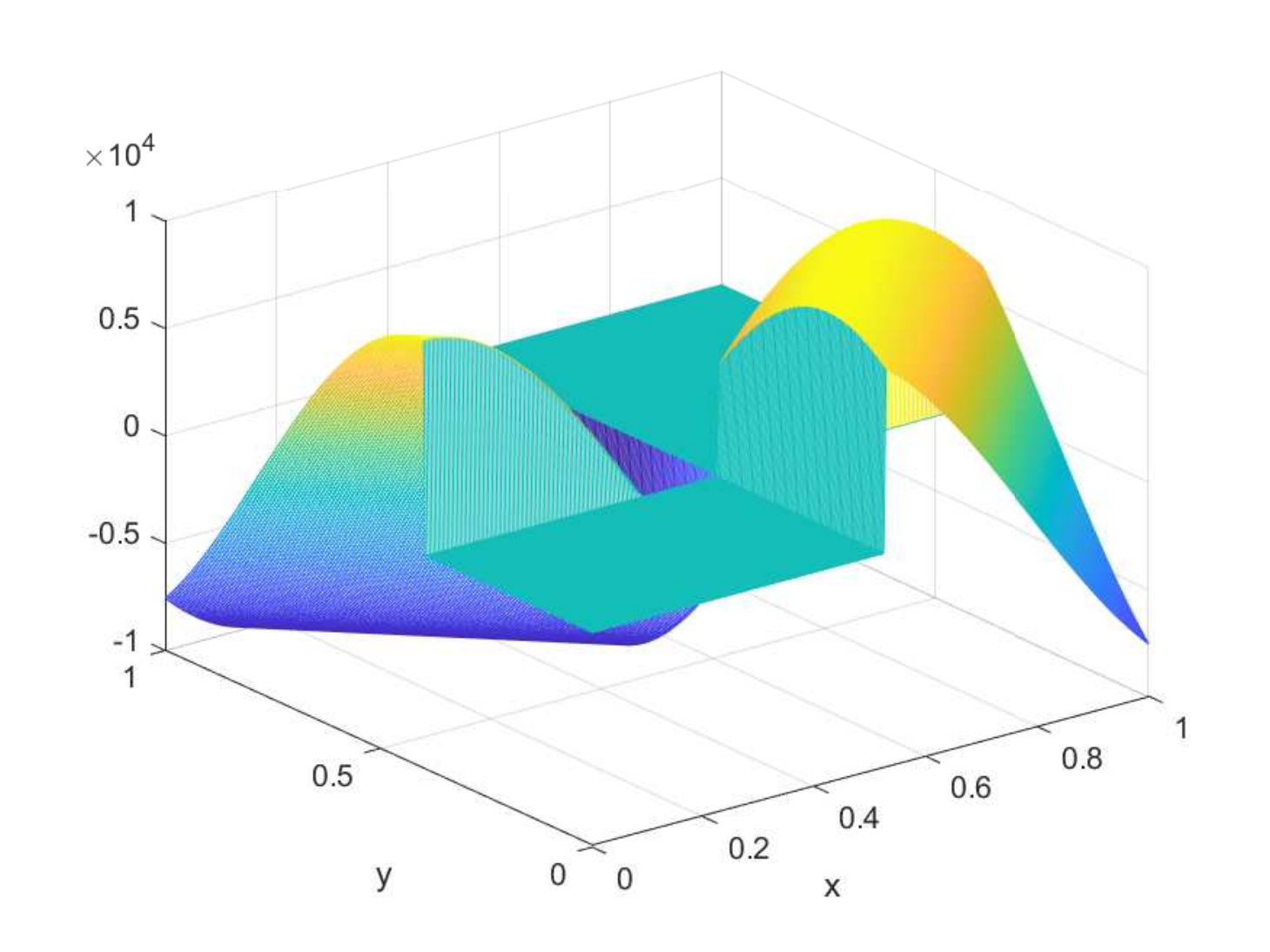}
	\end{subfigure}
	\begin{subfigure}[b]{0.3\textwidth}
		 \includegraphics[width=5.5cm,height=4.5cm]{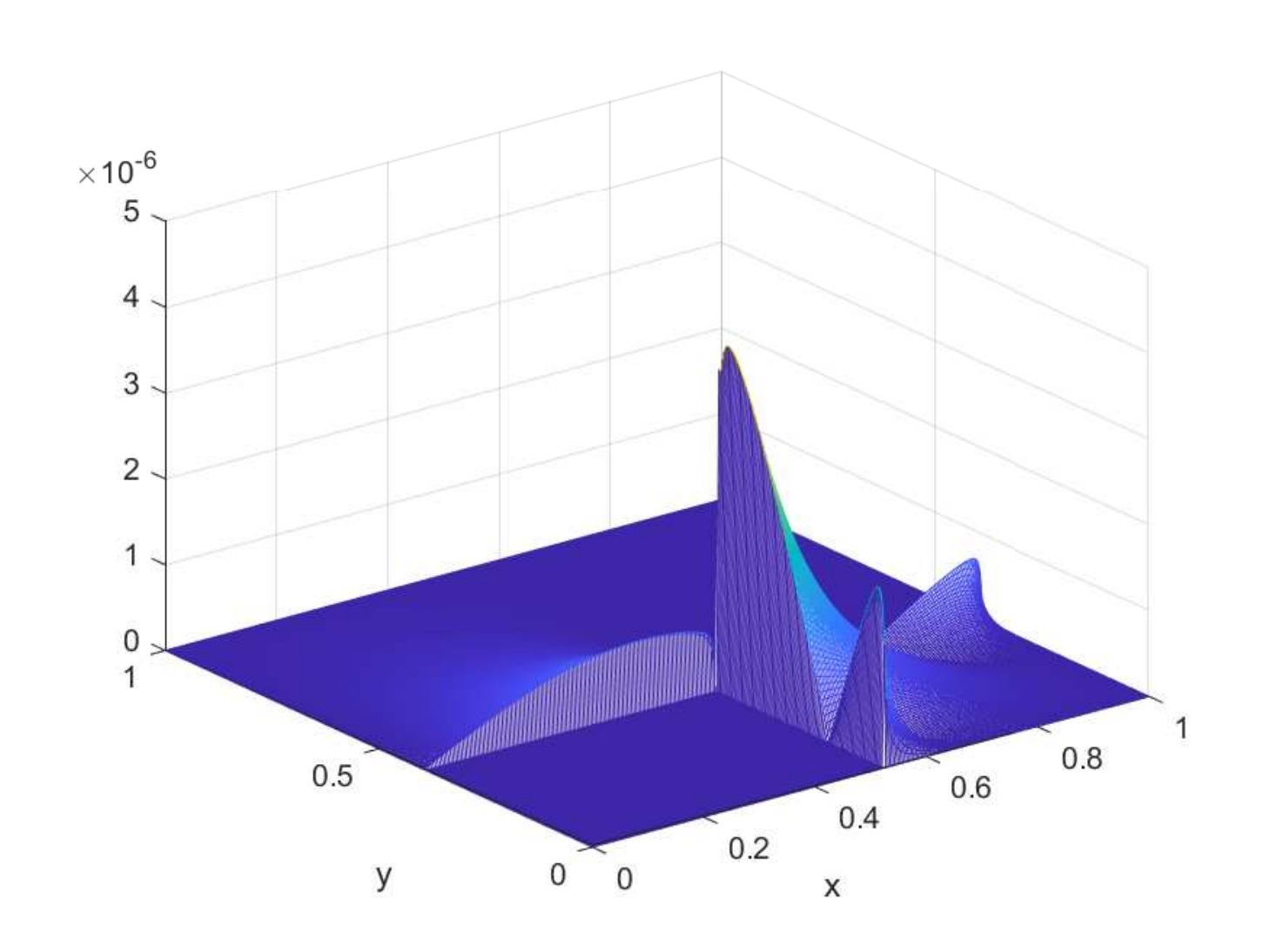}
	\end{subfigure}
	\caption
	{\cref{Intersect:ex6}:  the coefficient $a(x,y)$ (left),  the numerical solution $(u_h)_{i,j}$ (middle), and $|(u_h)_{i,j}-u(x_i,y_j)|$ (right), at all grid points $(x_i,y_j)$ on $\Omega$ with $h=2^{-8}$,  for the scheme  in \cref{Intersect:thm:regular,theorem:side3:w,Intersect:thm:cross1:w1w2}.}
	\label{Intersect:Numeri:fig6}
\end{figure}

%%%%%%%%%%%%%%%%%%%%%%%%%%%%%%%%%%%%%%%%%%%%%
\begin{example}\label{Intersect:ex6}
	\normalfont
	Let $\Omega=(0,1)^2$ and $(\xi,\zeta)=(\pi/6,\pi/8)$.
	The functions  in \eqref{intersect:3} are given by
	%%%
	\begin{align*}
		&a_{1}=a \chi_{\Omega_{1}}=10^{-4},\qquad  a_{2}=a \chi_{\Omega_{2}}=10^{5}, \qquad a_{3}=a \chi_{\Omega_{3}}=10^{-4},\qquad  a_{4}=a \chi_{\Omega_{4}}=10^{6},\\
		&  {  u_{1}=u \chi_{\Omega_{1}}=\sin(4(x+y))/a_1,\qquad  u_{2}=u \chi_{\Omega_{2}}=\cos(2(x-y))/a_2,}\\
		&  {   u_{3}=u \chi_{\Omega_{3}}=\sin(4(x-y))/a_3,\qquad  u_{4}=u \chi_{\Omega_{4}}=\cos(2(x+y))/a_4,}
	\end{align*}
and	$g, f_{p}$, $\phi_p$, $\psi_p$ for $p=1,2,3,4$ in \eqref{intersect:3} are obtained by plugging the above functions into \eqref{intersect:3}.
		Note that $\phi_p \ne 0$  for $p=1,2,3,4$.
	The numerical results are presented in \cref{Intersect:Numeri:table2} and \cref{Intersect:Numeri:fig6}.	
\end{example}

%%%%%%%%%%%%%%%%%%%%%%%%%%%%%%%%%%%%%%%%%%%%
%%%%%%%%%%%%%%%%%%%%%%%%%%%%%%%%%%%%%%%%%%%%%
%%%%%%%%%%%%%%%%%%%%%%%%%%%%%%%%%%%%%%%%%%%%%
%%%%%%%%%%%%%%%%%%%%%%%%%%%%%%%%%%%%%%%%%%%%%

\section{Conclusion}\label{Intersect:sec:Conclu}

This paper is aimed at the development of compact finite difference schemes of high orders for the cross-interface problems in \eqref{intersect:3}, under the assumptions (A1)-(A3).  The main contributions of this paper can be summarized as follows:

\begin{enumerate}
	\item[(1)] We derive  compact  9-point  finite difference schemes that are sixth-order accurate, if the internal interfaces coincide with some grid lines, and fourth-order accurate  otherwise.  The schemes use uniform meshes and all formulas of the schemes are constructed explicitly for all grid points which can be easily implemented.
	
	\item[(2)] In case that the internal interfaces of coefficient jumps are matched by grid lines, the resulting scheme satisfies the M-matrix property and the discrete maximum principle, that allows us to prove that it has the sixth-order convergence rate.
	If the interfaces are not matched by grid lines, the resulting linear system does not satisfy the M-matrix property, and we were unable to theoretically prove its  fourth/fifth-order  convergence rate.  The fourth/fifth-order convergence rate was verified only numerically.
	
	\item[(3)]  In the latter case, we derive a compact,  9-point, third-order accurate scheme satisfying the M-matrix property, and we prove its third-order convergence rate using  the discrete maximum principle.
\end{enumerate}

\section{Proofs of Results Stated in \cref{Intersect:sec:sixord,error:analysis:intersect:sec}}
\label{Intersect:sec:proofs}
\begin{figure}[h]
	\centering
	\hspace{12mm}	
	\begin{tikzpicture}[scale = 0.8]
		\node  [scale=1](A) at (0,7) {${(0,0)}$};
		\node  [scale=1](A) at (0,6) {${(0,1)}$};
		\node  [scale=1] (A) at (0,5) {${(0,2)}$};
		\node  [scale=1](A) at (0,4) {${(0,3)}$};
		\node  [scale=1](A) at (0,3) {${(0,4)}$};
		\node  [scale=1](A) at (0,2) {${(0,5)}$};
		\node  [scale=1](A) at (0,1) {${(0,6)}$};
		\node  [scale=1](A) at (0,0) {${(0,7)}$};
		%%%%%%%%%%%%%%%%%%%%%%%%%%%%%%%%%%%%%
		\node  [scale=1](A) at (2,7) {${(1,0)}$};
		\node  [scale=1](A) at (2,6) {${(1,1)}$};
		\node  [scale=1](A) at (2,5) {${(1,2)}$};
		\node  [scale=1](A) at (2,4) {${(1,3)}$};
		\node  [scale=1](A) at (2,3) {${(1,4)}$};
		\node  [scale=1](A) at (2,2) {${(1,5)}$};
		\node  [scale=1](A) at (2,1) {${(1,6)}$};	 
		%%%%%%%%%%%%%%%%%%%%%%%%%%%%%%%%%%%%%
		\node  [scale=1](A) at (4,7) {${(2,0)}$};
		\node  [scale=1](A) at (4,6) {${(2,1)}$};
		\node  [scale=1](A) at (4,5) {${(2,2)}$};
		\node  [scale=1](A) at (4,4) {${(2,3)}$};
		\node  [scale=1](A) at (4,3) {${(2,4)}$};
		\node  [scale=1](A) at (4,2) {${(2,5)}$};
		%%%%%%%%%%%%%%%%%%%%%%%%%%%%%%%%%%%%%
		\node  [scale=1](A) at (6,7) {${(3,0)}$};
		\node  [scale=1](A) at (6,6) {${(3,1)}$};
		\node  [scale=1](A) at (6,5) {${(3,2)}$};
		\node  [scale=1](A) at (6,4) {${(3,3)}$};
		\node  [scale=1](A) at (6,3) {${(3,4)}$};
		%%%%%%%%%%%%%%%%%%%%%%%%%%%%%%%%%%%%%
		\node  [scale=1](A) at (8,7) {${(4,0)}$};
		\node  [scale=1](A) at (8,6) {${(4,1)}$};
		\node  [scale=1](A) at (8,5) {${(4,2)}$};
		\node  [scale=1](A) at (8,4) {${(4,3)}$};
		%%%%%%%%%%%%%%%%%%%%%%%%%%%%%%%%%%%%%
		\node  [scale=1](A) at (10,7) {${(5,0)}$};
		\node  [scale=1](A) at (10,6) {${(5,1)}$};
		\node  [scale=1](A) at (10,5) {${(5,2)}$};
		%%%%%%%%%%%%%%%%%%%%%%%%%%%%%%%%%%%%%
		\node  [scale=1](A) at (12,7) {${(6,0)}$};
		\node  [scale=1](A) at (12,6) {${(6,1)}$};
		%%%%%%%%%%%%%%%%%%%%%%%%%%%%%%%%%%%%%
		\node  [scale=1](A) at (14,7) {${(7,0)}$};      	
		\draw[line width=1.5pt, red]  plot [tension=0.8]
		coordinates {(-0.75,7.5) (-0.75,-0.7) (2.7,0.6) (2.7,7.5) (-0.75,7.5)};	
		\draw[line width=1.5pt, blue]  plot [tension=0.8]
		coordinates {(3,7.5)  (3,1.0) (14.7,6.8) (14.7,7.5) (3,7.5)};
		 %%%%%%%%%%%%%%%%%%%%%%%%%%%%%%%%%%%%%%%%%
		 \draw[decorate,decoration={brace,mirror,amplitude=4mm},xshift=0pt,yshift=10pt,ultra thick] (-0.75,-1) -- node [black,midway,yshift=0.6cm]{} (2.7,-1);
		\node (A) at (1.3,-1.5)	 {$\{{(m,n)} \; :\; (m,n)\in \ind_{7}^{ 1}\}$};
		 %%%%%%%%%%%%%%%%%%%%%%%%%%%%%%%%%%%%%%%%%
		 \draw[decorate,decoration={brace,mirror,amplitude=4mm},xshift=0pt,yshift=10pt,ultra thick] (3,0.5) -- node [black,midway,yshift=0.6cm]{} (14.7,0.5);
		\node (A) at (9,0)	 {$\{{(m,n)} \; :\; (m,n)\in \ind_{7}^{ 2}\}$};
	\end{tikzpicture}
	\caption
	{Red trapezoid: $\{{(m,n)} \; :\; (m,n)\in \ind_{7}^{1}\}$. Blue trapezoid: $\{{(m,n)} \; :\; (m,n)\in \ind_{7}^{2}\}$. Note that $\ind_{7}^{ 1} \cup\ind_{7}^{2}=\ind_{7}:=\{(m,n)\in \N_0^2 \; :\; m+n\le 7\}$.}
	\label{Intersect:fig:umnV}
\end{figure}

\begin{figure}[htbp]
	\centering	
	\begin{subfigure}[b]{0.3\textwidth}
		\begin{tikzpicture}[scale = 1]
			\draw	(-2.5, -2.5) -- (-2.5, 2.5) -- (2.5, 2.5) -- (2.5, -2.5) --(-2.5,-2.5);	
			\draw   (0, -2.5) -- (0, 2.5);	
			\draw   (-2.5, 0) -- (2.5, 0);	
			\node (A) at (-1.5,1.5) {$u_1$};
			\node (A) at (1.7,1.5) {$u_2$};
			\node (A) at (1.7,-1.5) {$u_3$};
			\node (A) at (-1.5,-1.5) {$u_4$};
			 %%%%%%%%%%%%%%%%%%%%%%%%%%%%%%%%%%%%%%%%%%
			 %%%%%%%%%%%%%%%%%%%%%%%%%%%%%%%%%%%%%%%%%%%%%%%%%
			\node (A) at (-0.3,-1.5) {$\Gamma_2$};
			\node (A) at (-0.3,1.4) {$\Gamma_1$};
			\node (A) at (1.6,0.3) {$\Gamma_3$};
			\node (A) at (-1.6,0.3) {$\Gamma_4$};
			%%%%%%%%%%%%%%%%%%%%
			\node (A)[scale=0.8] at (0.6,-0.3) {$(x_i^*,y_j^*)$};
			\node at (0,0)[circle,fill,inner sep=2pt,color=black]{};
			%%%%%%%%%%%%%%%%%%%
			\draw[ultra thick][->] (-1, -1) to[bend left] (-1, 1);
			\draw[ultra thick][->] (1, -1) to[bend right] (1, 1);
			\draw[ultra thick][->] (-1, 1.5) to[bend left] (1, 1.5);
		\end{tikzpicture}
	\end{subfigure}
	 %%%%%%%%%%%%%%%%%%%%%%%%%%%%%%%%%%%%%%%%%%%%%
	\begin{subfigure}[b]{0.3\textwidth}
		\begin{tikzpicture}[scale = 1]
			\draw	(-2.5, -2.5) -- (-2.5, 2.5) -- (2.5, 2.5) -- (2.5, -2.5) --(-2.5,-2.5);	
			\draw   (0, -2.5) -- (0, 2.5);	
			\draw   (-2.5, 0) -- (2.5, 0);	
			\node (A) at (-1.5,1.5) {$u_1$};
			\node (A) at (1.7,1.5) {$u_2$};
			\node (A) at (1.7,-1.5) {$u_3$};
			\node (A) at (-1.5,-1.5) {$u_4$};
			 %%%%%%%%%%%%%%%%%%%%%%%%%%%%%%%%%%%%%%%%%%
			 %%%%%%%%%%%%%%%%%%%%%%%%%%%%%%%%%%%%%%%%%%%%%%%%%
			\node (A) at (-0.3,-1.5) {$\Gamma_2$};
			\node (A) at (-0.3,1.4) {$\Gamma_1$};
			\node (A) at (1.6,0.3) {$\Gamma_3$};
			\node (A) at (-1.5,0.3) {$\Gamma_4$};
			%%%%%%%%%%%%%%%%%%%%
			\node (A)[scale=0.8] at (0.6,-0.3) {$(x_i^*,y_j^*)$};
			\node at (0,0)[circle,fill,inner sep=2pt,color=black]{};
			%%%%%%%%%%%%%%%%%%%
			\draw[ultra thick][->] (-1, -1.5) to[bend right] (1, -1.5);
			\draw[ultra thick][->] (1, -1) to[bend right] (1, 1);
		%%%%	\draw[ultra thick][->] (-1, 1.5) to[bend left] (1, 1.5);
		\end{tikzpicture}
	\end{subfigure}
	%%%%%%%%%%
	\caption{}
	\label{transsmi:4:parts}
\end{figure}

For the proofs of the theorems in this paper we first need to establish
 some auxiliary identities about the solution $u$ of the elliptic cross-interface problem in \eqref{intersect:3}.
Let $(x_i^*,y_j^*) \in  \Gamma:=\Gamma_1 \cup \Gamma_2 \cup \Gamma_3 \cup \Gamma_4\cup\{(\xi,\zeta)\}$ such that $x_i^*\in (x_i-h,x_i+h)$ and $y_j^*\in (y_j-h,y_j+h)$.
Recall that for $p=1,\ldots,4$,
\[
u_p:=u\chi_{\Omega_p},\quad f_p:=f\chi_{\Omega_p},\quad
u^{(m,n)}_p:=\frac{\partial^{m+n}u_p}{\partial x^m \partial y^n}(x_i^*,y_j^*),\quad f^{(m,n)}_p:=\frac{\partial^{m+n}f_p}{\partial x^m \partial y^n}(x_i^*,y_j^*).
\]
As in \cite{FengHanMinev21,FHM21Helmholtz},
we can derive from $-\nabla \cdot(a\nabla u)=f$ that
	\be \label{intersect:uderiv:relation}
	\begin{split}
		& u_p^{(m,n)}=(-1)^{\lfloor\frac{m}{2}\rfloor}
		u_p^{(\odd(m),n+m-\odd(m))}+
			\frac{1}{a_p} \sum_{\ell=1}^{\lfloor m/2\rfloor}
	 (-1)^{\ell} f_p^{(m-2\ell,n+2\ell-2)},
\qquad \forall\; (m,n)\in \ind_{M}^{2},
\\
		& u_p^{(m,n)}=(-1)^{\lfloor\frac{n}{2}\rfloor}
		u_p^{(n+m-\odd(n),\odd(n))}+
			\frac{1}{a_p} \sum_{\ell=1}^{\lfloor n/2\rfloor}
	(-1)^{\ell} f_p^{(m+2\ell-2,n-2\ell)},
\qquad \forall\; (n,m)\in \ind_{M}^{2},
	\end{split}
	\ee
	%%%%
where the floor function $\lfloor x\rfloor$ is defined to be the largest integer less than or equal to $x\in \R$, and
\be\label{intersect:oddm}
\odd(m):=\begin{cases}
	0, &\text{if $m$ is even},\\
	1, &\text{if $m$ is odd}.
\end{cases}
\ee
Similarly to \cite[(2.5)-(2.10)]{FengHanMinev21} and
\cite[(2.9)-(2.13)]{FHM21Helmholtz}, using
the above identities, and using Taylor expansions at the point $(x_i^*,y_j^*)$, we obtain
	{\small{
	\begin{align}
u_p(x+x_i^*,y+y_j^*)&=
		\sum_{(m,n)\in \ind_{M}^{1}}
		u_p^{(m,n)} G_{M,m,n}(x,y)  + \frac{1}{a_p} \sum_{(m,n)\in \ind_{M-2}}
		f_p^{(m,n)} H_{M,m,n}(x,y)+\bo(h^{M+1}),
\label{Intersect:u:approx:key:V}\\
u_p(x+x_i^*,y+y_j^*)&=
		\sum_{(n,m)\in \ind_{M}^{1}}
		u_p^{(m,n)} G_{M,n,m}(y,x) +  	 \frac{1}{a_p} \sum_{(m,n)\in \ind_{M-2}}
			f_p^{(m,n)} H_{M,n,m}(y,x)+\bo(h^{M+1}), \label{Intersect:u:approx:key:H}
	\end{align}
}}
	for $(x,y) \in(-2h,2h)^2$, where the bivariate polynomials $G$ and $H$ are defined in \eqref{Intersect:GmnV}.

From $[u]=\phi_1$ and $\left[a \nabla  u \cdot \vec{n}\right]=\psi_1$ on $\Gamma_1$ in \eqref{intersect:3},
we obtain
\be\label{u1u2:jump:1}
u_1^{(0,n)}=u_2^{(0,n)}-\phi_1^{(n)}, \quad u_1^{(1,n)}=\frac{a_2}{a_1}u_2^{(1,n)}-\frac{1}{a_1}\psi_1^{(n)}, \quad n = 0,1,2,\dots, M.
\ee
The above identities will be frequently used in the following proofs.

	\begin{proof}[Proof of \cref{Intersect:thm:side1}]
Note that $(x_i,y_j)=(x_i^*,y_j^*) \in \Gamma_1$.
Since the discrete operator $\mathcal{L}_h u$ in \eqref{Lu} involves both $u_1$ and $u_2$ (see the first panel of \cref{fig:intersect:side1}), using jump conditions in \eqref{intersect:3} across the interface $\Gamma_1$, we can replace all $\{u_1^{(m,n)} \; :\; (m,n)\in \ind_M^1 \}$ on the right-hand side of \eqref{Intersect:u:approx:key:V} with $p=1$ by $\{ u_2^{(m,n)} \; :\; (m,n)\in \ind^1_M\}$. More precisely, using identities in \eqref{u1u2:jump:1} to replace $\{  u_1^{(m,n)} \; :\; (m,n) \in \ind^1_M \}$ in
\eqref{Intersect:u:approx:key:V} with $p=1$, we obtain
{\small{
\be\label{intersect:u1xy:side1:ir:thm1}	
\begin{split}
u_1(x+x_i^*,y+y_j^*)
	&=\sum_{n=0}^{M}
u_1^{(0,n)} G_{M,0,n}(x,y) +\sum_{n=0}^{M-1}
u_1^{(1,n)} G_{M,1,n}(x,y)+ \sum_{(m,n)\in \ind_{M-2}}
\frac{f_1^{(m,n)}}{a_1}  H_{M,m,n}(x,y),\\
&=\sum_{n=0}^{M}
\big( u_2^{(0,n)}	 -\phi_1^{(n)} \big) G_{M,0,n}(x,y)
+\sum_{n=0}^{M-1}
\Big(\frac{a_2}{a_1}u_2^{(1,n)}-\frac{\psi_1^{(n)}}{a_1} \Big) G_{M,1,n}(x,y)\\
&\quad +\frac{1}{a_1}\sum_{(m,n)\in \ind_{M-2}}
f_1^{(m,n)} H_{M,m,n}(x,y),\\
&=	 \sum_{(m,n)\in \ind_{M}^{1}}
u_2^{(m,n)} \Big(\frac{a_2}{a_1}\Big)^m G_{M,m,n}(x,y) - \sum_{n=0}^{M}
\phi_1^{(n)} G_{M,0,n}(x,y) \\
&\quad- \frac{1}{a_1}
\sum_{n=0}^{M-1}\psi_1^{(n)} G_{M,1,n}(x,y)   + \frac{1}{a_1} \sum_{(m,n)\in \ind_{M-2}}
f_1^{(m,n)} H_{M,m,n}(x,y)+\bo(h^{M+1}),
\end{split}
\ee}}
%%%
for  $(x,y) \in (-2h,2h)^2$.	
Using \eqref{Intersect:u:approx:key:V} with $p=2$, \eqref{intersect:u1xy:side1:ir:thm1}, and
%\eqref{intersect:u2xy:side1:ir:thm1} with
$(x_{i}^*,y_{j}^*)=(x_{i},y_{j})$,
we deduce from the definition of $\mathcal{L}_h u$ in \eqref{Lu} that
%and $u_{p}:=u \chi_{\Omega_{p}}$ for $p=1,2,3,4$ that
%

	\be\label{Lhu:Gamma1}
\begin{split}
	\mathcal{L}_h u  =&\begin{aligned}	
		&C_{-1,1}u_1(x_i-h,y_j+h)&
		+&C_{0,1}u_2(x_i,y_j+h)&
		+&C_{1,1}u_2(x_i+h,y_j+h)\\
		+&C_{-1,0} u_1(x_i-h,y_j)&
		+&C_{0,0}u_2(x_i,y_j)&
		+&C_{1,0}u_2(x_i+h,y_j)\\
		+&C_{-1,-1}u_1(x_i-h,y_j-h)&
		+&C_{0,-1}u_2(x_i,y_j-h)&
		+&C_{1,-1}u_2(x_i+h,y_j-h)\\
	\end{aligned}\\
	=& \sum_{(m,n)\in \ind_{M}^{1}}
	u_2^{(m,n)}  I_{m,n}+ \frac{1}{a_1}\sum_{(m,n)\in \ind_{M-2}} f_1^{(m,n)}
	H^{-,0}_{M,m,n} +\frac{1}{a_2}\sum_{(m,n)\in \ind_{M-2}} f_2^{(m,n)}
	H^{+,0}_{M,m,n}\\
	&-\sum_{n=0}^{M} \phi_1^{(n)} G^0_{M,0,n}
	-\frac{1}{a_1} \sum_{n=0}^{M-1} \psi_1^{(n)}  G^0_{M,1,n}
	+\bo(h^{M+1}),\qquad {as} \quad h \to 0,
\end{split}
\ee
where $H^{-,0}_{M,m,n}, H^{+,0}_{M,m,n}, G^0_{M,m,n}$ are defined in \eqref{GHw:v}, and
\be \label{Intersect:Imn:Proof:Thm1:Irre}
I_{m,n}:=
\Big(\frac{a_2}{a_1}\Big)^m \sum_{\ell=-1}^1 C_{-1,\ell}   G_{M,m,n}(-h, \ell h)
+\sum_{k=0}^1 \sum_{\ell=-1}^1 C_{k,\ell}  G_{M,m,n}(kh, \ell h).
\ee
%
%and $H^{\pm,0}_{M,m,n},
%G^0_{M,m,n}$ are defined in \eqref{GHw:v} (noting that $w=0$ due to $(x_i,y_j)=(x_i^*,y_j^*)$).
%
%
Then the conditions
\be\label{intersect:Imn:Conditions}
I_{m,n}=\bo(h^{M+1}), \qquad \mbox{ for all } \quad (m,n)\in \ind_{M}^{1},
\ee
can be equivalently rewritten as a system of linear equations on the unknowns $\{C_{k,\ell}\}_{k,\ell=-1,0,1}$. By calculation, we observe that
\eqref{intersect:Imn:Conditions} has a nontrivial solution $\{C_{k,\ell}\}_{k,\ell=-1,0,1}$ if and only if $M\le 7$. Moreover, for $M=7$, up to a multiplicative constant for normalization, \eqref{C:Gamma1} is the unique solution to
\eqref{intersect:Imn:Conditions}. Therefore, for the solution  $\{C_{k,\ell}\}_{k,\ell=-1,0,1}$ in \eqref{C:Gamma1} with $M=7$,
we conclude that
\[
h^{-1}	\mathcal{L}_h (u-u_h)=
h^{-1} \sum_{(m,n)\in \ind_{7}^{1}}
u_2^{(m,n)}  I_{m,n}=\bo(h^{7}),\qquad h\to 0,
\]
which proves the seventh order of consistency of
the compact finite difference scheme in \cref{Intersect:thm:side1}.
\end{proof}

%%%%%%%%%%%%%%%%%%%%%%%%
%\begin{proof}[Proof of \cref{Intersect:thm:side3}]
%By symmetry, the proof is similar to the proof of \cref{Intersect:thm:side1}.
%\end{proof}
%%%%%%%%%%%%%%%%%%%%%%%
%%
%
%
%
%%
\begin{proof}[Proof of \cref{Intersect:thm:cross1}]
Note that $(x_i,y_j)=(x_i^*,y_j^*)=(\xi,\zeta)$.
Since the discrete operator $\mathcal{L}_h u$ in \eqref{Lu} involves $u_p$ with $p=1,2,3,4$ (see the third panel of \cref{fig:intersect:side1}), using jump conditions in \eqref{intersect:3} across interfaces $\Gamma_1$, $\Gamma_3$ and $\Gamma_4$ (see the first panel of \cref{transsmi:4:parts}), we can replace all $\{ u_3^{(m,n)} \; :\; (n,m)\in \ind_M^1\}$ and $\{ u_4^{(m,n)} \; :\; (n,m)\in \ind_M^1\}$ on the right-hand side of \eqref{Intersect:u:approx:key:H} with $p=3,4$ by $\{  u_2^{(m,n)} \; :\; (m,n)\in \ind^1_M \}$ from the following \eqref{intersect:cross:U3:0}--\eqref{u4xy:Final:intersect:V1}. Then the compact 9-point scheme with a consistency order seven is derived from \eqref{LShu:proof:intersect}--\eqref{trunca:order:7}.

Similarly to \eqref{intersect:u1xy:side1:ir:thm1},  $[u]=\phi_p$ and $\left[a \nabla  u \cdot \vec{n}\right]=\psi_p$ on $\Gamma_p$ with $p=3,4$ imply that:
{\small{
\be\label{intersect:cross:U3:0}
\begin{split}
	u_3(x+x_i^*,y+y_j^*)
	&=	\sum_{(n,m)\in \ind_{M}^{1}}
	u_2^{(m,n)} \Big(\frac{a_2}{a_3}\Big)^n G_{M,n,m}(y,x)-\sum_{m=0}^{M}
	\phi_3^{(m)} G_{M,0,m}(y,x)\\& \quad - \frac{1}{a_3} \sum_{m=0}^{M-1} \psi_3^{(m)} G_{M,1,m}(y,x) +\frac{1}{a_3} \sum_{(m,n)\in \ind_{M-2}}
	f_3^{(m,n)} H_{M,n,m}(y,x)+\bo(h^{M+1}),
\end{split}
\ee}}
%%%%%
%
%	
%%
%%
{\small{
\be\label{intersect:cross:U4:0}
\begin{split}
	u_4(x+x_i^*,y+y_j^*)
	&=	\sum_{(n,m)\in \ind_{M}^{1}}
	u_1^{(m,n)} \Big(\frac{a_1}{a_4}\Big)^n G_{M,n,m}(y,x) -\sum_{m=0}^{M} \phi_4^{(m)} G_{M,0,m}(y,x)\\
	&\quad -\frac{ 1}{a_4} \sum_{m=0}^{M-1} \psi_4^{(m)} G_{M,1,m}(y,x) +\frac{ 1}{a_4} \sum_{(m,n)\in \ind_{M-2}}	 f_4^{(m,n)}H_{M,n,m}(y,x)+\bo(h^{M+1}),
\end{split}
\ee}}
%%%%%
%%%
for  $x,y\in (-2h,2h)$.
On the other hand,  \eqref{intersect:oddm} and \eqref{u1u2:jump:1} lead to:
\[
u_1^{(\odd(m),n+m-\odd(m))}=\begin{cases}
	u_2^{(0,m+n)}-\phi_1^{(m+n)}, &\text{if $m$ is even},\\
	 \frac{a_2}{a_1}u_2^{(1,m+n-1)}-\frac{1}{a_1}\psi_1^{(m+n-1)}, &\text{if $m$ is odd},
\end{cases}
\qquad \mbox{for all }(m,n)\in \ind_{M}^{2},
\]	
i.e.,
\be\label{intersect:u1:V2}
\begin{split}
	u_1^{(\odd(m),n+m-\odd(m))}&=	 \Big(\frac{a_2}{a_1}\Big)^{\odd(m)}u_2^{(\odd(m),m+n-\odd(m))}-\odd(m+1)\phi_1^{(m+n)}\\
	& \quad-\frac{\odd(m)}{a_1}\psi_1^{(m+n-1)}, \qquad \mbox{for all }(m,n)\in \ind_{M}^{2}.
\end{split}
\ee
%%%%%%%%
\eqref{intersect:uderiv:relation} implies that:
{\small{	
		\be\label{intersect:uimnV}
		 u_p^{(m,n)}=(-1)^{\lfloor\frac{m}{2}\rfloor}
		u_p^{(\odd(m),n+m-\odd(m))}+
		\sum_{\ell=1}^{\lfloor m/2\rfloor}
		\frac{(-1)^{\ell}}{a_p} f_p^{(m-2\ell,n+2\ell-2)},\quad \forall\; (m,n)\in \ind_{M}^{2} \mbox{ and }p=1,2.
		\ee}}
From \eqref{intersect:u1:V2} and \eqref{intersect:uimnV} with $p=1$, we observe that:	
\be\label{U1:intersect:V2:Step1}
\begin{split}
	 u_1^{(m,n)}
	 &=(-1)^{\lfloor\frac{m}{2}\rfloor}\Big(\frac{a_2}{a_1}\Big)^{\odd(m)}u_2^{(\odd(m),m+n-\odd(m))}-(-1)^{\lfloor\frac{m}{2}\rfloor} \odd(m+1)\phi_1^{(m+n)}
	\\
	&\quad  -(-1)^{\lfloor\frac{m}{2}\rfloor} \frac{\odd(m)}{a_1}\psi_1^{(m+n-1)}+
	\frac{1}{a_1}	 \sum_{\ell=1}^{\lfloor m/2\rfloor}
	(-1)^{\ell} f_1^{(m-2\ell,n+2\ell-2)}, \qquad \mbox{for all }(m,n)\in \ind_{M}^{2}.
\end{split}
\ee
By \eqref{Intersect:Sk}, we have:
\be\label{left:set:4}
\{ (m,n): (n,m)\in \ind_{M}^{1} \} \setminus \{ (m,n): (m,n)\in \ind_{M}^{2}\}=\{(0,0),(0,1),(1,0),(1,1)\}.
\ee
Note that for $m=0,1$, the summation $\sum_{\ell=1}^{\lfloor m/2\rfloor}$ in \eqref{intersect:uimnV} and \eqref{U1:intersect:V2:Step1}  is empty.
 So \eqref{intersect:uimnV} with $p=2$ and \eqref{left:set:4} result in:
\be\label{U2:intersect:H1}
u_2^{(m,n)}=(-1)^{\lfloor\frac{m}{2}\rfloor}
u_2^{(\odd(m),n+m-\odd(m))}+ \frac{1}{a_2}
\sum_{\ell=1}^{\lfloor m/2\rfloor}
(-1)^{\ell} f_2^{(m-2\ell,n+2\ell-2)}, \qquad \forall\; (n,m)\in \ind_{M}^{1}.
\ee
%%
%%%
%%
From \eqref{u1u2:jump:1},
\be\label{intersect:u1u200}
\begin{split}
	 &u_1^{(0,0)}=u_2^{(0,0)}-\phi_1^{(0)}, \quad u_1^{(0,1)}=u_2^{(0,1)}-\phi_1^{(1)},\quad u_1^{(1,0)}=\frac{a_2}{a_1}u_2^{(1,0)}-\frac{\psi_1^{(0)}}{a_1}, \quad u_1^{(1,1)}=\frac{a_2}{a_1}u_2^{(1,1)}-\frac{\psi_1^{(1)}}{a_1}.
\end{split}
\ee
So  \eqref{U1:intersect:V2:Step1}, \eqref{left:set:4} and \eqref{intersect:u1u200} lead to:
\be\label{U1:intersect:H1}
\begin{split}
	u_1^{(m,n)}
	 &=(-1)^{\lfloor\frac{m}{2}\rfloor}\Big(\frac{a_2}{a_1}\Big)^{\odd(m)}u_2^{(\odd(m),m+n-\odd(m))} - (-1)^{\lfloor\frac{m}{2}\rfloor} \odd(m+1)\phi_1^{(m+n)}
	\\
	&\quad -(-1)^{\lfloor\frac{m}{2}\rfloor} \frac{\odd(m)}{a_1}\psi_1^{(m+n-1)} +
	 	\frac{1}{a_1} \sum_{\ell=1}^{\lfloor m/2\rfloor}
	(-1)^{\ell}  f_1^{(m-2\ell,n+2\ell-2)},\qquad \forall\; (n,m)\in \ind_{M}^{1}.
\end{split}
\ee
\eqref{intersect:cross:U3:0} and \eqref{U2:intersect:H1} imply that:
{\footnotesize{
\be\label{U3xy:Intersect:H1}
\begin{split}
	 &u_3(x+x_i^*,y+y_j^*)+\bo(h^{M+1})\\
	&=	\sum_{(n,m)\in \ind_{M}^{1}}
	(-1)^{\lfloor\frac{m}{2}\rfloor}
	u_2^{(\odd(m),n+m-\odd(m))} \Big(\frac{a_2}{a_3}\Big)^n G_{M,n,m}(y,x) + \frac{1}{a_2} \sum_{(n,m)\in \ind_{M}^{1}}
	\sum_{\ell=1}^{\lfloor m/2\rfloor}
	(-1)^{\ell} f_2^{(m-2\ell,n+2\ell-2)}  \Big(\frac{a_2}{a_3}\Big)^n \\
	& \quad \times G_{M,n,m}(y,x)  -\sum_{m=0}^{M} \phi_3^{(m)} G_{M,0,m}(y,x) -\frac{1}{a_3} \sum_{m=0}^{M-1} \psi_3^{(m)} G_{M,1,m}(y,x) +\frac{1}{a_3} \sum_{(m,n)\in \ind_{M-2}}
	f_3^{(m,n)} H_{M,n,m}(y,x),
\end{split}
\ee}}
%%
%%%%%
%%
\eqref{intersect:cross:U4:0} and \eqref{U1:intersect:H1} imply that:
{\footnotesize{
\be\label{U4xy:Intersect:H1}
\begin{split}
	 &u_4(x+x_i^*,y+y_j^*)+\bo(h^{M+1})\\
	&=	\sum_{(n,m)\in \ind_{M}^{1}}
	 (-1)^{\lfloor\frac{m}{2}\rfloor}\Big(\frac{a_2}{a_1}\Big)^{\odd(m)}u_2^{(\odd(m),m+n-\odd(m))}\Big(\frac{a_1}{a_4}\Big)^n G_{M,n,m}(y,x)  - \frac{1}{a_4} \sum_{m=0}^{M-1} \psi_4^{(m)} G_{M,1,m}(y,x)\\
	&\quad -\sum_{(n,m)\in \ind_{M}^{1}} \Big( \frac{\odd(m+1)}{	 (-1)^{\lfloor\frac{m}{2}\rfloor}}\phi_1^{(m+n)} + \frac{\odd(m)}{(-1)^{\lfloor\frac{m}{2}\rfloor} a_1}\psi_1^{(m+n-1)} \Big) \Big(\frac{a_1}{a_4}\Big)^n G_{M,n,m}(y,x) -\sum_{m=0}^{M} \phi_4^{(m)} G_{M,0,m}(y,x)
	\\
	&\quad + 	\frac{1}{a_1} \sum_{(n,m)\in \ind_{M}^{1}}
	\sum_{\ell=1}^{\lfloor m/2\rfloor}
   (-1)^{\ell}  f_1^{(m-2\ell,n+2\ell-2)}  \Big(\frac{a_1}{a_4}\Big)^n G_{M,n,m}(y,x) +  \frac{	1 }{a_4} \sum_{(m,n)\in \ind_{M-2}}
   f_4^{(m,n)} H_{M,n,m}(y,x).
\end{split}
\ee
}}
By  \eqref{intersect:oddm} and the definition of $\ind_{M}^{1}$ in \eqref{Intersect:Sk}, we have:
\begin{align*}
	&\sum_{(n,m)\in \ind_{M}^{1}}
	 (-1)^{\lfloor\frac{m}{2}\rfloor}\Big(\frac{a_2}{a_1}\Big)^{\odd(m)}u_2^{(\odd(m),m+n-\odd(m))}  G_{M,n,m}(y,x)\\
	&=\sum_{v=0}^{\lfloor M/2\rfloor}
	(-1)^{v}u_2^{(0,2v)}  G_{M,0,2v}(y,x) +\sum_{v=0}^{\lfloor (M-1)/{2}\rfloor}
	(-1)^{v}u_2^{(0,2v+1)}  G_{M,1,2v}(y,x) \\
	& \quad +\frac{a_2}{a_1}\sum_{v=0}^{\lfloor (M+1)/2\rfloor-1}
	(-1)^{v}u_2^{(1,2v)}  G_{M,0,2v+1}(y,x) +\frac{a_2}{a_1}\sum_{v=0}^{\lfloor M/{2}\rfloor-1}
	(-1)^{v}u_2^{(1,2v+1)}  G_{M,1,2v+1}(y,x),\\
	&=\sum_{n=0}^{M}
	(-1)^{\lfloor\frac{n}{2}\rfloor}u_2^{(0,n)}  G_{M,\odd(n),n-\odd(n)}(y,x) +\frac{a_2}{a_1}\sum_{n=0}^{M-1}
	(-1)^{\lfloor\frac{n}{2}\rfloor}u_2^{(1,n)}  G_{M, \odd(n),n+1-\odd(n)}(y,x),
\end{align*}
i.e.:
{\footnotesize{
\be\label{u2:H1TOV1:intersect}
\begin{split}
	\sum_{(n,m)\in \ind_{M}^{1}}
	 (-1)^{\lfloor\frac{m}{2}\rfloor}\Big(\frac{a_2}{a_1}\Big)^{\odd(m)}u_2^{(\odd(m),m+n-\odd(m))}  G_{M,n,m}(y,x)=
	\sum_{(m,n)\in \ind_{M}^{1}}
	u_2^{(m,n)} (-1)^{\lfloor\frac{n}{2}\rfloor}\Big(\frac{a_2}{a_1}\Big)^m G_{M, \odd(n),z}(y,x),
\end{split}
\ee}}
where $z:=n+m-\odd(n)$.
%%%%%
Then
%%
%%
%%%%%%%%%
%%
%%
%%%%%%%%%
%%
%%
%%
%%%%%%%%%
%%
%%
%%
\eqref{U3xy:Intersect:H1} and \eqref{u2:H1TOV1:intersect} imply that:
{\footnotesize{
		 \be\label{u3xy:Final:intersect:V1}
		\begin{split}
			 &u_3(x+x_i^*,y+y_j^*)+\bo(h^{M+1})\\
			&=		\sum_{(m,n)\in \ind_{M}^{1}}
			u_2^{(m,n)} (-1)^{\lfloor\frac{n}{2}\rfloor}\Big(\frac{a_2}{a_3}\Big)^{\odd(n)} G_{M, \odd(n),z}(y,x)  +	 \frac{	 1 }{a_3} \sum_{(m,n)\in \ind_{M-2}}
			f_3^{(m,n)}H_{M,n,m}(y,x) \\ & \quad + \sum_{(n,m)\in \ind_{M}^{1}}
			\sum_{\ell=1}^{\lfloor m/2\rfloor}
			\frac{(-1)^{\ell}}{a_2} f_2^{(m-2\ell,n+2\ell-2)}  \frac{a_2^n}{a_3^n} G_{M,n,m}(y,x)
			-\sum_{m=0}^{M} \phi_3^{(m)} G_{M,0,m}(y,x)  -\sum_{m=0}^{M-1} \frac{\psi_3^{(m)}}{a_3} G_{M,1,m}(y,x),
		\end{split}
		\ee
}}
\eqref{U4xy:Intersect:H1} and \eqref{u2:H1TOV1:intersect} imply that:
{\footnotesize{
		 \be\label{u4xy:Final:intersect:V1}
		\begin{split}
			 &u_4(x+x_i^*,y+y_j^*)+\bo(h^{M+1})\\
			&=		\sum_{(m,n)\in \ind_{M}^{1}}
			u_2^{(m,n)} (-1)^{\lfloor\frac{n}{2}\rfloor}\Big(\frac{a_2}{a_1}\Big)^m  \Big(\frac{a_1}{a_4}\Big)^{\odd(n)} G_{M, \odd(n),z}(y,x) +\frac{1}{a_4}\sum_{(m,n)\in \ind_{M-2}}
				f_4^{(m,n)} H_{M,n,m}(y,x) \\
			&\quad -\sum_{(n,m)\in \ind_{M}^{1}}  \Big( \frac{\odd(m+1)}{	 (-1)^{\lfloor\frac{m}{2}\rfloor}} \phi_1^{(m+n)} +\frac{\odd(m)}{	 (-1)^{\lfloor\frac{m}{2}\rfloor} a_1}\psi_1^{(m+n-1)} \Big) \Big(\frac{a_1}{a_4}\Big)^n G_{M,n,m}(y,x)-\sum_{m=0}^{M} \phi_4^{(m)} G_{M,0,m}(y,x)
			\\
			&\quad +	\frac{1}{a_1} \sum_{(n,m)\in \ind_{M}^{1}}
			\sum_{\ell=1}^{\lfloor m/2\rfloor}
		(-1)^{\ell} f_1^{(m-2\ell,n+2\ell-2)}  \Big(\frac{a_1}{a_4}\Big)^n G_{M,n,m}(y,x)  -\frac{1}{a_4}  \sum_{m=0}^{M-1} \psi_4^{(m)} G_{M,1,m}(y,x),
		\end{split}
		\ee
}}
for $(x,y) \in (-2h,2h)^2$ and  $z:=n+m-\odd(n)$.

Now,
	by  \eqref{Intersect:u:approx:key:V} with $p=2$, \eqref{intersect:u1xy:side1:ir:thm1}, \eqref{u3xy:Final:intersect:V1}, \eqref{u4xy:Final:intersect:V1} with  $(x_{i}^*,y_{j}^*)=(x_{i},y_{j})$ and \eqref{Lu}, we have  (note that $w_1=w_2=0$ due to $(x_i^*,y_j^*)=(x_i,y_j)$, and see the third panel of \cref{fig:intersect:side1} for an illustration):
			 \be\label{LShu:proof:intersect}
			\begin{split}
				\mathcal{L}_h u  =&
				\begin{aligned}	
					 &C_{-1,1}u_1(x_i-h,y_j+h)&
					 +&C_{0,1} u_2(x_i,y_j+h)&
					 +&C_{1,1} u_2(x_i+h,y_j+h)  \\
					+&C_{-1,0} u_1(x_i-h,y_j)&
					 +&C_{0,0} u_2(x_i,y_j)&
					 +&C_{1,0} u_2(x_i+h,y_j)\\
					 +&C_{-1,-1}u_4(x_i-h,y_j-h)&
					 +&C_{0,-1} u_3(x_i,y_j-h)&
					 +&C_{1,-1}u_3(x_i+h,y_j-h)\\
				\end{aligned}\\
				=&	 \sum_{(m,n)\in \ind_{M}^{1}}
				u_2^{(m,n)}  I_{m,n}+ \sum_{(n,m)\in \ind_{M}^{1}} F^{1,0,0}_{M,m,n} +\sum_{(m,n)\in \ind_{M-2}} F^{0,0}_{M,m,n} +  \sum_{(n,m)\in \ind_{M}^{1}} \Phi^{1,0,0}_{M,m,n} \\
				&  + \sum_{m=0}^{M}  \Phi^{0,0}_{M,m}  +  \sum_{(n,m)\in \ind_{M}^{1}} \Psi^{1,0,0}_{M,m,n}+\sum_{m=0}^{M-1}  \Psi^{0,0}_{M,m}
				+\bo(h^{M+1}),\quad \mbox{as} \quad h \to 0,
			\end{split}
			\ee
where:
	{\footnotesize{
	\be\label{intersect:Imn:Cross}
	\begin{split}
		I_{m,n}:= &
			\Big(\frac{a_2}{a_1}\Big)^m    G^{1,w_1,w_2}_{M,m,n}+ G^{2,w_1,w_2}_{M,m,n} +   (-1)^{\lfloor\frac{n}{2}\rfloor} \Big(\frac{a_2}{a_3}\Big)^{\odd(n)} G^{3,w_1,w_2}_{M,z,\odd(n)}+(-1)^{\lfloor\frac{n}{2}\rfloor} \Big(\frac{a_2}{a_1}\Big)^m  \Big(\frac{a_1}{a_4}\Big)^{\odd(n)}  G^{4,w_1,w_2}_{M,z,\odd(n)},
	\end{split}
	\ee}}
	{\footnotesize{
	\be\label{intersect:FH1mn:Cross}
	\begin{split}
		F^{1,w_1,w_2}_{M,m,n} :=  \frac{1}{a_2} 	 \Big(\frac{a_2}{a_3}\Big)^n \sum_{s=1}^{\lfloor m/2\rfloor}
				(-1)^{s} f_2^{(m-2s,n+2s-2)}   G^{3,w_1,w_2}_{M,m,n} +\frac{1}{a_1}  \Big(\frac{a_1}{a_4}\Big)^n 	 \sum_{s=1}^{\lfloor m/2\rfloor}
		(-1)^{s} f_1^{(m-2s,n+2s-2)}    G^{4,w_1,w_2}_{M,m,n},
	\end{split}
	\ee}}
	\be\label{intersect:FMmn:Cross}
	\begin{split}
		F^{w_1,w_2}_{M,m,n}	:=&\sum_{p=1}^4
		\frac{f_p^{(m,n)}}{a_p} H^{p,w_1,w_2}_{M,m,n},
	\end{split}
	\ee
		 \be\label{intersect:PhiH1Mmn:Cross}
	\begin{split}
		\Phi^{1,w_1,w_2}_{M,m,n} 	:=
	 	\frac{	\phi_1^{(m+n)} \odd(m+1)}{   - (-1)^{\lfloor\frac{m}{2}\rfloor}   }  \Big(\frac{a_1}{a_4}\Big)^n  G^{4,w_1,w_2}_{M,m,n} , \quad 		 \Phi^{w_1,w_2}_{M,m}	:= -\phi_1^{(m)}  G^{1,w_1,w_2}_{M,0,m} -\sum_{p=3}^4 \phi_p^{(m)}  G^{p,w_1,w_2}_{M,m,0},
	\end{split}
	\ee
	\be\label{intersect:PsiH1Mmn:Cross}
	\begin{split}
		\Psi^{1,w_1,w_2}_{M,m,n} 	:=
	 	\frac{	\psi_1^{(m+n-1)} \odd(m)}{  -  (-1)^{\lfloor\frac{m}{2}\rfloor}    a_1}  \Big(\frac{a_1}{a_4}\Big)^n  G^{4,w_1,w_2}_{M,m,n}, \quad
	 			\Psi^{w_1,w_2}_{M,m}	 := -\frac{\psi_1^{(m)}}{a_1}  G^{1,w_1,w_2}_{M,1,m} -\sum_{p=3}^4 \frac{\psi_p^{(m)}}{a_p}  G^{p,w_1,w_2}_{M,m,1},
	\end{split}
	\ee
{\footnotesize{
		\be\label{Gmn:w1:w2}
		\begin{split}
			&	 G^{1,w_1,w_2}_{M,m,n} := \sum_{\ell=0}^{1} C_{-1,\ell}	
			 G_{M,m,n}((w_1-1)h,(w_2+\ell)h),\quad
			G^{2,w_1,w_2}_{M,m,n} :=   \sum_{k=0}^{1} \sum_{\ell=0}^{1} C_{k,\ell}
			 G_{M,m,n}((w_1+k)h,(w_2+\ell)h), \\
			&  G^{3,w_1,w_2}_{M,m,n} :=   \sum_{k=0}^{1}  C_{k,-1}  G_{M,n,m}((w_2-1)h,(w_1+k)h),\quad
			G^{4,w_1,w_2}_{M,m,n} :=  C_{-1,-1}
			 G_{M,n,m}((w_2-1)h,(w_1-1)h),
		\end{split}
		\ee}}
{\footnotesize{
		\be\label{Hmn:w1:w2}
		\begin{split}
			&	 H^{1,w_1,w_2}_{M,m,n} := \sum_{\ell=0}^{1} C_{-1,\ell}	
			 H_{M,m,n}((w_1-1)h,(w_2+\ell)h),\quad
			H^{2,w_1,w_2}_{M,m,n} :=   \sum_{k=0}^{1} \sum_{\ell=0}^{1} C_{k,\ell}
			 H_{M,m,n}((w_1+k)h,(w_2+\ell)h), \\
			&  H^{3,w_1,w_2}_{M,m,n} :=   \sum_{k=0}^{1}  C_{k,-1}  H_{M,n,m}((w_2-1)h,(w_1+k)h),\quad
			H^{4,w_1,w_2}_{M,m,n} :=  C_{-1,-1}
			 H_{M,n,m}((w_2-1)h,(w_1-1)h),
		\end{split}
		\ee}}
	and  $z:=n+m-\odd(n)$.

	Note  that $I_{m,n}$  in \eqref{intersect:Imn:Cross} satisfies
	 \be\label{intersect:Imn:Conditions:Cross}
	I_{m,n}=\bo(h^{M+1}), \qquad \mbox{ for all } \quad (m,n)\in \ind_{M}^{1}.
	\ee
Then the conditions in \eqref{intersect:Imn:Conditions:Cross}
can be equivalently rewritten as a system of linear equations for the unknowns $\{C_{k,\ell}\}_{k,\ell=-1,0,1}$.
\eqref{intersect:Imn:Conditions:Cross} has a nontrivial solution $\{C_{k,\ell}\}_{k,\ell=-1,0,1}$ if and only if $M\le 7$. Moreover, for $M=7$, up to a multiplicative constant for normalization, \eqref{Intersect:CKL:cross1} is the unique solution to
\eqref{intersect:Imn:Conditions:Cross}. Therefore, for the solution  $\{C_{k,\ell}\}_{k,\ell=-1,0,1}$ in \eqref{Intersect:CKL:cross1} with $M=7$,
we conclude that:
\be\label{trunca:order:7}
h^{-1}	\mathcal{L}_h (u-u_h)=
h^{-1} \sum_{(m,n)\in \ind_{7}^{1}}
u_2^{(m,n)}  I_{m,n}=\bo(h^{7}),\qquad h\to 0,
\ee
which proves the seventh order of consistency  of
the FDM in \cref{Intersect:thm:cross1}.	
\end{proof}
\begin{proof}[Proof of \cref{theorem:side3:w}]
	Since $(x_i,y_j)=(x_i^*+wh,y_j^*)$, we only need to	
	replace $(x_i,y_j)$ by  $(x_i^*+wh,y_j)$ on the left-hand side of \eqref{Lhu:Gamma1}. Then the rest proof is same as the proof in \cref{Intersect:thm:side1}.
\end{proof}
\begin{proof}[Proof of \cref{Intersect:thm:cross1:w1w2}]
	%%%
	Note that $(x_i-w_1h,y_j-w_2h)=(x_i^*,y_j^*)=(\xi,\zeta)$, and
 the discrete operator $\mathcal{L}_h u$ in \eqref{Lu} involves $u_p$ with $p=1,2,3,4$ (see the first panel of \cref{fig:intersect:cross:point:w}).
 	Since \cref{Intersect:thm:cross1} does not use jump conditions in \eqref{intersect:3} across the interface $\Gamma_2$ (see the first panel of \cref{transsmi:4:parts}),
if we only extend the derivation of the special case in \cref{Intersect:thm:cross1} to the general case in \cref{Intersect:thm:cross1:w1w2},  the compact 9-point scheme with a consistency order four could not be derived. In order to use jump conditions across the interface $\Gamma_2$, we derive $\tilde{u}_4$ in \eqref{U4xy:Intersect:H1:w}. More precisely, by  jump conditions across interfaces $\Gamma_2$ and $\Gamma_3$ (see the second panel of \cref{transsmi:4:parts}), we can replace all $\{ u_4^{(m,n)} \; :\; (m,n)\in \ind_M^1\}$ on the right-hand side of \eqref{Intersect:u:approx:key:V} with $p=4$ by $\{  u_2^{(m,n)} \; :\; (m,n)\in \ind^1_M \}$ from the following \eqref{U4xy:Intersect:H1:case2}--\eqref{U4xy:Intersect:H1:w}. Then the compact 9-point scheme with a consistency order four is derived in the rest of the proof.

Firstly,	%%
	%%%%%
	%%
similarly to \eqref{intersect:cross:U4:0}--\eqref{u4xy:Final:intersect:V1}, we have:
	{\footnotesize{
			 \be\label{U4xy:Intersect:H1:case2}
			\begin{split}
				 &\tilde{u}_4(x+x_i^*,y+y_j^*)+\bo(h^{M+1})\\
				&=	\sum_{(m,n)\in \ind_{M}^{1}}
				 (-1)^{\lfloor\frac{n}{2}\rfloor}\Big(\frac{a_2}{a_3}\Big)^{\odd(n)}u_2^{(m+n-\odd(n),\odd(n))} \Big(\frac{a_3}{a_4}\Big)^m G_{M,m,n}(x,y)  +\sum_{(m,n)\in \ind_{M-2}}
			\frac{	f_4^{(m,n)} }{a_4}H_{M,m,n}(x,y) \\
				&\quad-	\sum_{(m,n)\in \ind_{M}^{1}}
				\Big( (-1)^{\lfloor\frac{n}{2}\rfloor} \odd(n+1)\phi_3^{(m+n)}
				 +(-1)^{\lfloor\frac{n}{2}\rfloor} \frac{\odd(n)}{a_3}\psi_3^{(m+n-1)} \Big) \Big(\frac{a_3}{a_4}\Big)^m G_{M,m,n}(x,y)  \\
				&\quad+	\sum_{(m,n)\in \ind_{M}^{1}}
				\sum_{\ell=1}^{\lfloor n/2\rfloor}
				\frac{(-1)^{\ell}}{a_3} f_3^{(m+2\ell-2,n-2\ell)} \Big(\frac{a_3}{a_4}\Big)^m G_{M,m,n}(x,y) - \sum_{n=0}^{M}
				\phi_2^{(n)} G_{M,0,n}(x,y) -\sum_{n=0}^{M-1} \frac{\psi_2^{(n)}}{a_4} G_{M,1,n}(x,y).
			\end{split}
			\ee
	}}

		On the other hand:
	\[
	u_2^{(m+n-\odd(n),\odd(n))}
	=\begin{cases}
		u_2^{(n,0)}, &\text{if $n$ is even and } m=0,\\
		u_2^{(n+1,0)}, &\text{if $n$ is even and } m=1,\\
		u_2^{(n-1,1)}, &\text{if $n$ is odd and } m=0,\\
		u_2^{(n,1)}, &\text{if $n$ is odd and } m=1,
	\end{cases}
	\qquad \forall\;  (m,n)\in \ind_{M}^{1}.
	\]
	By \eqref{U2:intersect:H1}, for any $(m,n)\in \ind_{M}^{1}$ we have:
	\[
	\begin{split}
		 &u_2^{(m+n-\odd(n),\odd(n))}=\begin{cases}
			 (-1)^{\lfloor\frac{n}{2}\rfloor}
			u_2^{(0,n)}+
			\sum_{\ell=1}^{\lfloor \frac{n}{2} \rfloor}
			\frac{(-1)^{\ell}}{a_2} f_2^{(n-2\ell,2\ell-2)}, &\text{if $n$ is even and } m=0,\\
			 %%%%%%%%%%%%%%%%%%%%%%%%%%%%%
			 (-1)^{\lfloor\frac{n+1}{2}\rfloor}
			u_2^{(1,n)}+
			\sum_{\ell=1}^{\lfloor \frac{n+1}{2}\rfloor}
			\frac{(-1)^{\ell}}{a_2} f_2^{(n+1-2\ell,2\ell-2)}, &\text{if $n$ is even and } m=1,\\
			 %%%%%%%%%%%%%%%%%%%%%%%%%%%
			 (-1)^{\lfloor\frac{n-1}{2}\rfloor}
			u_2^{(0,n)}+
			\sum_{\ell=1}^{\lfloor \frac{n-1}{2}\rfloor}
			\frac{(-1)^{\ell}}{a_2} f_2^{(n-1-2\ell,1+2\ell-2)}, &\text{if $n$ is odd and } m=0,\\
			%%%%%%%%%%%%%%%%%%%%%%%%%
			 (-1)^{\lfloor\frac{n}{2}\rfloor}
			u_2^{(1,n)}+
			\sum_{\ell=1}^{\lfloor \frac{n}{2} \rfloor}
			\frac{(-1)^{\ell}}{a_2} f_2^{(n-2\ell,1+2\ell-2)}, &\text{if $n$ is odd and } m=1.
		\end{cases}
	\end{split}
	\]
	Now we observe that:
	\be\label{u2V1:w1:w2}
	\begin{split}
		u_2^{(m+n-\odd(n),\odd(n))}=
		(-1)^{\lfloor \frac{q}{2} \rfloor}
		u_2^{(m,n)} + 	\frac{1}{a_2}
		\sum_{\ell=1}^{\lfloor \frac{q}{2} \rfloor}
		(-1)^{\ell} f_2^{(q-2\ell,\odd(n)+2\ell-2)},
		\quad \forall\;  (m,n)\in \ind_{M}^{1},
	\end{split}
	\ee
	with $q:=n-(-1)^{m}\odd(n+m)$.
	So	
		\eqref{U4xy:Intersect:H1:case2} and \eqref{u2V1:w1:w2} imply that:
		{\footnotesize{
				 \be\label{U4xy:Intersect:H1:w}
				\begin{split}
					 &\tilde{u}_4(x+x_i^*,y+y_j^*)+\bo(h^{M+1})\\
					&=	\sum_{(m,n)\in \ind_{M}^{1}}
					 (-1)^{\lfloor\frac{n}{2}\rfloor+\lfloor \frac{q}{2} \rfloor}\Big(\frac{a_2}{a_3}\Big)^{\odd(n)} \Big(\frac{a_3}{a_4}\Big)^m 	
					u_2^{(m,n)} G_{M,m,n}(x,y)  +  	 \frac{1 }{a_4} \sum_{(m,n)\in \ind_{M-2}}
				f_4^{(m,n)} H_{M,m,n}(x,y) \\
					&\quad+	 \frac{1}{a_2} \sum_{(m,n)\in \ind_{M}^{1}}
					 \Big(\frac{a_2}{a_3}\Big)^{\odd(n)} \Big(\frac{a_3}{a_4}\Big)^m 	 \sum_{\ell=1}^{\lfloor \frac{q}{2} \rfloor}
					 (-1)^{\ell+\lfloor\frac{n}{2}\rfloor} f_2^{(q-2\ell,\odd(n)+2\ell-2)} G_{M,m,n}(x,y)  \\
					&\quad-	 \sum_{(m,n)\in \ind_{M}^{1}}
					\Big( \frac{ \odd(n+1)}{(-1)^{\lfloor\frac{n}{2}\rfloor}}\phi_3^{(m+n)}
					 + \frac{\odd(n)}{a_3 (-1)^{\lfloor\frac{n}{2}\rfloor}}\psi_3^{(m+n-1)} \Big) \Big(\frac{a_3}{a_4}\Big)^m G_{M,m,n}(x,y)  - \sum_{m=0}^{M}
					 \phi_2^{(m)} G_{M,0,m}(x,y) \\
					&\quad+	 \sum_{(m,n)\in \ind_{M}^{1}}
					 \sum_{\ell=1}^{\lfloor \frac{n}{2}\rfloor}
					 \frac{(-1)^{\ell}}{a_3} f_3^{(m+2\ell-2,n-2\ell)} \Big(\frac{a_3}{a_4}\Big)^m G_{M,m,n}(x,y) -\sum_{m=0}^{M-1} \frac{\psi_2^{(m)}}{a_4} G_{M,1,m}(x,y),
				\end{split}
				\ee
		}}
		with $q:=n-(-1)^{m}\odd(n+m)$.

	Now,
	by  \eqref{Lu}, \eqref{Intersect:u:approx:key:V} with $p=2$, \eqref{intersect:u1xy:side1:ir:thm1}, \eqref{u3xy:Final:intersect:V1}, \eqref{u4xy:Final:intersect:V1} and \eqref{U4xy:Intersect:H1:w} with
	$(x_i,y_j)=(x_i^*+w_1h,y_j^*+w_2h)$,
	we have that:
	{\small{
			\[
			\begin{split}
				\mathcal{L}_h u =&\sum_{\ell=0}^{1} {c}_{-1,\ell}u_1(x_i-h,y_j+\ell h)  +\sum_{k=0}^{1}\sum_{\ell=0}^{1} {c}_{k,\ell} u_2(x_i+kh,y_j+\ell h) +\sum_{k=0}^{1} {c}_{k,-1}u_3(x_i+kh,y_j- h) \\
				&+   \big({c}_{-1,-1} u_4(x_i-h,y_j- h) +\tilde{c}_{-1,-1}\tilde{u}_4(x_i-h,y_j- h)\big) =	 \sum_{(m,n)\in \ind_{M}^{1}}
				u_2^{(m,n)}  \tilde{I}_{m,n}+ \sum_{(n,m)\in \ind_{M}^{1}} \tilde{F}^{1,w_1,w_2}_{M,m,n}
				\\
				 %%%%%%%%%%%%%%%%%%%%%%%%%%%%%%%
				&
				+ \sum_{(m,n)\in \ind_{M-2}} \tilde{F}^{w_1,w_2}_{M,m,n}
				+   \sum_{(n,m)\in \ind_{M}^{1}} \tilde{\Phi}^{1,w_1,w_2}_{M,m,n}+ \sum_{m=0}^{M}  \tilde{\Phi}^{w_1,w_2}_{M,m}
				+  \sum_{(n,m)\in \ind_{M}^{1}} \tilde{\Psi}^{1,w_1,w_2}_{M,m,n}+ \sum_{m=0}^{M-1}  \tilde{\Psi}^{w_1,w_2}_{M,m}+\bo(h^{M+1}),
			\end{split}
			\]
	}}
	as 	$h \to 0$, where:
	{\small{
			\be\label{Imn:Final:one}
			\begin{split}
				 \tilde{I}_{m,n}:=\hat{I}_{m,n}+(-1)^{\lfloor\frac{n}{2}\rfloor +\lfloor \frac{q}{2} \rfloor }
				 \Big(\frac{a_2}{a_3}\Big)^{\odd(n)} \Big(\frac{a_3}{a_4}\Big)^m 	
				 \tilde{G}^{4,w_1,w_2}_{M,n,m},\quad 	 \tilde{F}^{w_1,w_2}_{M,m,n}	:=	 \hat{F}^{w_1,w_2}_{M,m,n}+
				 \frac{f_4^{(m,n)}}{a_4}  \tilde{H}^{4,w_1,w_2}_{M,m,n},
			\end{split}
			\ee}}
	{\footnotesize{
	 \be\label{intersect:FH1mn:Cross:tilde}
\begin{split}
	\tilde{F}^{1,w_1,w_2}_{M,m,n} :=\hat{F}^{1,w_1,w_2}_{M,m,n} +   \Big(\frac{a_3}{a_4}\Big)^m  \Big[ \Big(\frac{a_2}{a_3}\Big)^{\odd(n)}
	 	 \sum_{s=1}^{\lfloor \frac{q}{2}\rfloor}
	 \frac{ (-1)^{\lfloor\frac{n}{2}\rfloor} }{ a_2(-1)^{s}  }f_2^{(q-2s,\odd(n)+2s-2)}+ 			 \sum_{s=1}^{\lfloor \frac{n}{2} \rfloor}
	\frac{(-1)^{s}}{a_3} f_3^{(m+2s-2,n-2s)} \Big] \tilde{G}^{4,w_1,w_2}_{M,m,n},
\end{split}
\ee}}
		 \be\label{intersect:PhiH1Mmn:Cross:tilde}
		\begin{split}
			 \tilde{\Phi}^{1,w_1,w_2}_{M,m,n} 	 := \hat{\Phi}^{1,w_1,w_2}_{M,m,n}+
			\frac{	\phi_3^{(m+n)} \odd(n+1)}{   - (-1)^{\lfloor\frac{n}{2}\rfloor}   }    \Big(\frac{a_3}{a_4}\Big)^m \tilde{G}^{4,w_1,w_2}_{M,m,n},\quad 	 \tilde{\Phi}^{w_1,w_2}_{M,m}	:=	 \hat{\Phi}^{w_1,w_2}_{M,m} - \phi_2^{(m)} \tilde{G}^{4,w_1,w_2}_{M,m,0},
		\end{split}
		\ee
\be\label{intersect:PsiMm:Cross:tilde}
\begin{split}
				 \tilde{\Psi}^{1,w_1,w_2}_{M,m,n} 	 :=	 \hat{\Psi}^{1,w_1,w_2}_{M,m,n}+
	\frac{	\psi_3^{(m+n-1)} \odd(n)}{  -  (-1)^{\lfloor\frac{n}{2}\rfloor}    a_3}   \Big(\frac{a_3}{a_4}\Big)^m \tilde{G}^{4,w_1,w_2}_{M,m,n},\quad
	\tilde{\Psi}^{w_1,w_2}_{M,m}	:= 	 \hat{\Psi}^{w_1,w_2}_{M,m}- \frac{\psi_2^{(m)}}{a_4} \tilde{G}^{4,w_1,w_2}_{M,m,1},
\end{split}
\ee
{\small{
		\be\label{tilde:Gmn:w1:w2}
		\begin{split}
			 \tilde{G}^{4,w_1,w_2}_{M,m,n} :=  \tilde{c}_{-1,-1}
			 G_{M,n,m}((w_1-1)h,(w_2-1)h),\quad \tilde{H}^{4,w_1,w_2}_{M,m,n} :=  \tilde{c}_{-1,-1}
			 H_{M,m,n}((w_1-1)h,(w_2-1)h),
		\end{split}
		\ee}}
	every $c_{k,\ell}, \tilde{c}_{-1,-1}\in \R$, $q:=n-(-1)^{m}\odd(n+m)$,  $\hat{I}_{m,n}$, $\hat{F}^{1,w_1,w_2}_{M,m,n}$, $\hat{F}^{w_1,w_2}_{M,m,n}$, $\hat{\Phi}^{1,w_1,w_2}_{M,m,n}$, $\hat{\Phi}^{w_1,w_2}_{M,m}$, $\hat{\Psi}^{1,w_1,w_2}_{M,m,n}$, $\hat{\Psi}^{w_1,w_2}_{M,m}$ are obtained by replacing $C_{k,\ell}$ by $c_{k,\ell}$ in
	 \eqref{intersect:Imn:Cross}--\eqref{intersect:PsiH1Mmn:Cross} for $k,\ell=-1,0,1$.
	%
	%%
%
%%
%%
%%
%%
%%
%%
%%
%%
%%%%%
%%
%%
 We consider
 \be\label{sum:Ckl}
{C}_{k,\ell}
:=\begin{cases}
	c_{-1,-1}+\tilde{c}_{-1,-1}, &\text{if}\quad  k=\ell=-1,\\
	{c}_{0,0}:=1, &\text{if}\quad  k=\ell=0,\\
	{c}_{k,\ell}, &\text{else},
\end{cases}
\ee
 and
	\be\label{solve:Imn:last}
	 \tilde{I}_{m,n}=\bo(h^{M+1}), \qquad \mbox{ for all } \quad (m,n)\in \ind_{M}^{1}.
	\ee
	where $ \tilde{I}_{m,n}$ is defined in \eqref{Imn:Final:one}.
Then the conditions in \eqref{solve:Imn:last}
can be equivalently rewritten as a system of linear equations on the unknowns $\{c_{k,\ell}\}_{k,\ell=-1,0,1} \cup \{\tilde{c}_{-1,-1}\}$. By calculation, we observe that
\eqref{solve:Imn:last} has a nontrivial solution $\{C_{k,\ell}\}_{k,\ell=-1,0,1}$ defined in \eqref{sum:Ckl} if and only if $M\le 4$. Furthermore,
	we observe that
	$\{{C}_{k,\ell} \}_{k,\ell=-1,0,1}$ defined in \eqref{sum:Ckl} is uniquely determined by solving \eqref{solve:Imn:last} with $M=4$.
	The rest of the proof is similar as the proof of \cref{Intersect:thm:cross1}.
\end{proof}
\begin{proof}[Proof of \cref{intersect:thm:convergence}]
	Clearly, all the $\{C_{k,\ell}\}_{k,\ell=-1,0,1}$ in \cref{Intersect:thm:regular,Intersect:thm:side1,Intersect:thm:cross1} satisfy the sign condition \eqref{intersect:sign:condition},
	and the summation condition \eqref{intersect:sum:condition}.
	For simplicity, we assume $\Omega:=(0,1)^2$ and $h:=1/N$ with $N\in \N$. We define $\Omega_h:=\Omega \cap (h\Z^2)$, $\partial \Omega_h:=\partial \Omega \cap (h\Z^2)$, $\overline{\Omega}_h:=\overline{\Omega} \cap (h\Z^2)$, and $(x_i,y_j):=(ih,jh)$. So  $\overline{\Omega}_h:=\{(x_i,y_j):\ 0\le i,j \le N\}$ and we also define $V(\overline{\Omega}_h):=\{(v)_{i,j}:\ 0\le i,j \le N\}$ with $(v)_{i,j}\in \mathbb{R}$,
	and for any $v\in V(\overline{\Omega}_h)$,  $(v)_{i,j}$ represents the value of $v$ at the point $(x_i,y_j)$.	
	Recall that $\Gamma:=\Gamma_1 \cup \Gamma_2 \cup \Gamma_3 \cup \Gamma_4\cup\{(\xi,\zeta)\}$,
	so we define that:
	\be\label{Deltahuhij}
	(\Delta_h u_h)_{i,j}:=
	\begin{cases}
		-h^{-2}\mathcal{L}_h u_h, &\text{if } (x_i,y_j)\in \Omega_h \mbox{ and } (x_i,y_j)\in \Omega \setminus \Gamma,\\
		-h^{-1}\mathcal{L}_h u_h, &\text{if } (x_i,y_j)\in \Omega_h \mbox{ and } (x_i,y_j)\in \Gamma,
	\end{cases}
	\ee
	where $\mathcal{L}_h$ is defined in \cref{Intersect:thm:regular} for $(x_i,y_j)\in \Omega \setminus \Gamma$, and
	$\mathcal{L}_h$ is defined in \cref{Intersect:thm:side1,Intersect:thm:cross1} for $ (x_i,y_j)\in \Gamma$.
	Therefore,  using FDMs in \cref{Intersect:thm:regular,Intersect:thm:side1,Intersect:thm:cross1}, we find $u_h \in
	V(\overline{\Omega}_h)$ satisfying:
	\be\label{DeltahuhF}
	\Delta_h u_h=\tilde{F}:=
	\begin{cases}
		\tfrac{-1}{a(x_i,y_j)}F, &\text{if } (x_i,y_j)\in \Omega_h \mbox{ and } (x_i,y_j)\in \Omega \setminus \Gamma,\\
		\tfrac{-1}{h}F, &\text{if } (x_i,y_j)\in \Omega_h \mbox{ and } (x_i,y_j)\in \Gamma,
	\end{cases} \quad \mbox{on}\quad \Omega_h \quad\mbox{with}\quad u_h=g\quad \mbox{on}\quad
	\partial \Omega_h,
	\ee
	where $F$ is the right-hand side of  FDM in \cref{Intersect:thm:regular,Intersect:thm:side1,Intersect:thm:cross1}.

Using \eqref{intersect:sign:condition} and \eqref{intersect:sum:condition},
we now prove the discrete maximum principle:
for any $v\in V(\overline{\Omega}_h)$ satisfying $\Delta_h v\ge 0$ on $\Omega_h$, we must have $\max_{(x_i,y_j)\in \Omega_h} v(x_i,y_j)\le \max_{(x_i,y_j)\in \partial \Omega_h} v(x_i,y_j)$, where $\Delta_h$ is defined in \eqref{Deltahuhij}.

%%%%%%%%%%%%%%
Suppose that $\max\limits_{(x_i,y_j)\in \Omega_h} (v)_{i,j}> \max\limits_{(x_i,y_j)\in \partial \Omega_h} (v)_{i,j}$.
Take $(x_m,y_n)\in \Omega_h$ where $v$ achieves its maximum in $\Omega_h$. Because all the stencils satisfying \eqref{intersect:sign:condition} and \eqref{intersect:sum:condition}, we have:
\[
\sum_{k,\ell\in\{-1,0,1\} \atop k\ne 0, \  \ell\ne 0} -C_{k,\ell} (v)_{m+k,n+\ell} \le C_{0,0} (v)_{m,n}.
\]
By
\[
0\le h^s(\Delta_h v)_{m,n}=-C_{0,0} (v)_{m,n}-\sum_{k,\ell\in\{-1,0,1\} \atop k\ne 0, \  \ell\ne 0} C_{k,\ell} (v)_{m+k,n+\ell},
\]
where $s=1,2$,	we have
\[
C_{0,0} (v)_{m,n}\le \sum_{k,\ell\in\{-1,0,1\} \atop k\ne 0, \  \ell\ne 0} -C_{k,\ell} (v)_{m+k,n+\ell} \le C_{0,0} (v)_{m,n}.
\]
Thus, equality holds throughout and $v$ achieves its maximum at all its nearest neighbors of $(x_m,y_n)$. Applying the same argument to the neighbors in $\Omega_h$ and repeat this argument, we conclude that $v$ must be a constant contradicting our assumption.
This proves the discrete maximum principle.
%%%%%%%%%%%%%%%%%	

	%
	%
	%
	Let $U_h:=\{ u(x_i,y_j)\}_{(x_i,y_j)\in \overline{\Omega}_h}$. By \cref{Intersect:thm:regular,Intersect:thm:side1,Intersect:thm:cross1}, we have:
	$\Delta_h U_h=\tilde{F}+R$, where:
	\be\label{Rxiyj}
	|R(x_i,y_j)|\le
	\begin{cases}
		C h^6, &\text{if } (x_i,y_j)\in \Omega_h \mbox{ and } (x_i,y_j)\in \Omega \setminus \Gamma,\\
		C h^7, &\text{if } (x_i,y_j)\in \Omega_h \mbox{ and } (x_i,y_j)\in  \Gamma,
	\end{cases}
	\ee
	where $C>0$ is independent of $h$.
	Define $E_h:=U_h-u_h$ on $\overline{\Omega}_h$. By \eqref{DeltahuhF} and $\Delta_h U_h=\tilde{F}+R$,
	\be\label{intersect:trunR}
	\Delta_h E_h=\Delta_h U_h-\Delta_h u_h=R \quad \mbox{on}\quad \Omega_h \quad \mbox{with}\quad E_h=0\quad \mbox{on}\quad \partial \Omega_h.
	\ee
		By \eqref{intersect:sign:condition}, \eqref{intersect:sum:condition} and \eqref{Deltahuhij}, we have:
	\be\label{Deltah1}
	\Delta_h 1=0 \quad \mbox{on}\quad \Omega_h.
	\ee
	We define the comparison function $\theta:=\frac{1}{24}(x-\xi)^2+\frac{1}{24}(y-\zeta)^2$ and $\Theta:=\{\theta(x_i,y_j):\ 0\le i,j \le N\}$. For $0<\xi,\zeta<1$, we have  $0\le \Theta \le \frac{1}{12}$ on $[0,1]^2$. By $\{C_{k,\ell}\}_{k,\ell=-1,0,1}$ in \cref{Intersect:thm:regular,Intersect:thm:side1,Intersect:thm:cross1},  \eqref{Luh}, and \eqref{Lu}, we have
\be\label{Lh:Theta:ij:order6}
(\mathcal{L}_h \Theta)_{i,j}=	 \sum_{k=-1}^1\sum_{\ell=-1}^1  C_{k,\ell}\theta(x_i+kh,y_j+\ell h)	=
\begin{cases}
	-h^{2}, &\text{if } (x_i,y_j)\in \Omega_h \mbox{ and } (x_i,y_j)\in \Omega \setminus \Gamma,\\
-h^2\frac{a_1+a_2}{2a_2}, &\text{if } (x_i,y_j)\in \Omega_h \mbox{ and } (x_i,y_j)\in  \Gamma_1,\\
-h^2\frac{a_4+a_3}{2a_3}, &\text{if } (x_i,y_j)\in \Omega_h \mbox{ and } (x_i,y_j)\in  \Gamma_2,\\
-h^2\frac{a_2+a_3}{2a_2}, &\text{if } (x_i,y_j)\in \Omega_h \mbox{ and } (x_i,y_j)\in  \Gamma_3,\\
-h^2\frac{a_1+a_4}{2a_1}, &\text{if } (x_i,y_j)\in \Omega_h \mbox{ and } (x_i,y_j)\in  \Gamma_4,\\
-h^2\frac{(a_1+a_2)(a_2+a_3)}{4a_2^2}, &\text{if } (x_i,y_j)\in \Omega_h \mbox{ and } (x_i,y_j)=(\xi,\zeta).
\end{cases}
\ee
 \eqref{Deltahuhij}	leads to
		{\small{
			\be\label{order6:proof:Theta}
			(\Delta_h \Theta)_{i,j}=
			\begin{cases}
				-h^{-2}(\mathcal{L}_h \Theta)_{i,j}=1, &\text{if } (x_i,y_j)\in \Omega_h \mbox{ and } (x_i,y_j)\in \Omega \setminus \Gamma,\\
				-h^{-1}(\mathcal{L}_h \Theta)_{i,j}=\frac{h}{2}(1+\frac{a_1}{a_2})> \frac{h}{2}, &\text{if } (x_i,y_j)\in \Omega_h \mbox{ and } (x_i,y_j)\in  \Gamma_1,\\
				-h^{-1}(\mathcal{L}_h \Theta)_{i,j}=\frac{h}{2}(1+\frac{a_4}{a_3})> \frac{h}{2}, &\text{if } (x_i,y_j)\in \Omega_h \mbox{ and } (x_i,y_j)\in  \Gamma_2,\\
				-h^{-1}(\mathcal{L}_h \Theta)_{i,j}=\frac{h}{2}(1+\frac{a_3}{a_2})> \frac{h}{2}, &\text{if } (x_i,y_j)\in \Omega_h \mbox{ and } (x_i,y_j)\in  \Gamma_3,\\
			   -h^{-1}(\mathcal{L}_h \Theta)_{i,j}=\frac{h}{2}(1+\frac{a_4}{a_1})> \frac{h}{2}, &\text{if } (x_i,y_j)\in \Omega_h \mbox{ and } (x_i,y_j)\in  \Gamma_4,\\
			   -h^{-1}(\mathcal{L}_h \Theta)_{i,j}=\frac{h}{4}(1+\frac{a_1}{a_2}+\frac{a_3}{a_2}+\frac{a_1a_3}{a_2^2})> \frac{h}{4}, &\text{if } (x_i,y_j)\in \Omega_h \mbox{ and } (x_i,y_j)=(\xi,\zeta).					 
			\end{cases}
			\ee}}
From $C>0$ in \eqref{Rxiyj}, we observe that
	\be\label{4CH6}
	4Ch^6(\Delta_h \Theta)_{i,j}\ge 	
	\begin{cases}
		C h^6, &\text{if } (x_i,y_j)\in \Omega_h \mbox{ and } (x_i,y_j)\in \Omega \setminus \Gamma,\\
		C h^7, &\text{if } (x_i,y_j)\in \Omega_h \mbox{ and } (x_i,y_j)\in  \Gamma.
	\end{cases}
	\ee
Note that $E_h:=U_h-u_h$. We deduce that
 \eqref{intersect:trunR} implies
	\be\label{order6:proof:compare}
	(\Delta_h (E_h+4C h^6 \Theta))_{i,j}
	=(\Delta_h E_h)_{i,j}+4C h^6 (\Delta_h \Theta)_{i,j}
= R(x_i,y_j)+4C h^6 (\Delta_h \Theta)_{i,j}, \mbox{ for } (x_i,y_j)\in \Omega_h.
	\ee
		By \eqref{Rxiyj}, \eqref{4CH6}, and \eqref{order6:proof:compare},
		\[
			(\Delta_h (E_h+4C h^6 \Theta))_{i,j}\ge 0 \quad \mbox{for} \quad (x_i,y_j)\in \Omega_h.
		\]
	By the discrete maximum principle of $\Delta_h$ on $\Omega_h$, $C> 0$  in \eqref{Rxiyj},  $E_h=0$ on $\partial \Omega_h$, and  $0\le \Theta \le \frac{1}{12}$ on  $\overline{\Omega}$, we obtain that
	\be\label{use:max:order6:inequ}
	\begin{split}
		\max_{(x_i,y_j)\in \Omega_h} (E_h)_{i,j}
		&\le
		\max_{(x_i,y_j)\in \Omega_h} (E_h+4C h^6 \Theta)_{i,j} \le
		\max_{(x_i,y_j)\in \partial \Omega_h} (E_h+4C h^6 \Theta)_{i,j}\\
		&\le \max_{(x_i,y_j)\in \partial \Omega_h} (E_h)_{i,j}+
		4C h^6 \times \max_{(x_i,y_j)\in \partial \Omega_h}(\Theta)_{i,j}
		=\frac{C}{3}  h^6.
	\end{split}
	\ee
	A similar argument can be applied to $-E_h$. Hence, $\|E_h\|_\infty \le \frac{C}{3} h^6$.
	Thus \eqref{intersect:order6:formula} is proved. Finally, the sign condition \eqref{intersect:sign:condition}, the summation condition \eqref{intersect:sum:condition}, and the Dirichlet boundary condition of \eqref{intersect:3} together imply the M-matrix property.
\end{proof}

%%%%%%%%%%%%\noindent \emph{The right hands in \cref{Intersect:thm:cross1:w1w2:order2}.}
%
%
%
%
\begin{proof}[Proof of  \cref{theorem:side1:w:order3}]
	The derivation \eqref{max:sign:Ckl:side} and \eqref{LhGamma1uh:order3} is straightforward by the proof of \cref{theorem:side3:w} with $M=3$. The summation condition \eqref{intersect:sum:condition} for $\{C_{k,\ell}\}_{k,\ell=-1,0,1}$ in \eqref{max:sign:Ckl:side} for any $\rho \in \R$ can be verified easily.

	For the $\{C_{k,\ell}\}_{k,\ell=-1,0,1}$ in \eqref{max:sign:Ckl:side},
	we can check that all $r_p$ in \eqref{Ckl:side:splus} satisfy $r_p>0$ for $p=1,2,\dots,5$ and $w\in(0,1)$, all $s_p$ in \eqref{Ckl:side:splus} satisfy $s_p<0$ for $p=1,2,\dots,8$ and $w\in(0,1)$.
	So
	\eqref{max:sign:Ckl:side} satisfies the sign condition \eqref{intersect:sign:condition},
	if and only if $\rho\in \R$ satisfies:
	\be\label{rho:proof}
	 \max\{\underline{b}_1,\underline{b}_2,\underline{b}_3\} \le \rho \le \min\{\overline{b}_1,\overline{b}_2\},
	\ee
	which is just
\[
\underline{b}_1:=-\frac{t_3+t_4\alpha^2+t_5\alpha}{s_1+s_2\alpha^2+s_3\alpha}, \quad \underline{b}_2:=-\frac{t_6\alpha^2+t_7\alpha}{s_4\alpha^2+s_5\alpha}, \quad \underline{b}_3:=-\frac{t_8+t_9\alpha^2+t_{10}\alpha}{s_6+s_7\alpha^2+s_8\alpha},
\]
	\[
	\overline{b}_1:=0, \qquad  \overline{b}_2:=-\frac{t_1\alpha^2+t_2\alpha}{r_4\alpha^2+r_5\alpha}.
	\]
We now show that the interval in \eqref{rho:range} (i.e., the interval in \eqref{rho:proof}) is nonempty.	
In particular,  by a direct calculation, we obtain that:
	\[
	\begin{cases}
		 \max\{\underline{b}_1,\underline{b}_2,\underline{b}_3\}<-0.2\quad \mbox{and} \quad -0.018<\min\{\overline{b}_1,\overline{b}_2\}, &\text{if } \alpha \in (0,1] \text{ and } w\in(0,1/2],\\	
		 \max\{\underline{b}_1,\underline{b}_2,\underline{b}_3\}<-0.04 \quad \mbox{and} \quad  \min\{\overline{b}_1,\overline{b}_2\}= 0, &\text{if } \alpha \in (0,1] \text{ and } w\in[1/2,1),\\
		 \max\{\underline{b}_1,\underline{b}_2,\underline{b}_3\}<0 \quad \mbox{and} \quad  \min\{\overline{b}_1,\overline{b}_2\}= 0,  &\text{if } \alpha \in [1,+\infty) \text{ and } w\in(0,1/2],\\
		 \max\{\underline{b}_1,\underline{b}_2,\underline{b}_3\}<-0.1 \quad \mbox{and} \quad-0.012<\min\{\overline{b}_1,\overline{b}_2\}, &\text{if } \alpha \in [1,+\infty) \text{ and } w\in[1/2,1).

	\end{cases}
	\]
	Thus, the interval in \eqref{rho:range} is nonempty and	there exists $\rho\in \R$ such that $\{C_{k,\ell}\}_{k,\ell=-1,0,1}$ in \eqref{max:sign:Ckl:side}  satisfies the sign condition \eqref{intersect:sign:condition} for any positive $a_1,a_2$ and $w\in(0,1)$.
\end{proof}

\begin{proof}[Proof of  \cref{Intersect:thm:cross1:w1w2:order2}]
	The derivation of $\mathcal{L}_h u_h=F$ is straightforward by the proof of \cref{Intersect:thm:cross1} with $M=2$ and $(x_i,y_j)=(x_i^*+ w_1h,y_j^*+ w_2h)$.
The summation condition \eqref{intersect:sum:condition} for $\{C_{k,\ell}\}_{k,\ell=-1,0,1}$ in \eqref{max:sign:Ckl:intersect} can be verified easily.

	For the $\{C_{k,\ell}\}_{k,\ell=-1,0,1}$ in \eqref{max:sign:Ckl:intersect},
	we can check that all $r_p$ in \eqref{Ckl:cross:tplus} satisfy $r_p>0$ for $p=1,2,\dots,12$ and $(w_1,w_2)\in(0,1)^2$, all $s_p$ in \eqref{Ckl:cross:tplus} satisfy $s_p<0$ for $p=1,2,\dots,4$ and $(w_1,w_2)\in(0,1)^2$. Thus, $\{C_{k,\ell}\}_{k,\ell=-1,0,1}$ in \eqref{max:sign:Ckl:intersect}  satisfies the sign condition \eqref{intersect:sign:condition}  for any positive $a_1,a_2,a_3,a_4$ and $(w_1,w_2)\in(0,1)^2$.
\end{proof}	
\begin{proof}[Proof of \cref{intersect:thm:convergence:order3}]
	Recall that $\Gamma:=\Gamma_1 \cup \Gamma_2 \cup \Gamma_3 \cup \Gamma_4\cup\{(\xi,\zeta)\}$.	 Let
	$ \Omega_{h,\Gamma}:=\Omega_{h,\Gamma_1,\Gamma_2} \cup \Omega_{h,\Gamma_3,\Gamma_4} \cup \Omega_{h,\xi,\zeta}$,
	$\Omega_{h,\Gamma_1,\Gamma_2}:=\{ (x_i,y_j): (x_i\pm wh,y_j)=(x_i^*,y_j^*)\in \Gamma_1\cup \Gamma_2,  0<w<1\}$,
	$\Omega_{h,\Gamma_3,\Gamma_4}:=\{ (x_i,y_j): (x_i,y_j\pm wh)=(x_i^*,y_j^*)\in \Gamma_3\cup \Gamma_4,  0<w<1\}$,
	$\Omega_{h,\xi,\zeta}:=\{ (x_i,y_j): (x_i\pm w_1h,y_j\pm w_2h)=(\xi,\zeta), 0<w_1,w_2<1 \}$.
	Then
	 \eqref{intersect:order6:formula:order3} is obtained similarly by the proof of \cref{intersect:thm:convergence} by the following replacements:\\
	Replace \eqref{Deltahuhij} by
	\[
	(\Delta_h u_h)_{i,j}:=
	\begin{cases}
		-h^{-2}\mathcal{L}_h u_h, &\text{if } (x_i,y_j)\in \Omega_h \setminus \Omega_{h,\Gamma},\\
		-h^{-1}\mathcal{L}_h u_h, &\text{if } (x_i,y_j)\in \Omega_{h,\Gamma_1,\Gamma_2} \cup \Omega_{h,\Gamma_3,\Gamma_4},\\
		-h^{-1}\mathcal{L}_h u_h, &\text{if }  (x_i,y_j)\in \Omega_{h,\xi,\zeta},
	\end{cases}
	\]
	where $\mathcal{L}_h$ is defined in \cref{Intersect:thm:regular} for $(x_i,y_j)\in \Omega_h \setminus \Omega_{h,\Gamma}$, $\mathcal{L}_h$ is defined in \cref{theorem:side1:w:order3} for  $(x_i,y_j)\in \Omega_{h,\Gamma_1,\Gamma_2} \cup \Omega_{h,\Gamma_3,\Gamma_4}$, and
	$\mathcal{L}_h$ is defined in \cref{Intersect:thm:cross1:w1w2:order2} for $(x_i,y_j)\in \Omega_{h,\xi,\zeta}$.\\
	Replace \eqref{Rxiyj} by:
	\be\label{Trun:General}
	|R(x_i,y_j)|\le
	\begin{cases}
		C h^3, &\text{if } (x_i,y_j)\in \Omega_h \setminus \Omega_{h,\Gamma},\\
		C h^2, &\text{if } (x_i,y_j)\in \Omega_{h,\Gamma_1,\Gamma_2} \cup \Omega_{h,\Gamma_3,\Gamma_4},\\
		C h^2, &\text{if } (x_i,y_j)\in \Omega_{h,\xi,\zeta}.
	\end{cases}
	\ee
We define the comparison function $\theta_{p}:=\theta \chi_{\Omega_{p}}$ for $p=1,2,3,4$, and choose $\theta_{1}=\frac{1}{24}\big((x-\xi)^2+(y-\zeta)^2\big)+10$, $\theta_{2}=\frac{1}{24}\big((x-\xi)^2+(y-\zeta)^2\big)+8$, $\theta_{3}=\frac{1}{24}\big((x-\xi)^2+(y-\zeta)^2\big)+9$, $\theta_{4}=\frac{1}{24}\big((x-\xi)^2+(y-\zeta)^2\big)+15$.
	For $0<\xi,\zeta<1$, we have  $0\le \Theta \le \frac{1}{12}+15=\frac{181}{12}$ on $[0,1]^2$.
Similar as in \eqref{Lh:Theta:ij:order6} and \eqref{order6:proof:Theta}, by $\{C_{k,\ell}\}_{k,\ell=-1,0,1}$  in \cref{theorem:side1:w:order3}, we have
\be\label{interface:Deltah:Theta:ij}
\begin{split}
(\Delta_h \Theta)_{i,j}=	-h^{-1}(\mathcal{L}_h \Theta)_{i,j}  =h^{-1} C_1  + hC_2,
\end{split}
\ee
where
{\small{
\be\label{Thm34:proof:C1C2}
\begin{split}
& C_1=4\frac{       -6\big((w+1)\alpha^2+(1-w)\alpha \big) \rho+ 3w\alpha^2-3(w-1)\alpha }{r_2\alpha^2+r_3\alpha+r_1 }, \quad C_2= \frac{    (p_1\alpha^2+p_2\alpha+p_3)   \rho+ q_1\alpha^2+q_2\alpha+q_3  }{2(r_2\alpha^2+r_3\alpha+r_1) },
\end{split}
\ee}}
%%%%
\[
\begin{split}
&p_1 = -4w^3-2w^2-2,\quad p_2 = 8w^3-4w^2-4,\quad p_3 = -4w^3+6w^2-2,\\
&q_1 = 2w^3-w^2+w,\quad q_2 = -4w^3+4w^2-w+1,\quad q_3 = 2w^3-3w^2+1,
\end{split}
\]
%%%
%%%
and $r_1,r_2,r_3$ are defined in \eqref{Ckl:side:splus}, $\alpha=\frac{a_1}{a_2}$, $(x_i,y_j)\in \Omega_h, (x_i- wh,y_j)=(x_i^*,y_j^*)\in \Gamma_1$, and  $0<w<1$. Since $r_1, r_2, r_3, q_1, q_2, q_3$ are all positive and  $p_1, p_2, p_3$ are all negative for $0<w<1$, we have that the coefficients $C_1$ and $C_2$ of $h^{-1} $ and $h$ in $(\Delta_h \Theta)_{i,j}$ in \eqref{interface:Deltah:Theta:ij} are both positive for $0<w<1$  if and only if
\be\label{proof:rho:less}
\rho< \min \Big\{ \frac{        w\alpha+(1-w)}{ 2(w+1)\alpha+2(1-w)}, -\frac{q_1\alpha^2+q_2\alpha+q_3}{p_1\alpha^2+p_2\alpha+p_3} \Big \}.
\ee
For any choice of $\rho$ satisfying $\rho\le 0$, it is straightforward to observe that \eqref{proof:rho:less} holds for all $w\in (0,1)$ and for all $\alpha \in (0,\infty)$.

Similar as \eqref{Lh:Theta:ij:order6} and \eqref{order6:proof:Theta}, by $\{C_{k,\ell}\}_{k,\ell=-1,0,1}$  in \cref{Intersect:thm:cross1:w1w2:order2}, we have
\be\label{intersection:Deltah:Theta:ij}
(\Delta_h \Theta)_{i,j}=	-h^{-1}(\mathcal{L}_h \Theta)_{i,j} = h^{-1}  C_3+hC_4,
\ee
where
		\be\label{Thm34:proof:C3C4}
		\begin{split}
C_3=\frac{-a_{1}a_{3}e_{1}-2a_{1}a_{2}r_{2}-a_{2}a_{3}r_{3}}{a_{1}a_{2}s_{4}+a_{1}a_{3}s_{3}+a_{2}^2s_{2}+a_{2}a_{3}s_{1}}, \qquad C_4=\frac{1}{12}\frac{-(a_{2}r_{2}+a_{3}r_{1})(a_{1}r_{4}+a_{2}r_{3})
}{a_{1}a_{2}s_{4}+a_{1}a_{3}s_{3}+a_{2}^2s_{2}+a_{2}a_{3}s_{1}},
		\end{split}
		\ee
and $s_1, s_2, s_3, s_4, r_1, r_2, r_3, r_4$ are defined in \eqref{Ckl:cross:tplus},  $e_{1}=2w_{1}^2+4w_{2}^2-w_{1}-2w_{2}+3$,  $(x_i,y_j)\in \Omega_h$, and $(x_i- w_1h,y_j- w_2h)=(\xi,\zeta), 0<w_1,w_2<1$.  Since $r_1, r_2, r_3, r_4, e_{1}$ are all positive and $s_1, s_2, s_3, s_4$ are all negative for $0<w_1,w_2<1$, we have  that the coefficients $C_3$ and $C_4$ of $h^{-1} $ and $h$ in $(\Delta_h \Theta)_{i,j}$ in \eqref{intersection:Deltah:Theta:ij} are both positive for $0<w_1,w_2<1$.
Now, we can replace \eqref{order6:proof:Theta} by
		\[
		(\Delta_h \Theta)_{i,j}=
		\begin{cases}
			-h^{-2}(\mathcal{L}_h \Theta)_{i,j}=1, &\text{if } (x_i,y_j)\in \Omega_h \setminus \Omega_{h,\Gamma},\\
			-h^{-1}(\mathcal{L}_h \Theta)_{i,j}= h^{-1}  C_1+hC_2, &\text{if } (x_i,y_j)\in \Omega_h \mbox{ and } (x_i,y_j)\in \Omega_{h,\Gamma_1,\Gamma_2} \cup \Omega_{h,\Gamma_3,\Gamma_4},\\
			-h^{-1}(\mathcal{L}_h \Theta)_{i,j}= h^{-1} C_3+hC_4, &\text{if } (x_i,y_j)\in \Omega_h \mbox{ and }  (x_i,y_j)\in \Omega_{h,\xi,\zeta},				
		\end{cases}
		\]
	where $C_1,C_2,C_3,C_4>0$. Choose $C_5:=\max\{C,\frac{C}{C_1},\frac{C}{C_3}\}$, where $C$ is the positive constant in \eqref{Trun:General}. 	Then we can
	replace
	\eqref{4CH6}  by:
				\[
C_5 h^3 (\Delta_h \Theta)_{i,j}
	= \begin{cases}
	C_5h^3\ge Ch^3, &\text{if } (x_i,y_j)\in \Omega_h \setminus \Omega_{h,\Gamma},\\
C_5 h^3( h^{-1}  C_1+hC_2) \ge 	C_5C_1h^2\ge Ch^2, &\text{if }  (x_i,y_j)\in \Omega_{h,\Gamma_1,\Gamma_2} \cup \Omega_{h,\Gamma_3,\Gamma_4},\\
	C_5 h^3(h^{-1}  C_3+hC_4) \ge C_5C_3h^2\ge Ch^2, &\text{if }  (x_i,y_j)\in \Omega_{h,\xi,\zeta}.
	\end{cases}
	\]
Note that $E_h:=U_h-u_h$.
	So we can replace
	\eqref{order6:proof:compare}  by:
	{\footnotesize{
			\[
			(\Delta_h (E_h+C_5 h^3 \Theta))_{i,j}
			=(\Delta_h E_h)_{i,j}+C_5 h^3 (\Delta_h \Theta)_{i,j}
			\ge \begin{cases}
				R(x_i,y_j)+Ch^3\ge 0, &\text{if } (x_i,y_j)\in \Omega_h \setminus \Omega_{h,\Gamma},\\
				R(x_i,y_j)+Ch^2\ge 0, &\text{if }  (x_i,y_j)\in \Omega_{h,\Gamma_1,\Gamma_2} \cup \Omega_{h,\Gamma_3,\Gamma_4},\\
				R(x_i,y_j)+Ch^2\ge 0, &\text{if }  (x_i,y_j)\in \Omega_{h,\xi,\zeta}.
			\end{cases}
			\]
	}}
Finally, replace \eqref{use:max:order6:inequ} by:
	\begin{align}
		\max_{(x_i,y_j)\in \Omega_h} (E_h)_{i,j}
		&\le
		\max_{(x_i,y_j)\in \Omega_h} (E_h+C_5 h^3 \Theta)_{i,j} \le
		\max_{(x_i,y_j)\in \partial \Omega_h} (E_h+C_5 h^3 \Theta)_{i,j}\le  \frac{181}{12} C_5  h^3.
	\end{align}
	A similar argument can be applied to $-E_h$. Hence, $\|E_h\|_\infty \le  \frac{181}{12} C_5  h^3$. 	Thus \eqref{intersect:order6:formula:order3}  is proved. Finally,  \eqref{intersect:sign:condition}, \eqref{intersect:sum:condition}, and the Dirichlet boundary condition of \eqref{intersect:3} together imply the M-matrix property.
	This completes the proof of \cref{intersect:thm:convergence:order3}.
\end{proof}

\end{document}